\documentclass[11pt]{amsart}
\usepackage{amsmath, amsthm, amssymb, amscd, amsfonts}
\usepackage[all]{xy}
\usepackage[OT2,OT1]{fontenc}
\usepackage{newlfont}
\usepackage[dvips]{graphics}
\usepackage[latin1]{inputenc}
\usepackage{eucal}
\usepackage{latexsym}
\usepackage[dvips, dvipsnames, usenames]{color}
\input xy
\xyoption{all}
\usepackage{ifthen}

\newcommand{\VGamma}{\widehat{\Gamma}}

\def\de{\Delta}

\def\bdem{\begin{proof}}
\def\edem{\end{proof}}
\def\N{\mathbb{N}}

\def\G{\mathbb{G}}

\def\cA{\mathcal{A}}

\def\cK{\mathcal{K}}

\def\cP{\mathcal{P}}

\def\Z{\mathbb{Z}}

\def\bH{\mathfrak{H}}

\def\bB{\mathfrak{B}}

\def\bI{\mathfrak{I}}

\def\bS{\mathfrak{S}}
\def\bP{\mathfrak{P}}

\def\bG{\mathfrak{G}}

\def\bX{\mathfrak{X}}

\def\bq{\mathfrak{q}}

\def\t{\otimes}

\def\kk{\mathsf{k}}
\newcommand{\hlb}{{\mathcal H}}

\def\yd{{}^{H}_{H}\mathcal{YD}}
\def\gyd{{}^{\kk \Gamma}_{\kk \Gamma}\mathcal{YD}}

\def\xx{\mathbb{X}}

\def\gr{\emph{gr}}
\def\td{\triangleright}

\def\unon{\left\{ 1,\ldots ,\theta \right\}}
\def\zt{\Z^{\theta}}

\def\dc{\de(\chi)}

\newcommand\id{{\operatorname{id}}}
\newcommand\ad{\operatorname{ad}}

\newcommand\ord{{\operatorname{ord}}}

\newcommand\Aut{{\operatorname{Aut}}}
\newcommand{\e}{\mathbf e}
\newcommand{\f}{\mathbf f}
\newcommand{\ub}{\mathbf u}
\newcommand{\vb}{\mathbf v}
\newcommand{\wb}{\mathbf w}
\newcommand{\zb}{\mathbf z}
\newcommand{\ot}{{\otimes}}
\newcommand{\otv}{{1\le i \le \theta}}
\newcommand{\otvz}{{1\le i, j \le \theta}}

\newcommand{\qf}{\widetilde{q}}

\newcommand{\mf}{\widetilde{m}}
\newcommand{\si}{s_{i,F}}
\newcommand{\sE}{s_{i,E}}
\newcommand{\wo}{W_0(\chi)}

\newcommand{\cI }{\mathcal{I}}

\newcommand{\Dchaintwo}[4]{
\rule[-3\unitlength]{0pt}{8\unitlength}
\begin{picture}(14,5)(0,3)
\put(1,2){\ifthenelse{\equal{#1}{l}}{\circle*{2}}{\circle{2}}}
\put(2,2){\line(1,0){10}}
\put(13,2){\ifthenelse{\equal{#1}{r}}{\circle*{2}}{\circle{2}}}
\put(1,5){\makebox[0pt]{\scriptsize #2}}
\put(7,4){\makebox[0pt]{\scriptsize #3}}
\put(13,5){\makebox[0pt]{\scriptsize #4}}
\end{picture}}
\newcommand{\Dchainthree}[6]{
\rule[-3\unitlength]{0pt}{8\unitlength}
\begin{picture}(26,5)(0,3)
\put(1,2){\ifthenelse{\equal{#1}{l}}{\circle*{2}}{\circle{2}}}
\put(2,2){\line(1,0){10}}
\put(13,2){\ifthenelse{\equal{#1}{m}}{\circle*{2}}{\circle{2}}}
\put(14,2){\line(1,0){10}}
\put(25,2){\ifthenelse{\equal{#1}{r}}{\circle*{2}}{\circle{2}}}
\put(1,5){\makebox[0pt]{\scriptsize #2}}
\put(7,4){\makebox[0pt]{\scriptsize #3}}
\put(13,5){\makebox[0pt]{\scriptsize #4}}
\put(19,4){\makebox[0pt]{\scriptsize #5}}
\put(25,5){\makebox[0pt]{\scriptsize #6}}
\end{picture}}
\newcommand{\Dtriangle}[7]{
\rule[-3\unitlength]{0pt}{12\unitlength}
\begin{picture}(18,7)(0,3)
\put(4,4){\ifthenelse{\equal{#1}{l}}{\circle*{2}}{\circle{2}}}
\put(5,4){\line(1,0){8}}
\put(14,4){\ifthenelse{\equal{#1}{r}}{\circle*{2}}{\circle{2}}}
\put(4.4472,4.8944){\line(1,2){4.1056}}
\put(9,14){\ifthenelse{\equal{#1}{t}}{\circle*{2}}{\circle{2}}}
\put(13.5528,4.8944){\line(-1,2){4.1056}}
\put(2,3){\makebox[0pt][r]{\scriptsize #2}}
\put(9,17){\makebox[0pt]{\scriptsize #3}}
\put(16,3){\makebox[0pt][l]{\scriptsize #4}}
\put(6,9){\makebox[0pt][r]{\scriptsize #5}}
\put(12.5,9){\makebox[0pt][l]{\scriptsize #6}}
\put(9,1){\makebox[0pt]{\scriptsize #7}}
\end{picture}}
\newcommand{\cDtriangle}[7]{
\rule[-10\unitlength]{0pt}{\unitlength}
\begin{picture}(18,7)(0,10)
\put(4,4){\ifthenelse{\equal{#1}{l}}{\circle*{2}}{\circle{2}}}
\put(5,4){\line(1,0){8}}
\put(14,4){\ifthenelse{\equal{#1}{r}}{\circle*{2}}{\circle{2}}}
\put(4.4472,4.8944){\line(1,2){4.1056}}
\put(9,14){\ifthenelse{\equal{#1}{t}}{\circle*{2}}{\circle{2}}}
\put(13.5528,4.8944){\line(-1,2){4.1056}}
\put(2,3){\makebox[0pt][r]{\scriptsize #2}}
\put(9,17){\makebox[0pt]{\scriptsize #3}}
\put(16,3){\makebox[0pt][l]{\scriptsize #4}}
\put(6,9){\makebox[0pt][r]{\scriptsize #5}}
\put(12.5,9){\makebox[0pt][l]{\scriptsize #6}}
\put(9,1){\makebox[0pt]{\scriptsize #7}}
\end{picture}}

\newcommand{\longDchaintwo}[4]{
\rule[-3\unitlength]{0pt}{8\unitlength}
\begin{picture}(18,5)(0,3)
\put(1,2){\ifthenelse{\equal{#1}{l}}{\circle*{2}}{\circle{2}}}
\put(2,2){\line(1,0){14}}
\put(17,2){\ifthenelse{\equal{#1}{r}}{\circle*{2}}{\circle{2}}}
\put(1,5){\makebox[0pt]{\scriptsize #2}}
\put(9,4){\makebox[0pt]{\scriptsize #3}}
\put(17,5){\makebox[0pt]{\scriptsize #4}}
\end{picture}}

\newlength{\mpb}

\newcommand{\Dchainofur}[8]{
\rule[-3\unitlength]{0pt}{5\unitlength}
\begin{picture}(38,5)(0,3)
\put(1,2){\ifthenelse{\equal{#1}{1}}{\circle*{2}}{\circle{2}}}
\put(2,2){\line(1,0){10}}
\put(13,2){\ifthenelse{\equal{#1}{2}}{\circle*{2}}{\circle{2}}}
\put(14,2){\line(1,0){10}}
\put(25,2){\ifthenelse{\equal{#1}{3}}{\circle*{2}}{\circle{2}}}
\put(26,2){\line(1,0){10}}
\put(37,2){\ifthenelse{\equal{#1}{4}}{\circle*{2}}{\circle{2}}}
\put(1,5){\makebox[0pt]{\scriptsize #2}}
\put(7,4){\makebox[0pt]{\scriptsize #3}}
\put(13,5){\makebox[0pt]{\scriptsize #4}}
\put(19,4){\makebox[0pt]{\scriptsize #5}}
\put(25,5){\makebox[0pt]{\scriptsize #6}}
\put(31,4){\makebox[0pt]{\scriptsize #7}}
\put(37,5){\makebox[0pt]{\scriptsize #8}}
\end{picture}}

\newcommand{\bigDchainofur}[8]{
\rule[-3\unitlength]{0pt}{5\unitlength}
\begin{picture}(42,5)(0,3)
\put(1,2){\ifthenelse{\equal{#1}{1}}{\circle*{2}}{\circle{2}}}
\put(2,2){\line(1,0){14}}
\put(17,2){\ifthenelse{\equal{#1}{2}}{\circle*{2}}{\circle{2}}}
\put(18,2){\line(1,0){10}}
\put(29,2){\ifthenelse{\equal{#1}{3}}{\circle*{2}}{\circle{2}}}
\put(30,2){\line(1,0){10}}
\put(41,2){\ifthenelse{\equal{#1}{4}}{\circle*{2}}{\circle{2}}}
\put(1,5){\makebox[0pt]{\scriptsize #2}}
\put(9,4){\makebox[0pt]{\scriptsize #3}}
\put(17,5){\makebox[0pt]{\scriptsize #4}}
\put(23,4){\makebox[0pt]{\scriptsize #5}}
\put(29,5){\makebox[0pt]{\scriptsize #6}}
\put(35,4){\makebox[0pt]{\scriptsize #7}}
\put(41,5){\makebox[0pt]{\scriptsize #8}}
\end{picture}}

\newcommand{\Dthreefork}[8]{
\rule[-9\unitlength]{0pt}{12\unitlength}
\begin{picture}(28,12)(0,9)
\put(2,10){\ifthenelse{\equal{#1}{l}}{\circle*{2}}{\circle{2}}}
\put(3,10){\line(1,0){10}}
\put(14,10){\ifthenelse{\equal{#1}{m}}{\circle*{2}}{\circle{2}}}
\put(15,10){\line(1,1){7}}
\put(15,10){\line(1,-1){7}}
\put(22,18){\ifthenelse{\equal{#1}{t}}{\circle*{2}}{\circle{2}}}
\put(22,2){\ifthenelse{\equal{#1}{b}}{\circle*{2}}{\circle{2}}}
\put(2,12){\makebox[0pt]{\scriptsize #2}}
\put(8,11){\makebox[0pt]{\scriptsize #3}}
\put(14,12){\makebox[0pt]{\scriptsize #4}}
\put(19,16){\makebox[0pt][r]{\scriptsize #5}}
\put(19,4){\makebox[0pt][r]{\scriptsize #6}}
\put(24,17){\makebox[0pt][l]{\scriptsize #7}}
\put(24,2){\makebox[0pt][l]{\scriptsize #8}}
\end{picture}}

\newcommand{\Drightofway}[9]{
\rule[-9\unitlength]{0pt}{12\unitlength}
\begin{picture}(28,12)(0,9)
\put(2,10){\ifthenelse{\equal{#1}{l}}{\circle*{2}}{\circle{2}}}
\put(3,10){\line(1,0){10}}
\put(14,10){\ifthenelse{\equal{#1}{m}}{\circle*{2}}{\circle{2}}}
\put(15,10){\line(1,1){7}}
\put(15,10){\line(1,-1){7}}
\put(22,18){\ifthenelse{\equal{#1}{t}}{\circle*{2}}{\circle{2}}}
\put(22,2){\ifthenelse{\equal{#1}{b}}{\circle*{2}}{\circle{2}}}
\put(22,3){\line(0,1){14}}
\put(2,12){\makebox[0pt]{\scriptsize #2}}
\put(8,11){\makebox[0pt]{\scriptsize #3}}
\put(14,12){\makebox[0pt]{\scriptsize #4}}
\put(19,16){\makebox[0pt][r]{\scriptsize #5}}
\put(19,4){\makebox[0pt][r]{\scriptsize #6}}
\put(24,18){\makebox[0pt][l]{\scriptsize #7}}
\put(23,10){\makebox[0pt][l]{\scriptsize #8}}
\put(24,1){\makebox[0pt][l]{\scriptsize #9}}
\end{picture}}

\numberwithin{equation}{section} \theoremstyle{plain}
\newtheorem{theorem}{Theorem}[section]
\newtheorem{cor}[theorem]{Corollary}
\newtheorem{lema}[theorem]{Lemma}

\newtheorem{prop}[theorem]{Proposition}

\theoremstyle{definition}
\newtheorem{definition}[theorem]{Definition}

\newtheorem{obs}[theorem]{Remark}

\def\pf{\begin{proof}}
\def\epf{\end{proof}}

\begin{document}

\title[On Nichols algebras with standard braiding]{On Nichols algebras with standard braiding}
\author{Iv\'an Ezequiel Angiono}
\address{Facultad of Matem\'atica, Astronom\'\i a y F\'\i sica
\newline \indent
Universidad Nacional of C\'ordoba
\newline
\indent CIEM -- CONICET
\newline
\indent (5000) Ciudad Universitaria, C\'ordoba, Argentina}
 \email{angiono@mate.uncor.edu}
\date{\today}
\thanks{ 2000 {\it Mathematics Subject Classification:} Primary 17B37;
Secondary: 16W20,16W30 \newline \indent {\it Key words and
phrases:} quantized enveloping algebras, Nichols algebras,
automorphisms of non-commutative algebras. }

\begin{abstract}
The class of standard braided vector spaces, introduced by
Andruskiewitsch and the author in \texttt{arXiv:math/0703924v2} to
understand the proof of a theorem of Heckenberger \cite{H2}, is
slightly more general than the class of braided vector spaces of
Cartan type. In the present paper, we classify standard braided
vector spaces with finite-dimensional Nichols algebra. For any
such braided vector space, we give a PBW-basis, a closed formula
of the dimension and a presentation by generators and relations of
the associated Nichols algebra.
\end{abstract}

\setlength{\unitlength}{1mm} \settowidth{\mpb}{$q_0\in k^\ast
\setminus \{-1,1\}$,}

\maketitle

\setcounter{tocdepth}{2} \tableofcontents

\section*{Introduction}

A breakthrough in the development of the theory of Hopf algebras
was the discovery of quantized enveloping algebra by Drinfeld and
Jimbo \cite{Dr,Ji}. This special class of Hopf algebras was
intensively studied by many authors and from many points of view.
In particular, finite-dimensional analogues of quantized
enveloping algebras were introduced and investigated by Lusztig
\cite{L1,L2}.

About ten years ago, a classification program of pointed Hopf
algebras was launched by Andruskiewitsch and Schneider \cite{AS1},
see also \cite{AS5}. The success of this program depends on
finding solutions to several questions, among them:

\begin{quote}
\cite[Question 5.9]{A} Given a braided vector space of diagonal
type $V$, such that the entries of its matrix are roots of 1,
compute the dimension of the associated Nichols algebra $\bB(V)$.
If it is finite, give a nice presentation of $\bB(V)$.
\end{quote}

Partial answers to this question were given in \cite{AS2, H2} for
the class of braided vector spaces of Cartan type. These answers
were already crucial to prove a classification theorem for
finite-dimensional Hopf algebras whose group is abelian with prime
divisors of the order great than 7 \cite{AS6}. Later, a complete
answer to the first part of \cite[Question 5.9]{A} was given in
\cite{H3}.

The notion of standard braided vector space, a special kind of
diagonal braided vector space, was introduced in \cite{AA}, see
Definition \ref{standardbraiding} below. This class includes
properly the class of braided vector spaces of Cartan type.

The purpose of this paper is to develop from scratch the theory of
standard braided vector spaces. Here are our main contributions:

\begin{itemize}
    \item We give a complete classification of standard braided vector
    spaces with finite-dimensional Nichols algebras. As usual, we
    may assume the connectedness of the corresponding braiding.
    It turns out that standard braided vector spaces are of Cartan type
    when the associated Cartan matrix is of type $C$, $ D$, $E$ or $F$, see Proposition \ref{solocartan}.
    For types $A$, $B$, $G$ there are standard braided vector spaces not of
    Cartan type; these are listed in Propositions \ref{casosAn},
    \ref{casosBn} and \ref{casosG2}. Those of type $A_2$ and $B_2$
    appeared already in \cite{Gr}. Our classification does not rely on
    \cite{H3}, but we can identify our examples in the tables of
    \cite{H3}.

\medbreak
    \item We describe a concrete PBW-basis of the Nichols algebra
    of a standard braided vector space as in the previous point; this follows from the
    general theory of Kharchenko \cite{Kh} together with Theorem 1 of \cite{H2}. As an application, we
    give closed formulas for the dimension of these
    Nichols algebras.

\medbreak    \item We present a concrete set of defining relations
of the Nichols algebras of standard braided vector spaces as in
the previous points. This is an answer to the second part of
\cite[Question 5.9]{A} in the standard case. We note
    that this seems to be new even for Cartan type, for some values of
    the roots of 1 appearing in the picture. Essentially, these relations are
    either quantum Serre  relations or powers of root vectors; but
    in some cases, there are some substitutes of the quantum Serre  relations
    due to the smallness of the intervening root vectors. Some of these substitutes
    can be recognized already in the relations in \cite{AD}.
\end{itemize}

\medbreak Here is the plan of this article. In section
\ref{section:pbw}, we collect different tools that will be used in
the following sections. Namely, we recall the definition of Lyndon
words and give some properties about them, such as the Shirshov
decomposition, in subsection \ref{subsection:lyndon}. In
\ref{subsection:bvs}, we discuss the notions of hyperletter and
hyperword following \cite{Kh} (they are called superletter and
superword in loc. cit.); these are certain specific iterations of
braided commutators applied to Lyndon words. Next, in subsection
\ref{subsection:pbw}, a PBW basis is given for any quotient of the
tensor algebra of a diagonal braided vector space $V$ by a Hopf
ideal using these hyperwords. This applies in particular to
Nichols algebras.

In section \ref{section:transf}, after some technical
preparations, we present a transformation of a braided graded Hopf
algebra into another, with different space of degree one. This
generalizes an analogous transformation for Nichols algebras given
in \cite[Prop. 1]{H2} -- see Subsection
\ref{subsection:transformation}.

In section \ref{section:standard} we classify standard braided
vector spaces with finite dimensional Nichols algebra. In
subsection \ref{subsection:weylgroupoid}, we prove that if the set
of PBW generators is finite, then the associated generalized
Cartan matrix is of finite type. So in subsection
\ref{subsection:classification} we obtain all the standard
braidings associated to Nichols algebras of finite dimension.

Section \ref{section:nichols} is devoted to PBW-bases of Nichols
algebras of standard braided vector spaces with finite Cartan
matrix. In subsection \ref{subsection:PBWbases} we prove that
there is exactly one PBW generator whose degree corresponds with
each positive root associated to the finite Cartan matrix. We give
a set of PBW-generators in subsection \ref{subsection:generators},
following a nice presentation from \cite{LR}. As a consequence, we
compute the dimension in Subsection \ref{subsection:dimension}.

The main result of this paper is the explicit presentation by
generators and relations of Nichols algebras of standard braided
vector spaces with finite Cartan matrix, given in section
\ref{section:genrel}. This result relies on the explicit PBW-basis
and the transformation described in Subsection
\ref{subsection:transformation}. In subsection
\ref{subsection:relations}, we state some relations for Nichols
algebras of standard braidings, and prove some facts about the
coproduct. Subsections \ref{subsection:An}, \ref{subsection:Bn}
and \ref{subsection:presentationG2} contain the explicit
presentation for types $A_{\theta}$, $B_{\theta}$ and $G_2$,
respectively. For this, we establish relations among the elements
of the PBW-basis, inspired in \cite{AD} and \cite{Gr}. We finally
prove the presentation in the case of Cartan type in
\ref{subsection:presentation}. To our knowledge, this is the first
self-contained exposition of Nichols algebras of braided vector
spaces of Cartan type.

\medbreak \textbf{Notation.} We fix an algebraically closed field
$\kk$ of characteristic 0; all vector spaces, Hopf algebras and
tensor products are considered over $\kk$.

Given $n \in \N$ and $q \in \kk$, $q \notin \cup_{0 \leq j \leq n}
\G_j$, we denote
$$ \binom{n}{j}_q = \frac{(n)q!}{(k)_q! (n-k}_q!, \quad \mbox{where }(n)_q!= \prod_{j=1}^n (k)_q, \quad \mbox{and } (k)_q= \sum_{j=0}^{k-1} q^j. $$
For each $n=(n_1, \ldots , n_{\theta}) \in \zt$, we set
$x^{n}=x_1^{n_1} \cdots x_{\theta}^{n_{\theta}} \in \kk [[x_1^{\pm
1}, \ldots, x_{\theta}^{\pm 1}]]$. Also we denote
\begin{eqnarray*}
\bq _h(\mathrm{t}) := \frac{\mathrm{t}^h-1}{\mathrm{t}-1} \in \kk
[\mathrm{t}], \quad h \in \N; \quad \bq _{\infty}(\mathrm{t}):=
\frac{1}{1-\mathrm{t}}= \sum_{s=0}^{\infty} \mathrm{t}^s \in \kk
[[\mathrm{t}]].
\end{eqnarray*}

For each $N \in \N$, $\G_N$ denotes the set of primitive $N$-th
roots of 1 in $\kk$.

For each $\theta \in \N$ and each $\zt$-graded vector spaces
$\bB$, we denote by $\hlb_{\bB}= \sum_{n \in \ \zt} \dim \bB^n$
the Hilbert series associated to $\bB$.

Let $C= \oplus_{n \in \N_0} C_{i+j}$ be a $\N_0$-graded coalgebra,
with projections $\pi_n: C \rightarrow C_n$. Given $i,j \geq 0$,
we denote by
$$\de_{i,j}:= (\pi_i \otimes \pi_j) \circ \de: C_{i+j}
\rightarrow C_i \otimes C_j, $$ the (i,j)-th component of the
comultiplication.

\section{PBW-basis}\label{section:pbw}

Let $A$ be an algebra, $P,S \subset A$ and $h: S \mapsto \N \cup
\{ \infty \}$. Let also $<$ be a linear order on $S$. Let us
denote by $B(P,S,<,h)$ the set
\begin{align*}
\big\{ &p\,s_1^{e_1}\dots s_t^{e_t}: t \in \N_0, \quad s_1>\dots
>s_t,\quad s_i \in S, \quad 0<e_i<h(s_i), \quad p \in P \big\}.
\end{align*}

If $B(P,S,<,h)$ is a linear basis of $A$, then we say that
$(P,S,<,h)$ is a set of \emph{PBW generators} with height $h$, and
that $B(P,S,<,h)$ is a \emph{PBW-basis} of $A$. Occasionally, we
shall simply say that $S$ is a PBW-basis of $A$.

\bigbreak In this Section, we describe-- following \cite{Kh}-- an
appropriate PBW-basis of a braided graded Hopf algebra $\bB =
\oplus_{n\in \N} \bB^n$ such that $\bB^1 \cong V$, where $V$ is a
braided vector space of diagonal type. This applies in particular,
to the Nichols algebra $\bB(V)$. In Subsection
\ref{subsection:lyndon} we recall the classical construction of
Lyndon words. Let $V$ be a vector space $V$ together with a fixed
basis. Then there is a basis of the tensor algebra $T(V)$ by
certain words satisfying a special condition, called Lyndon words.
Each Lyndon word has a canonical decomposition as a product of a
pair of smaller Lyndon words, called the Shirshov decomposition.

We briefly remind the notions of braided vector space $(V,c)$ of
diagonal type and Nichols algebra in Subsection
\ref{subsection:bvs}. Then we recall-- in Subsection
\ref{subsection:pbw}-- the definition of the hyperletter $[l]_c$,
for any Lyndon word $l$; this is the braided commutator of the
hyperletters corresponding to the words in the Shirshov
decomposition.  The hyperletters are a set of generators for a
PBW-basis of $T(V)$ and their classes form a PBW-basis of $\bB$.

\subsection{Lyndon words}\label{subsection:lyndon}

\

Let $\theta\in \N$.  Let $X$ be a set with $\theta$ elements  and
fix a numeration $x_1,\dots, x_{\theta}$ of $X$; this induces a
total order on $X$. Let $\xx$ be the corresponding vocabulary (the
set of words with letters in $X$) and consider the lexicographical
order on $\xx$.

\begin{definition}
An element $u \in \xx$, $u\neq 1$, is called a \emph{Lyndon word}
if $u$ is smaller than any of its proper ends; that is, if $u=vw$,
$v,w \in\xx - \left\{ 1 \right\}$, then $u<w$. The set of Lyndon
words is denoted by $L$.
\end{definition}

We shall need the following properties of Lyndon words.

\begin{enumerate}
    \item Let $u \in \xx-X$. Then $u$ is Lyndon if and only if for any
    representation $u=u_1 u_2$, with $u_1,u_2 \in \xx$ not empty, one has $u_1u_2=u < u_2u_1$.
    \item Any Lyndon word begins by its smallest letter.
    \item If $u_1,u_2 \in L, u_1<u_2$, then $u_1u_2 \in L$.
\end{enumerate}

\bigbreak The basic Theorem about Lyndon words, due to Lyndon,
says that any word $u \in \xx$ has a unique decomposition
\begin{equation}\label{descly}
u=l_1l_2\dots  l_r,
\end{equation}
with $l_i \in L$, $l_r \leq \dots \leq l_1$, as a product of non
increasing Lyndon words. This is called the \emph{Lyndon
decomposition} of $u \in \xx$; the $l_i \in L$ appearing in the
decomposition \eqref{descly} are called the \emph{Lyndon letters}
of $u$.

\bigbreak The lexicographical order of $\xx$ turns out to be the
same as the lexicographical order in the Lyndon letters. Namely,
if $v=l_1\dots l_r$ is the Lyndon decomposition of $v$, then $u<v$
if and only if:
\begin{enumerate}
    \item[(i)] the Lyndon decomposition of $u$ is $u=l_1\dots l_i$, for some $1\leq i <r$,
    or

    \item[(ii)] the Lyndon decomposition of $u$ is $u=l_1\dots l_{i-1}ll_{i+1}'\dots
    l_s'$, for some $1 \leq i <r$, $s \in \N$ and $l, l_{i+1}',\dots ,l_s'$ in $L$, with
    $l<l_i$.
\end{enumerate}

Here is another useful characterization of Lyndon words.

\begin{lema}
Let $u \in\xx-X$. Then $u \in L$ if and only if there exist
$u_1,u_2 \in L$ with $u_1<u_2$ such that $u=u_1u_2$.
\end{lema}

\pf See \cite[p.6, Shirshov Th.]{Kh}. \epf

\begin{definition}
Let $u \in L-X$. A decomposition $u=u_1u_2$, with $u_1,u_2 \in L$
such that $u_2$ is the smallest end among those proper non-empty
ends of $u$ is called the \emph{ Shirshov decomposition } of $u$.
\end{definition}

Let $u,v,w \in L$ be such that $u=vw$. Then $u=vw$ is the Shirshov
decomposition of $u$ if and only if either $v \in X$, or else if
$v=v_1v_2$ is the Shirshov decomposition of $v$, then $w \leq
v_2$.

\subsection{Braided vector spaces of diagonal
type and Nichols algebras}\label{subsection:bvs}

\

A braided vector space is a pair $(V,c)$, where $V$ is a vector
space and $c\in \Aut (V\ot V)$ is a solution of the braid
equation: $$(c\otimes \id) (\id\otimes c) (c\otimes \id) =
(\id\otimes c) (c\otimes \id) (\id\otimes c).$$ We extend the
braiding to $c:T(V)\ot T(V) \to T(V)\ot T(V)$ in the usual way. If
$x,y\in T(V)$, then the braided commutator is
\begin{equation}\label{eqn:braidedcommutator}
[x,y]_c := \text{multiplication } \circ \left( \id - c \right)
\left( x \ot y \right).
\end{equation}

Assume that $\dim V<\infty$ and pick a basis $X = \{x_1,\dots
,x_{\theta}\}$ of $V$; we may then identify $\kk \xx$ with $T(V)$.
We consider the following gradings of the algebra $T(V)$:

\begin{enumerate}
    \item[(i)] The usual $\N_0$-grading $T(V) = \oplus_{n\geq 0}T^n(V)$.
If $\ell$ denotes the length of a word in $\xx$, then $T^n(V)=
\oplus_{x\in\xx, \, \ell(x) = n}\kk x$.

    \item[(ii)] Let $\e_1, \dots, \e_\theta$ be the canonical basis of $\zt$.
    Then $T(V)$ is also $\zt$-graded, where the degree is
    determined by $\deg x_i = \e_i$, $\otv$.\label{paginatres}
\end{enumerate}

A braided vector space $(V,c)$ is of \emph{diagonal type} with
respect to the basis $x_1, \dots x_\theta$ if there exist
$q_{ij}\in \kk^{\times}$ such that $c(x_i \ot x_j)= q _{ij} x_j
\ot x_i$, $\otvz$. Let $\chi: \zt\times \zt \to \kk^{\times}$ be
the bilinear form determined by $\chi(\e_i, \e_j) = q_{ij}$,
$\otvz$. Then
\begin{equation}\label{braiding}
    c(u \ot v)= \chi( \deg u, \deg v ) v \ot u
\end{equation}
for any $u,v \in \xx$, where $q_{u,v} = \chi(\deg u, \deg v)\in
\kk^{\times}$. In this case, the braided commutator satisfies a
``braided" Jacobi identity as well as braided derivation
properties, namely
\begin{align}\label{idjac}
\left[\left[ u, v \right]_c, w \right]_c &= \left[u, \left[ v, w
\right]_c \right]_c
 - \chi( \alpha, \beta ) v \ \left[ u, w \right]_c + \chi( \beta, \gamma) \left[ u,
w \right]_c \ v,
 \\
\label{der}
    \left[ u,v \ w \right]_c &= \left[ u,v \right]_c w + \chi( \alpha, \beta ) v \ \left[ u,w \right]_c,
\\ \label{der2} \left[ u \ v, w \right]_c &= \chi( \beta, \gamma ) \left[ u,w \right]_c \ v + u \ \left[ v,w \right]_c,
\end{align}
for any homogeneous $u,v,w \in T(V)$, of degrees $\alpha, \beta, \gamma \in \N^{\theta}$, respectively.

\bigbreak We denote by $\yd$ the category of Yetter-Drinfeld
module over $H$, where $H$ is a Hopf algebra with bijective
antipode. Any $V\in \yd$ becomes a braided vector space \cite{Mo}.
If $H$ is the group algebra of a finite abelian group, then  any
$V\in \yd$ is a braided vector space of diagonal type. Indeed, $V
= \oplus_{g\in \Gamma, \chi\in \VGamma}V_{g}^{\chi}$, where
$V_{g}^{\chi} = V^{\chi} \cap V_{g}$, $V_{g} = \{v\in V \mid
\delta(v) = g\otimes v\}$, $V^{\chi} = \{v\in V \mid  g \cdot v =
\chi(g)v \text{ for all } g \in \Gamma\}$.  The braiding is given
by $ c(x\otimes y) = \chi(g) y\otimes x$,   for all $x\in V_{g}$,
$g \in \Gamma$, $y\in V^{\chi}$, $\chi \in \VGamma$.

Reciprocally, any braided vector space of diagonal type can be
realized as a Yetter-Drinfeld module over the group algebra of an
abelian group.

\bigbreak If $V\in \yd$, then the tensor algebra $T(V)$ admits a
unique structure of graded braided Hopf algebra in $\yd$ such that
$V \subseteq \cP(V)$. Following \cite{AS5}, we consider the class
$\bS$ of all the homogeneous two-sided ideals $I \subseteq T(V)$
such that

\begin{itemize}
    \item $I$ is generated by homogeneous elements of degree $\geq
    2$,
    \item $I$ is a Yetter-Drinfeld submodule of $T(V)$,
    \item $I$ is a Hopf ideal: $\Delta(I) \subset I\ot T(V) +
    T(V)\ot I$.
\end{itemize}

The Nichols algebra $\bB(V)$ associated to $V$ is the quotient of
$T(V)$ by the maximal element $I(V)$ of $\bS$.

\bigbreak Let $(V,c)$ be a braided vector space of diagonal type,
and assume that $q_{ij}=q_{ji}$ for all $i,j$. Let $\Gamma$ be the
free abelian group of rank $\theta$, with basis $g_1, \ldots,
g_{\theta}$, and define the characters $\chi_1, \ldots,
\chi_{\theta}$ of $\Gamma$ by
    \[ \chi_j(g_i)=q_{ij}, \quad 1 \leq i,j \leq \theta. \]
Consider $V$ as a Yetter-Drinfeld module over $\kk \Gamma$ by
defining $x_i \in V_{g_i}^{\chi_i}$.

We shall need the following proposition.

\begin{prop}\label{formabilineal} \cite[Prop. 1.2.3]{L3}, \cite[Prop. 2.10]{AS5}.
Let $a_1, \ldots,a_{\theta} \in \kk^{\times}$. There is a unique
bilinear form $( | ): T(V) \times T(V) \rightarrow \kk$ such that
$(1|1)=1$, and:
\begin{eqnarray}
    (x_i | x_j) &=& \delta_{ij}a_i, \quad \mbox{for all  } i,j; \label{bilinearprop1}
    \\ (x|yy') &=& (x_{(1)} | y) (x_{(2)} | y'), \quad \mbox{for all  } x,y,y' \in
    T(V);\label{bilinearprop2}
    \\ (xx'|y) &=& (x|y_{(1)}) (x'|y_{(2)}), \quad \mbox{for all  } x,x',y \in T(V). \label{bilinearprop3}
\end{eqnarray}
This form is symmetric and also satisfies
\begin{equation}
    (x|y)=0, \quad \mbox{for all  } x \in T(V)_g, \ y \in T(V)_h, \ g,h \in \Gamma, \ g \neq h. \label{bilinearprop4}
\end{equation}

The quotient $T(V)/I(V)$, where $$I(V):= \left\{ x \in T(V):
(x|y)=0, \ \forall y \in T(V) \right\}$$ is the radical of the
form, is canonically isomorphic to the Nichols algebra of $V$.
Thus, $(|)$ induces a non degenerate bilinear form on $\bB(V)$
denoted by the same name. \qed
\end{prop}

If $(V,c)$ is of diagonal type, then the ideal $I(V)$ is
$\zt$-homogeneous hence $\bB(V)$ is $\zt$-graded. See \cite{AS4}
for details. The following statement, that we include for later
reference, is well-known.

\begin{lema}\label{conditions}
Let $V$ a braided vector space of diagonal type, and consider its
Nichols algebra $\bB(V)$.
\begin{enumerate}
    \item[(a)] If $q_{ii}$ is a root of unit of order $N>1$, then $x_i^N=0$.
    \smallbreak

    \item[(b)] If $i \neq j$, then $(ad_c x_i)^{r}(x_j)=0$ if and only if
    \begin{center}
    $(r)!_{q_{ii}} \prod _{0 \leq k \leq r-1} (1-q_{ii}^k q_{ij}q_{ji})=0$.
    \end{center}
    \smallbreak

    \item[(c)] If $i\neq j$ and $q_{ij}q_{ji}=q_{ii}^r$, for some
    $r$ such that
    $0\leq -r < \ord(q_{ii})$, then $(ad_c x_i)^{1-r}(x_j)=0$.
    \qed
\end{enumerate}
\end{lema}

\subsection{PBW basis of a quotient of the tensor algebra by a Hopf
ideal}\label{subsection:pbw}

\

Let $(V,c)$ be a braided vector space with a basis $X = \{x_1,
\dots, x_\theta\}$; identify $T(V)$ with $\kk \xx$.  An important graded endomorphism $\left[ - \right]_c$ of $\kk \xx$ is given by
$$
\left[ u \right]_c := \begin{cases} u,& \text{if } u = 1
\text{ or }u \in X;\\
[\left[ v \right]_c, \left[ w \right]_c]_c,  & \text{if } u \in
L, \, \ell(u)>1 \text{ and }u=vw \\ &\qquad\text{ is the Shirshov decomposition};\\
\left[ u_1 \right]_c \dots  \left[ u_t \right]_c,& \text{if } u
\in \xx-L \\ &\qquad
\text{ with Lyndon decomposition  }u=u_1\dots u_t;\\
\end{cases}
$$

Let us now assume that $(V,c)$ is of diagonal type with respect to
the basis $x_1, \dots, x_\theta$, with matrix $(q_{ij})$.

\begin{definition} The \emph{hyperletter} corresponding to
$l \in L$ is the element $\left[ l \right]_c$. A \emph{hyperword}
is a word in hyperletters, and a \emph{monotone hyperword} is a
hyperword of the form $W=\left[u_1\right]_c^{k_1}\dots
\left[u_m\right]_c^{k_m}$, where $u_1>\dots >u_m$.
\end{definition}

\begin{obs}\label{corchete}
If $u \in L$, then $\left[ u \right]_c$ is a homogeneous
polynomial with coefficients in $\mathbb{Z} \left[q_{ij}\right]$
and $ \left[ u \right]_c\in u+ \kk \xx^{\ell(u)}_{>u}$.
\end{obs}

\medbreak  The hyperletters inherit the order from the Lyndon
words; this induces in turn an ordering in the hyperwords (the
lexicographical order on the hyperletters). Now, given monotone
hyperwords  $W,V$, it can be shown that
    \[W=\left[w_1\right]_c\dots \left[w_m\right]_c > V=\left[v_1\right]_c\dots  \left[v_t\right]_c, \]
where $w_1 \geq \dots  \geq w_r, v_1 \geq \dots  \geq v_s$, if and
only if     \[w=w_1\dots w_{m} > v=v_i\dots v_t. \] Furthermore,
the principal word  of the polynomial $W$, when  decomposed as sum
of monomials, is $w$ with coefficient 1.

\begin{theorem} \label{corch} (Rosso, see \cite{R2}).
Let $m,n \in L$, with $m<n$. Then the braided commutator
$\left[\left[m\right]_c, \left[n\right]_c \right]_c$ is a
$\mathbb{Z} \left[q_{ij}\right]$-linear combination of monotone
hyperwords $\left[l_1\right]_c \dots  \left[l_r\right]_c, l_i \in
L$, such that
\begin{itemize}
    \item the hyperletters of those hyperwords satisfy $n>l_i \geq mn$,
    \item $\left[mn\right]_c$ appears in the expansion with
non-zero coefficient,
    \item any hyperword appearing in this
decomposition satisfies $$\deg (l_1\dots l_r)= \deg(mn).$$ \qed
\end{itemize}
 \end{theorem}

A crucial result of Rosso describes the behavior of the coproduct
of $T(V)$ in the basis of hyperwords.

\begin{lema}\label{copro} \cite{R2}.
Let $u \in \xx$, and $u= u_1\dots u_r v^m, \ v, u_i \in L, v<u_r
\leq \dots  \leq u_1$ the Lyndon decomposition of $u$. Then
\begin{eqnarray*}
        \Delta \left(\left[ u \right]_c\right) &=& 1 \ot \left[ u \right]_c+ \sum ^{m}_{i=0} \binom{ n }{ i } _{q_{v,v}} \left[u_1\right]_c\dots  \left[u_r\right]_c \left[ v \right]_c ^i \ot \left[ v \right]_c^{n-i}
        \\ && + \sum_{ \substack{ l_1\geq \dots  \geq l_p >l, l_i \in L \\ 0\leq j \leq m } } x_{l_1,\dots ,l_p}^{(j)}
        \ot \left[l_1\right]_c\dots
        \left[l_p\right]_c\left[v\right]_c^j;
\end{eqnarray*}
here each $x_{l_1,\dots ,l_p}^{(j)}$ is $\zt$-homogeneous, and
$$\deg(x_{l_1,\dots ,l_p}^{(j)})+\deg(l_1\dots  l_p v^j)= \deg(u).$$
\qed
\end{lema}

As in \cite{U}, we consider another order in $\xx$; it is implicit
in \cite{Kh}.

\begin{definition}
Let $u,v \in \xx$. We say that $u \succ v$ if and only if either
$\ell(u)<\ell(v)$, or else $\ell(u)=\ell(v)$ and $u>v$
(lexicographical order). This $\succ$ is a total order, called the
\emph{deg-lex order}.
\end{definition}

Note that the empty word 1 is the maximal element for $\succ$.
Also, this order is invariant by right and left multiplication.

\medskip

Let now $I$ be a proper ideal  of $T(V)$, and set $R=T(V)/I$.
Let $\pi: T(V) \rightarrow R$ be the canonical projection. Let us
consider the subset of $\xx$:
    \[G_I:= \left\{ u \in \xx: u \notin \\ \kk \xx_{\succ u}+I  \right\}. \]
Notice that
\begin{enumerate}
    \item[(a)] If $u \in G_I$ and $u=vw$, then $v,w \in G_I$.
    \item[(b)] Any word $u \in G_I$  factorizes uniquely as a non-increasing product of Lyndon words in $G_I$.
\end{enumerate}

\begin{prop}\label{firstPBWbasis} \cite{Kh}, see also \cite{R2}.
The set $\pi(G_I)$ is a basis of $R$. \qed
\end{prop}

\noindent In what follows, $I$ is a Hopf ideal. We seek to find a PBW-basis by hyperwords of the quotient $R$ of $T(V)$. For this, we look at the set
\begin{equation}\label{setsi}
S_I:= G_I \cap L.
\end{equation}
We then define the function $h_I: S_I \rightarrow \left\{2,3,\dots
\right\}\cup \left\{ \infty \right\}$ by
\begin{equation}\label{defheight}
    h_I(u):= \min \left\{ t \in \N : u^t  \in \kk \xx_{\succ u^t} + I \right\}.
\end{equation}

The next result plays a fundamental role in this paper.

\begin{theorem}\label{basePBW} \cite{Kh}.
Keep the notation above. Then $$B_I':= B\left( \left\{1+I\right\}
, \left[ S_I \right]_c+I, <, h_I \right)$$
    is a PBW-basis of $H=T(V)/I$. \qed
\end{theorem}
See \cite{Kh} for proofs of the next consequences of the Theorem
\ref{basePBW}.

\begin{cor}\label{cor:primero}
A word $u$ belongs to $G_I$ if and only if the corresponding
hyperletter $\left[u\right]_c$  is not a linear combination,
module $I$, of hyperwords $\left[ w \right]_c$, $w \succ u$, where all the hyperwords belong to $B_I$. \qed
\end{cor}

\begin{prop}\label{altf}
In the conditions of the Theorem \ref{basePBW}, if $v \in S_I$ is
such that $h_I(v)< \infty$, then $q_{v,v}$ is a root of unit. In
this case, if $t$ is the order of $q_{v,v}$, then $h_I(v)=t$. \qed
\end{prop}

\begin{cor}\label{cor:segundo}
If $h_I(v):= h < \infty$, then $\left[ v \right]^{h}$ is a linear
combination of hyperwords $\left[ w \right]_c$, $w \succ u^t$.
\qed
\end{cor}

\section{Transformations of braided graded Hopf algebras}\label{section:transf}

In Subsection \ref{subsection:transformation}, we shall introduce
a transformation over certain graded braided Hopf algebras,
generalizing \cite[Prop. 1]{H2}. It is instrumental step in the
proof of Theorem \ref{presentation}, one of the main results of
this article.

\subsection{Preliminaries on braided graded Hopf algebras}\label{subsection:braidedHA}

\

Let $H$ be the group algebra of an abelian group $\Gamma$. Let
$V\in \yd$ with a basis $X = \{x_1, \dots, x_\theta\}$ such that
$x_i \in V^{\chi_i}_{g_i}$, $1\le i\le \theta$. Let $q_{ij} =
\chi_j(g_i)$, so that $c(x_i\ot x_j) = q_{ij} x_j\ot x_i$, $1\le
i,j\le \theta$.

\bigbreak We fix an ideal $I$ in the class $\bS$; \emph{we assume
that $I$ is $\zt$-homogeneous}. Let $\bB:=T(V)/I$: this is a
braided graded Hopf algebra, $\bB^0=\kk 1$ and $\bB^1=V$. By
definition of $I(V)$, there exists a canonical epimorphism of
braided graded Hopf algebras $\pi: \bB \rightarrow \bB(V)$. Let
$\sigma_i: \bB \rightarrow \bB$ be the algebra automorphism given
by the action of $g_i$.

\begin{prop}\label{derivations}(See for example \cite[2.8]{AS5}).
\begin{enumerate}
    \item For each $1 \leq i \leq \theta$,
    there exists a uniquely determined $(id,\sigma_i)$-derivation
    $D_i: \bB\rightarrow \bB$ with $D_i(x_j)=\delta_{i,j}$ for all $j$.

    \item $I=I(V)$ if and only if $\cap _{i=1}^{\theta} \ker D_i= \kk 1$.\qed
\end{enumerate}
\end{prop}

These operators are defined for each $x \in \bB^{k}, k \geq 1$ by the formula
    \[ \Delta_{n-1,1} (x)= \sum ^{\theta}_{i=1} D_i(x) \ot x_i. \]

Analogously, we can define operators $F_i: \bB \rightarrow \bB$ by $F_i(1)=0$,
    \[ \Delta_{1,n-1} (x)= \sum ^{\theta}_{i=1} x_i \ot F_i(x), \qquad x \in \oplus _{k >0}\bB^{k} . \]

Let $\chi$ be as in \ref{formabilineal}. Consider the action $\td$
of $\kk \zt$ on $\bB$ given by
\begin{equation}\label{actiontriangle}
    \e_i\td b = \chi (\ub , \e_i) b, \qquad  b \mbox{ homogeneous of degree }\ub \in \zt.
\end{equation}
Then, such operators $F_i$ satisfy $F_i(x_j)=\delta_{i,j}$ for all
$j$, and
    \[ F_i(b_1b_2) = F_i(b_1)b_2+ (\e_i \td b_1) F_i(b_2), \quad b_1,b_2 \in \bB. \]

Let $z_r^{(ij)}:=(ad_c x_i)^r(x_{j})$, $i,j \in \unon, i \neq j$
and $r \in \N$.

\begin{obs}
The operators $D_i$, $F_i$ satisfy
\begin{eqnarray}
    && D_i^L(x_i^n)=(n)_{q_{ii}}x_i^{n-1}, \label{13}
    \\ && D_i \left( (ad_c x_i)^r(x_{j_1}\ldots x_{j_s}) \right)=0, \ \forall r,s \geq 1, j_k \neq i, \label{14}
    \\ && D_j \left( z_r^{(ij)} \right)= \prod ^{r-1}_{k=0} \left( 1-q_{ii}^k q_{ij}q_{ji} \right) x_i^r, \forall r \geq 0, \label{15}
    \\ && F_i \left( z_m^{(ij)} \right) = (m)_{q_{ii}}(1-q_{ii}^{m-1}q_{ij}q_{ji}) z_{m-1}^{(ij)}, \label{paragenerar}
    \\ && F_j \left( z_m^{(ij)} \right) = 0, \quad m \geq 1. \label{paragenerar2}
\end{eqnarray}
The proof of the first three identities is as in \cite[Lemma
3.7]{AS4}; the proof of the last two is by induction on $m$.
\end{obs}

\medskip

For each pair $1 \leq i,j \leq \theta, i \neq j$, we define
\begin{eqnarray}
     M_{i,j}(\bB) &:=& \left\{ (ad_cx_i)^m(x_j): m \in \N_0 \right\};
     \\ m_{ij} &:=& \min \left\{ m \in \mathbb{N}: (m+1)_{q_{ii}}(1-q_{ii}^mq_{ij}q_{ji})=0 \right\}.
\end{eqnarray}
Then either $q_{ii}^{m_{ij}}q_{ij}q_{ji}=1$, or
$q_{ii}^{m_{ij}+1}=1$, if $q_{ii}^mq_{ij}q_{ji} \neq 1$ for all
$m=0,1, \ldots, m_{ij}$.

\bigbreak If $\bB=\bB(V)$, then we simply denote $M_{i,j} =
M_{i,j}(\bB(V))$. Note that $(\ad_c x_i)^{m_{ij}+1}x_j=0$ and
$(\ad_c x_i)^{m_{ij}}x_j \neq 0$, by Lemma \ref{conditions}, so
$$ \left| M_{i,j} \right|=m_{ij}+1. $$

\bigbreak

By Theorem \ref{basePBW}, the braided graded Hopf algebra $\bB$
has a PBW-basis consisting of homogeneous elements (with respect
to the $\zt$-grading). As in \cite{H2}, we can even assume that
\begin{itemize}
    \item[$\circledast$] The height of a PBW-generator $\left[ u \right],
\deg(u)=d$, is finite if and only if $2 \leq \ord(q_{u,u}) <
\infty$, and in such case, $h_{I(V)}(u)= \ord(q_{u,u})$.
\end{itemize}

This is possible because if the height of $\left[ u \right],
\deg(u)=d$, is finite, then $2 \leq ord(q_{u,u})=m< \infty$, by
Proposition \ref{altf}. And if $2 \leq \ord(q_{u,u})=m < \infty$,
but $h_{I(V)}(u)$ is infinite, we can add $\left[ u \right]^m$ to
the PBW basis: in this case, $h_{I(V)}(u)= \ord(q_{u,u})$, and
$q_{u^m,u^m}=q_{u,u}^{m^2}=1$.

Let $\Delta^+(\bB)\subseteq \N^n$ be the set of degrees of the
generators of the PBW-basis, counted with their multiplicities and
let also $\Delta(\bB)= \Delta^+(\bB) \cup \left(-
\Delta^+(\bB)\right)$: $\Delta^+(\bB)$ is independent of the
choice of the PBW-basis with the property $\circledast$ (see
\cite[Lemma 2.18]{AA} for a proof of this statement).

\subsection{Auxiliary results}\label{subsection:aux}
\

Let $I$ be $\zt$-homogeneous ideal in $\bS$ and $\bB=T(V)/I$ as in
Subsection \ref{subsection:braidedHA}. We shall use repeatedly the
following fact.

\begin{obs}\label{obs:nilpotency}If $x_i^N=0$ in $\bB$ with $N$ minimal
(this is called the order of nilpotency of $x_i$), then $q_{ii}$ is a root of
1 of order $N$. Hence $(ad_c x_i)^Nx_j=0$.
\end{obs}

The following result extends (18) in     the proof of
\cite[Proposition 1]{H2}.

\begin{lema}\label{keryil}
For each $i \in \unon$, let $\cK_i$ be the subalgebra generated by
$\cup _{j \neq i} M_{i,j}(\bB)$ and denote by $n_i$ the order of
$q_{ii}$. Then there are isomorphisms of graded vector spaces

\begin{itemize}
    \item $\ker (D_i) \cong \cK_i$, if $\ord \ q_{ii}$ is the order of nilpotency of $x_i$,
    or
    \item $\ker (D_i) \cong \cK_i \t \kk \left[ x_i^{n_i}
    \right]$, if $\ord \ q_{ii}< \infty$ but $x_i$ is not
    nilpotent.
\end{itemize}
Moreover,
\begin{equation}\label{espaciosgraduados}
\bB \cong \cK_i \t \kk \left[ x_i \right].
\end{equation}
\end{lema}

\bdem We assume for simplicity $i=1$ and consider the PBW basis
obtained in the Theorem \ref{basePBW}. Now, $x_1 \in S_I$, and it
is the least element of $S_I$, so each element of $B_I'$ is of the
form $\left[u_1 \right]^{s_1}\ldots \left[ u_k \right]^{s_k}
x_1^s$, with $u_k<\ldots <u_1, u_i \in S_I\setminus \left\{ x_1
\right\}, 0<s_i<h_I(u_i), 0 \leq s <h_I(x_1)$. Call
$S'=S_I\setminus \left\{ x_1 \right\}$, and
    \[ B_2:= B \left( 1+I, \left[ S' \right]_c+I, <, h_{I} |_{S'} \right), \]
that is, the PBW set generated by $\left[ S' \right]_c +I$, whose
height is the restriction of the height of the PBW basis
corresponding to $S'$. We have $$\bB \cong \kk B_2 \t \kk \left[
x_1 \right].$$

By \eqref{14}, any $(ad_c x_1)^r(x_j)\in \ker(D_1)$; as
$D_{1}$ is a skew-derivation, we have $\cK_1 \subseteq \ker
(D_1)$.

Now, if $v \in S', v=x_{j_1}\ldots  x_{j_s}, \quad j_1,\ldots ,j_s
\geq 2$, then $\left[ v \right]_c \in \cK_1$, because it is a
homogeneous polynomial in $x_{j_1},\ldots , x_{j_s}$, and each
$x_{j_p} \in \cK_1$.

Let $v\in L$ be a word in letters $x_2,\ldots ,x_{\theta}$, of degree $\vb \in N^{\theta}$. Then $x_1v \in L$, and
    \[ \left[ x_1v \right]_c=x_1 \left[ v \right]_c - \chi (\e_1,\vb) \left[ v \right]_c x_1 = \sum _{u \geq v , \deg (u)=\vb} \alpha_u (x_1u-\chi ( \e_1,\vb) x_1), \]
where $\alpha_u \in \kk$. If $u=x_{j_1} \ldots x_{j_s}, \quad j_1,
\ldots ,j_s \in \left\{ 2,\ldots ,\theta \right\}$, we have
\begin{eqnarray*}
    x_1u-q_{x_1,u}ux_1 &=& x_1x_{j_1} \ldots  x_{j_s}-q_{1j_1} \ldots  q_{1j_s}x_{j_1} \ldots  x_{j_s}x_1
    \\ &=&  ad_c(x_1)(x_{j_1})x_{j_2}\ldots x_{j_s}+ q_{1j_1}x_{j_1}(ad_c x_1)(x_{j_2})x_{j_3}\ldots  x_{j_s}
    \\ &&+\ldots + q_{1j_1} \ldots  q_{1j_{s-1}}x_{j_1} \ldots  x_{j_{s-1}}(ad_cx_1)(x_{j_s}).
\end{eqnarray*}
Then $x_1u-q_{x_1,u}ux_1 \in \cK_1$, so $\left[ x_1v \right]_c \in
\cK_1$.

Now let $v \notin L$ be a word in letters $x_2,\ldots ,
x_{\theta}$; consider $v=u_1\ldots u_p$ its Lyndon decomposition, where $u_p \leq \ldots \leq u_1$, $u_i \in L$, $p \geq 2$. The
Shirshov decomposition of $x_1v$ is $(x_1u_1\ldots u_{p-1}, u_p)$,
so
    \[ \left[ x_1v \right]_c = \left[ x_1u_1\ldots u_{p-1} \right]_c \left[u_p\right]_c - q_{x_1u_1\ldots u_{p-1}, u_p} \left[u_p\right]_c \left[ x_1u_1\ldots u_{p-1} \right]_c, \]
and by induction on $p$ we can prove that $\left[ x_1v \right]_c
\in \cK_1$, because each $\left[ u_p \right] \in \cK_1$, and we
proved already the case $p=1$.

We next prove, by induction on $t$, that $\left[ x_1^t u \right]_c
\in \cK_1, \ \forall t \in \N$, where $u$ is a word in letters
$x_2,\ldots ,x_{\theta}$: the case $t=1$ is the previous one. Then
we consider $t \geq 2$ and $\left[ x_1^{t-1} u \right]_c \in
\cK_1$. The Shirshov decomposition of $x_1^t u$ is $(x_1,
x_1^{t-1} u)$, so $\left[x_1^t u \right]= x_1 \left[x_1^{t-1}u \right] -q_{x_1,x_1^{t-1}u} \left[ x_1^{t-1}u \right]x_1$.

By induction hypothesis, $\left[x_1^{t-1}u \right]= \sum \alpha_i
B_{1}^{(i)}\ldots  B_{n_i}^{(i)}$, for some $\alpha_i \in \kk$,
and $B_{p}^{(i)} \in \cup _{j =2}^{\theta} M_{1,j}$. using that $(ad_c x_1)$ is an skew derivation,
\begin{align*}
     x_1B_{1}^{(i)}& \cdots  B_{n_i}^{(i)}- \chi (\e_1,(t-1)\e_1 + \ub) B_{1}^{(i)} \cdots B_{n_i}^{(i)}x_1
    \\ =& (ad_c x_1)(B_{1}^{(i)})B_2^{(i)} \cdots B_{n_i}^{(i)}
    \\ &+ \chi ( \e_1, \deg B_{1}^{(i)}) B_{1}^{(i)} (ad_c x_1) (B_2^{(i)})B_3^{(i)} \cdots  B_{n_i}^{(i)}+ \ldots
    \\ &+ \chi( \e_1, \sum_{j=1}^{n_i} \deg B_{j}^{(i)})  B_{1}^{(i)} \cdots  B_{n_i-1}^{(i)}(ad_c x_1)(B_{n_i}^{(i)}).
\end{align*}
Note that if $B_p^{(i)} \in M_{1,j_l}$, then $(ad_cx_1)(B_p^{(i)})
\in M_{1,j_l}$, so $\left[ x_1^tu \right] \in \cK_1$.

For the last case, let $u \in L \setminus \left\{x_1\right\}$ be a
word that begins with the letter $x_1$ (it is the least letter);
there exist $s \geq 1, t_1,\ldots ,t_s \geq 1$ and non empty words
$u_1,\ldots ,u_s$ in letters $x_2,\ldots ,x_{\theta}$ such that
    \[u=x_1^{t_1}u_1\ldots x_1^{t_s}u_s .\]
We prove that $\left[ u \right]_c \in \cK_1$ by induction on $s$,
where the case $s=1$ is as before. So for $s>1$, if $u=vw$, where
$(v,w)$ is the Shirshov decomposition of $u$, $w$ must begin with
the letter $x_1$, because $s>1$ and $w$ is the least proper end of
$u$. Then there exists $k \in \N, 1 \leq k <s$ such that
    \[v=x_1^{t_1}u_1\ldots x_1^{t_k}u_k, \quad w=x_1^{t_{k+1}}u_{k+1} \ldots x_1^{t_s}u_s .\]
By inductive hypothesis, $\left[v\right]_c,\left[w\right]_c \in
\cK_1$, and finally $$\left[u\right]_c=
\left[v\right]_c\left[w\right]_c- \chi (\deg v, \deg w ) \left[w\right]_c\left[v\right]_c
\in \cK_1 .$$

Then we prove that $L \setminus \left\{x_1\right\} \subseteq
\cK_1$, and $B_2$ is generated by $L\setminus \left\{x_1\right\}$;
that is, $\kk B_2 \subseteq \cK_1$, and $D_1(B_2)=0$.

If $u \in \ker (D_1)$, we can write $\left[u\right]_c= \sum_{w
\in B'_I} \alpha_w \left[w\right]_c$. If $w$ does not end with
$x_1$, then $w \in B_2$, and $D_1(\left[w\right]_c)=0$. But if
$w=u_w x_1^{t_w}, \quad \left[u_w\right]_c \in B_2,
0<t_w<h_I(x_1)$, we have
    \[D_1 \left( \left[w\right]_c \right)= (t_w)_{q_{11}^{-1}}
    \left[ u_w \right]_c x_1^{t_w-1}, \]
where $(t_w)_{q_{11}^{-1}} \neq 0$ if $n_i$ does not divide $t_w$.
Then
    \[0 = D_1(\left[u\right]_c)= \sum_{w \in B_I'/ t_w>0}
    \alpha_w (t_w)_{q_{11}^{-1}} \left[ u_w \right]_c x_1^{t_w-1},    \]

But $\left[ u_w \right]_c x_1^{t_w-1} \in B_2$, and $B_2$ is a
basis, so $\alpha_w=0$ for each $w$ such that $n_i$ does not
divide $t_w$. This concludes the proof. \edem

\subsection{Transformations of certain braided graded Hopf algebras}\label{subsection:transformation}
\

Let $I$ be $\zt$-homogeneous ideal in $\bS$ and $\bB=T(V)/I$ as in
the previous Subsections. \emph{We fix $i\in \{1, \dots,
\theta\}$.}

\begin{obs} $\ord q_{ii} = \min\{k\in \N: F_i^k=0\}$.
\end{obs}

\pf Note that, if $k\in \N$, then $F_i(x_i^k)
=(k)_{q_{ii}}x_i^{k-1}$, and for all $k \in \N$,
    \[F_i^k(x_i^k)= (k) _{q_{ii}^{-1}}!. \]
That is, if $F_i^k=0$, then $(k) _{q_{ii}^{-1}}!=0$. Hence $\ord
q_{ii} \leq \min\{k\in \N: F_i^k=0\}$. Reciprocally, if $q_{ii}$
is a root of 1 of order $k$, then $F_i^k(x_i^t)= 0$ for all $t\geq
k$ by the previous claim, and $F_i^k(x_i^t)= 0$ for all $t< k$ by
degree arguments. Since $F_i(x_j)= 0$ for $j\neq i$,  $F_i^k=0$.
\epf

We now extend some considerations in \cite[p. 180]{H2}. We
consider the Hopf algebra
\begin{align*}
H_i&:= \begin{cases}\kk \langle y, e_i, e_i^{-1} \vert e_iy -
q_{ii}^{-1}y e_i, y^{N_i} \rangle & \text{where $N_i$ is the order
of nilpotency},
\\& \qquad \text{of $x_i$ in $\bB$, if $x_i$ is nilpotent};
\\ \kk \langle y, e_i, e_i^{-1} \vert e_iy -
q_{ii}^{-1}y e_i \rangle &\text{ if $x_i$ is not nilpotent};
\end{cases}
\\\de& (e_i) = e_i \t e_i, \quad \de(y )= e_i
\t y + y \t 1. \end{align*}

Notice that $\Delta$ is well-defined by Remark
\ref{obs:nilpotency}. We also consider the action $\td$ of $H_i$
on $\bB$ given by
$$
e_i\td b = \chi (\ub, \e_i ) b, \qquad  y\td b = F_i(b),
$$
if $b$ is homogeneous of degree $\ub \in \N^{\theta}$, extending the previous one defined in \eqref{actiontriangle}. The action is
well-defined by Remark \ref{obs:nilpotency} and because
    \[ \left(e_iy \right)\td b = e_i \td \left(F_i (b)\right)
= q_{ii}^{-1} F_i (e_i \td b) = \left(q_{ii}^{-1}y e_i\right)
\td b, \ \forall b \in \bB. \] It is easy to see that $\bB$ is an
$H_i$-module algebra; hence we can form
    \[\cA_i := \bB \# H_i . \]
Also, if we denote explicitly by $\cdot$ the multiplication in $\cA_i$, we
have
\begin{equation}\label{nose}
\left( 1 \# y \right) \cdot \left( b \# 1\right) = (e_i  \td b \# 1) \cdot (1 \# y) + F_i(b) \# 1, \quad \forall b \in
\bB.
\end{equation}
As in \cite{H2}, $\cA_i$ is a left Yetter-Drinfeld module over
$\kk \Gamma$, where the action and the coaction are given  by

\begin{align*}
    g_k \cdot x_j \# 1 &= q_{kj} x_j \# 1, & \delta (x_j \# 1)&=g_j \t x_j \# 1,
    \\ g_k \cdot 1 \# y &= q_{ki}^{-1} 1 \# y, &  \delta (1 \# y)&=g_i^{-1} \t 1 \# y,
    \\  g_k \cdot 1 \# e_i &= 1 \# e_i, &  \delta (1 \# e_i)&=1 \t 1 \# e_i,
\end{align*}
for each pair $k,j \in \unon$. Also, $\cA_i$ is a $\kk
\Gamma$-module algebra.

\bigbreak We now prove a generalization of \cite[Proposition
1]{H2} in the more general context of our braided Hopf algebras
$\bB$. Although the general strategy of the proof is similar as in
\emph{loc. cit.}, many points need slightly different
argumentations here.

\bigskip
\begin{theorem}\label{transfnichols} Keep the notation above.
Assume that $M_{i,j}(\bB)$ is finite and
\begin{equation}\label{conditiontransformation}
    \left| M_{i,j}(\bB) \right|=m_{ij}+1, \quad  j \in \unon, j \neq i.
\end{equation}

\emph{(i)} Let $V_i$ be the vector subspace of $\cA_i$ generated
by
    $$\left\{ (ad_cx_i)^{m_{ij}}(x_j) \# 1: j \neq i \right\} \cup \left\{ 1\# y \right\}.$$
The subalgebra $s_i(\bB)$ of $\cA_i$ generated by $V_i$ is a
graded algebra such that $s_i(\bB)^1 \cong V_i$. There exist skew
derivations $Y_i: s_i(\bB) \to s_i(\bB)$ such that, for all
$b_1,b_2 \in s_i(\bB)$, and $l, j \in \unon, j \neq i$,
\begin{eqnarray}
    && Y_j \left( b_1b_2 \right)=b_1 Y_j(b_2)+ Y_j(b_2) \left( g_i^{-m_{ij}}g_j^{-1} \cdot b_2 \right) ,\label{cond1}
    \\ && Y_i \left( b_1b_2 \right)=b_1Y_i(b_2)+ Y_i(b_1) \left( g_i^{-1} \cdot b_1 \right), \label{cond2}
    \\ && Y_l ((ad_c x_i)^{m_{ij}}(x_j)\#1)= \delta_{lj}, \qquad Y_l(1 \#y)=\delta_{li}. \label{cond3}
\end{eqnarray}

\emph{(ii)} The Hilbert series of $s_i(\bB)$ satisfies
\begin{equation}
\hlb_{s_i(\bB)} = \left( \prod_{\alpha \in \Delta^+ (\bB)
\setminus \{\e_i\}} \bq_{h_{\alpha}}(X^{s_i(\alpha)}) \right)
\bq_{h_i}(x_i).
\end{equation}

Therefore, if $s_i(\bB)$ is a graded braided Hopf algebra,
\[\de ^{+} (s_i(\bB)) = \left\{ s_i \left( \de^{+} \left( \bB \right)\right)
\setminus \left\{-\e_i\right\} \right\} \cup \left\{\e_i\right\}.
\]

\emph{(iii)} If $\bB=\bB(V)$, then the algebra $s_i(\bB)$ is
isomorphic to the Nichols algebra $\bB(V_i)$.
\end{theorem}

\bdem We prove (i). Note that $V_i$ is a Yetter-Drinfeld submodule
over $\kk \Gamma$ of $\cA_i$. Now, $\cA_i \cong \bB \t H_i$ as
graded vector spaces. Let $\cK_i$ be the subalgebra generated by
$\cup _{j \neq i} M_{i,j}(\bB)$, as in Lemma \ref{keryil}. Then
$s_i(\bB) \subseteq \cK_i \t \kk \left[ y \right] $, since $F_i$
is a skew-derivation and $F_i \left( z^{(ij)}_{k} \right) =
(k)_{q_{ii}} (1-q_{ii}^{k-1}q_{ij}q_{ji}) z^{(ij)}_{k-1}$, by
\eqref{paragenerar}. From \eqref{nose},
    \[ \left( 1 \# y \right) \cdot \left( z^{(ij)}_{m_{ij}} \# 1 \right) = \left( z^{(ij)}_{m_{ij}} \# 1 \right) \cdot \left( 1 \# y \right) + F_i \left( z^{(ij)}_{m_{ij}} \right) \# 1. \]
Also, as $m_{ij} + 1= \left|M_{i,j}(\bB)\right|$, we have
$(m_{ij})_{q_{ii}} (1-q_{ii}^{m_{ij}-1}q_{ij}q_{ji}) \neq 0$, so
$z^{(ij)}_{m_{ij}-1} \# 1 \in s_i(\bB)$, and by induction each
$z^{(ij)}_k \# 1, k=0,\ldots ,m_{ij}-1$ is an element of
$s_i(\bB)$. Then $ \cK_i \t \kk \left[ y \right] \subseteq
s_i(\bB)$, and therefore
\begin{equation}\label{bi}
     s_i(\bB) = \cK_i \t \kk \left[ y \right].
\end{equation}

Thus, $s_i(\bB)$ is a graded algebra in $\gyd$ with
$s_i(\bB)^1=V_i$. We have to find the skew derivations $Y_l \in
End(s_i(\bB))$, $l=1,\ldots , \theta$. Set $Y_i := g_i^{-1} \circ
\ad(x_i \# 1)|_{s_i(\bB)}$. Then, for each $b \in \cK_i$ and each
$j \neq i$
\begin{eqnarray*}
    && \ad(x_i \# 1)(b \# 1) = (\ad_c x_i)(b) \# 1,
    \\ && \ad(x_i \# 1) \left( (\ad_c x_i)^{m_{ij}}(x_j)\#1
\right)=(\ad_c x_i)^{m_{ij}+1}(x_j)\#1=0.
\end{eqnarray*}
Also,
\begin{eqnarray*}
    Y_i(1 \# y) &=& g_i^{-1} \cdot \left( (x_i \# 1) \cdot (1 \# y)- \left( g_i \cdot (1 \# y)\right) \cdot (x_i \# 1) \right)
    \\ &=& g_i^{-1} \cdot \left( x_i \# y + 1- q_{ii} \left( q_{ii}^{-1} x_i \# y \right) \right) =1.
\end{eqnarray*}
Thus $Y_i\in End (s_i(\bB))$ satisfies \eqref{cond3}.

Therefore, $\ad (x_i \# 1)(b_1b_2)=\ad(x_i \# 1)(b_1) b_2 +(g_i
\cdot b_1) \ad(x_i \# 1)(b_2)$, for each pair $b_1, b_2 \in
s_i(\bB)$, so we conclude that $\ad(x_i \# 1)(s_i(\bB)) \subseteq
s_i(\bB)$, and $Y_i\in End (s_i(\bB))$ satisfies \eqref{cond2}.

\smallskip
Before proving that $Y_i$ satisfies \eqref{cond1}, we need to
establish some preliminary facts.  Let us fix $j \neq i$, and let
$z^{(ij)}_k= (ad_c x_i)^k(x_j)$ as before. We define inductively
$$ \hat{z}^{(ij)}_0:=D_j, \quad \hat{z}^{(ij)}_{k+1} := D_i \hat{z}^{(ij)}_{k}-q_{ii}^kq_{ij} \hat{z}^{(ij)}_{k+1} D_i \in End(\bB). $$
We calculate

\begin{eqnarray*}
    \lambda_{ij} &:=& \hat{z}^{(ij)}_{m_{ij}} \left( z^{(ij)}_{m_{ij}} \right) = \sum^{m_{ij}}_{s=0} a_s D_i^{m_{ij}-s}D_jD_i^{s} \left( z^{(ij)}_{m_{ij}} \right)
    \\ &=& (D_i)^{m_{ij}}(D_j)\left( z^{(ij)}_{m_{ij}} \right)= \alpha_{m_{ij}} \left( m_{ij} \right)_{q_{ii}}! \in \kk^{\times},
\end{eqnarray*}
where $a_s= (-1)^k \binom{m}{k}_{q_{ii}} q_{ii}^{k(k-1)/2}q_{ij}^k$.

Note that $\left( D_i \right)^{m_{ij}+1}D_j(b)=0, \forall b \in
M_{i,k}, k \neq i,j$, and
    \[ \left( D_i \right)^{m_{ij}+1}D_j (z^{(ij)}_r)= \left(D_i \right)^{m_{ij}+1}\left( q_{ji}^{-r}\alpha_r x_i^r \right)=0, \quad \forall r \leq m_{ij}, \]
so $\left(D_i \right)^{m_{ij}+1}D_j(\cK_i)=0$. This implies that,
for each $b \in \cK_i$, $\hat{z}^{(ij)}_{m_{ij}} (b) \in \cK_i$.
Then, we define $Y_j \in End(s_i(\bB))$ by
    \[ Y_j \left( b \# y^m \right):= q_{ii}^{mm_{ij}}q_{ji}^{m}\lambda_{ij}^{-1} \hat{z}^{(ij)}_{m_{ij}} (b) \# y ^m, \quad \ b \in \cK_i, m \in \N. \]
We have $Y_j(1\# y)=0$, and if $l \neq i,j$, $Y_j ((\ad_c x_i)^{m_{il}}(x_l)\#1)=0$. By the choice of $\lambda_{ij}$, $Y_j ((\ad_c x_i)^{m_{ij}}(x_j)\#1)=1$.

Now, using that $D_k(g_l \cdot b)=q_{kl}g_l \cdot (D_k(b))$, for each $b \in \bB$, and $k,l \in \unon$, we prove inductively that for $b_1,b_2 \in \cK_i$,
$$ \hat{z}^{(ij)}_{k}(b_1b_2)= b_1 \hat{z}^{(ij)}_{k}(b_2)+ \hat{z}^{(ij)}_{k}(b_1) (g_i^kg_j \cdot b_2). $$
Then,
\begin{eqnarray*}
Y_j\left( b_1 \# 1 \cdot b_2 \# 1 \right) &=& Y_j \left( b_1 b_2 \#
1 \right) = \lambda_{ij}^{-1}\hat{z}_{m_{ij}} (b_1 b_2) \# 1
\\ &=& b_2 \# 1 \cdot Y_j ( b_2 \# 1) +  Y_j(b_1 \# 1) \cdot \left( g_i^{m_{ij}}g_j \cdot (b_2 \# 1) \right).
\end{eqnarray*}
By induction on the degree we prove that $F_i$ commute with $D_i, D_j$, so
\begin{center}
$\hat{z}^{(ij)}_{m_{ij}}( F_i( b)) = F_i \left( \hat{z}^{(ij)}_{m_{ij}}(b) \right), \quad
\forall b \in \bB$.
\end{center}
Consider $b \in \cK_i \subseteq \ker (D_i)$,
\begin{eqnarray*}
    Y_j \left( b \# 1 \cdot 1 \# y \right) &=& Y_j \left( b \# y \right)= q_{ii}^{m_{ij}}q_{ji} \hat{z}^{(ij)}_{m_{ij}} (b) \# y
    \\ &=& b \# 1 \cdot Y_j \left( 1 \# y \right) +  Y_j \left( b \# 1 \right) \cdot \left( g_i^{m_{ij}}g_j \cdot (1 \# y) \right),
\end{eqnarray*}
where we use that $Y_j(1 \# y )=0$. Then as
$$b_1 \# 1 \cdot b_2 \# y^t = b_1 \# 1 \cdot b_2 \# 1 \cdot  \left( 1\# y \right)^t, $$
\eqref{cond1} is valid for products of this form. To prove it in the general case, note that
$$ (b_1 \# y^t) \cdot (b_2 \# y^s) = (b_1 \# 1) \cdot (1 \# y)^t \cdot (b_2 \# y^s) $$

At this point, we have to prove \eqref{cond1} for $b \in \cK_i
\ker (D_i)$, $s \in \N$:
\begin{eqnarray*}
    Y_j \left( 1 \# y \cdot b \# y^s \right) &=& Y_j \left( F_i(b) \# y^s +  \left( e_i \td b \# y \right) \cdot 1\# y \right)
    \\ &=& q_{ii}^{m_{ij}s}q_{ji}^s \lambda_{ij}^{-1} \hat{z}^{(ij)}_{m_{ij}} \left( F_i(b) \right)\# y^s
    \\ && + q_{ii}^{m_{ij}(s+1)}q_{ji}^{s+1} \lambda_{ij}^{-1} \cdot \hat{z}^{(ij)}_{m_{ij}} \left( e_i \td b \right) \# y^{s+1}
    \\ &=& F_i \left(  q_{ii}^{m_{ij}(s+1)}q_{ji}^{s+1} \lambda_{ij}^{-1} \hat{z}^{(ij)}_{m_{ij}} (b) \right) \# y^s
    \\ && +  q_{ii}^{m_{ij}}q_{ji} \left( e_i \td \left( q_{ii}^{m_{ij}s}q_{ji}^{s} \lambda_{ij}^{-1} \hat{z}^{(ij)}_{m_{ij}}(b) \right) \# y^s \right)
    \\ &=& \left( 1 \# y \right) \cdot Y_j \left( b \# y^s \right)
    \\ &=& 1 \# y  \cdot Y_j \left( b \# y^s \right)  +  Y_j \left( 1 \# y \right) \cdot \left( g_i^{m_{ij}}g_j \cdot b \# y^s \right),
\end{eqnarray*}
where we use that $ \hat{z}^{(ij)}_{m_{ij}}(e_i \td b) = q_{ii}^{m_{ij}}q_{ji} e_i \td ( \hat{z}^{(ij)}_{m_{ij}}(b))$.

\bigskip
To prove (ii), note that the algebra $H_i$ is $\zt$-graded, with
$$\deg y = -\e_i, \deg e_i^{\pm 1}=0.$$
Hence, the algebra $\cA_i$ is $\zt$-graded, because $\bB$ and $H_i$ are graded, and \eqref{nose} holds.

Hence, consider the abstract basis $\left\{u_j\right\}_{j \in
\unon}$ of $V_i$, with the grading $\deg u_j=\e_j$, $\bB(V_i)$ is
$\zt$-graded. Consider also the algebra homomorphism $\Omega:
T(V_i) \rightarrow s_i(\bB)$ given by
\[ \Omega(u_j):= \left\{ \begin{array}{lc} (ad_c x_i)^{m_{ij}}
    (x_j) & j \neq i \\ y & j=i. \end{array} \right. \]
By the first part of the Theorem, $\Omega$ is an epimorphism, so
it induces an isomorphism between $s_i(\bB)':= T(V_i)/\ker \Omega$
and $s_i(\bB)$, that we also denote $\Omega$. Note:
\begin{itemize}
    \item $\deg \Omega (u_j)= \deg \left( \left( ad_c x_i \right)^{m_{ij}} (x_j) \right) = \e_j+m_{ij}\e_i= s_i(\deg \ub_j)$, if $j \neq i$;
    \item $\deg \Omega (u_i)= \deg \left( y \right) = -\e_i= s_i(\deg \ub_i)$.
\end{itemize}
As $\Omega$ is an algebra homomorphism,  $\deg (\Omega(\ub))= s_i
(\deg(\ub))$, for all $\ub \in s_i(\bB)'$. As $s_i^2=\id$, $s_i(
\deg (\Omega(\ub))) = \deg(\ub)$, for all $\ub \in s_i(\bB)'$, and
$\bH_{s_i(\bB)'}= s_i (\bH_{s_i(\bB)})$.

From this point, the proof goes exactly as in \cite[Theorem
3.2]{AA}.
\bigskip

The statement in (iii) is exactly \cite[Prop. 1]{H2}. \edem

By Theorem \ref{transfnichols}, the initial braided vector space
with matrix $(q_{kj})_{1 \leq k,j \leq \theta}$ is transformed
into another braided vector space of diagonal type $V_i$, with
matrix $(\overline{q}_{kj})_{1 \leq k,j \leq \theta}$, where
$\overline{q}_{jk}=q_{ii}^{m_{ij}m_{ik}}q_{ik}^{m_{ij}}q_{ji}^{m_{ik}}q_{jk}, \ j,k \in\unon .$

If $j \neq i$, then $\overline{m_{ij}} = \min \left\{ m \in \N :
(m+1)_{\overline{q}_{ii}}
    \left(\overline{q}_{ii}^m\overline{q}_{ij} \overline{q}_{ji}=0 \right) \right\} \overset{\!}= m_{ij}$.
\bigskip

For later use the previous Theorem in Section
\ref{section:genrel}, we recall a result from \cite{AHS}, adapted
to diagonal braided vector spaces.

\begin{lema}\label{lemaAHS}
Let $V$ a diagonal braided vector space, and $I$ an ideal of $T(V)$. Call $\bB:= T(V)/I$, and assume that there exist $(id,\sigma_i)$-derivations $D_i: \bB\rightarrow \bB$ with $D_i(x_j)=\delta_{i,j}$ for all $j$. Then, $I \subseteq I(V)$.
\end{lema}

That is, the canonical surjective algebra morphisms from $T(V)$ onto $\bB$, $\bB(V)$ induce a surjective algebra morphism $\bB \rightarrow \bB(V)$.

\bdem
See \cite[Lemma 2.8(ii)]{AHS}
\edem

\section{Standard braidings}\label{section:standard}

In \cite{H3}, Heckenberger classifies diagonal braidings whose set
of PBW generators is finite. Standard braidings form an special
subclass, which includes properly braidings of Cartan type.

we first recall the definition of standard braiding from
\cite{AA}, and the notion of Weyl groupoid, introduced in
\cite{H2}. Then we present the classification of standard
braidings, and compare them with \cite{H3}.

As in Heckenberger's works, we use the notion of \emph{generalized
Dynkin diagram associated} to a braided vector space of diagonal
type, with matrix $(q_{ij})_{1 \leq i,j \leq \theta}$: this is a
graph with $\theta$ vertices, each of them labeled with the
corresponding $q_{ii}$, and an edge between two vertices $i,j$
labeled with $q_{ij}q_{ji}$ if this scalar is different from $1$.
So two braided vector spaces of diagonal type have the same
generalized Dynkin diagram if and only if they are twist
equivalent. We shall assume that the generalized Dynkin diagram is
connected, by \cite[Lemma 4.2]{AS2}.

Summarizing, the main result of this Section says:

\begin{theorem}
Any standard braiding is twist equivalent with some of the following
\begin{itemize}
    \item a braiding of Cartan type,
    \item a braiding of type $A_{\theta}$ listed in Proposition \ref{casosAn},
    \item a braiding of type $B_{\theta}$ listed in Proposition \ref{casosBn},
    \item a braiding of type $G_{2}$ listed in Proposition \ref{casosG2}.
\end{itemize}
\end{theorem}

The generalized Dynkin diagrams appearing in Propositions
\ref{casosAn} and \ref{casosBn} correspond to the rows 1,2,3,4,5,6
in \cite[Table C]{H3}. Also, the generalized Dynkin diagrams in
Proposition \ref{casosG2} are (T8) in \cite[Section 3]{H1}.
However, our classification does not rely on \cite{H3}.

\subsection{Definitions of Weyl groupoid and standard braidings}\label{subsection:weylgroupoid}
\

Let $E = (\e_1, \dots, \e_\theta)$ be the canonical basis of
$\Z^\theta$. Consider an arbitrary matrix $(q_{ij})_{1 \leq i,j \leq \theta}\in
(\kk ^{\times})^{\theta\times\theta}$, and fix once and for all the
bilinear form $\chi: \Z^\theta \times \Z^\theta \to\kk^{\times}$
determined by
\begin{equation}\label{eq-chi}
 \chi(\e_i, \e_j) = q_{ij}, \qquad 1 \leq i,j \leq \theta.
\end{equation}

If $F= (\f_1, \dots, \f_\theta)$ is another ordered basis of
$\Z^\theta$, then we set $\qf_{ij} = \chi(\f_i,\f_j)$, $\otvz$. We
call $(\qf_{ij})$ the \emph{braiding matrix with respect to the
basis $F$}. Fix $i\in \unon$. If $1 \leq i,j \leq \theta$, then we
consider the set
$$ \widetilde{M}_{ij}:= \{m\in \N_0: (m+1)_{\qf_{ii}} \, (\qf_{ii}^m\qf_{ij}\qf_{ji} - 1)
= 0 \}. $$

If this set is nonempty, then its minimal element is denoted
$\mf_{ij}$ (which of course depends on the basis $F$). Define also
$\mf_{ii} = 2$. Let $\si\in GL(\Z^{\theta})$ be the
pseudo-reflection given by $\si(\f_j) := \f_j + \mf_{ij}\f_i,
\quad j\in \unon. $

Let $G$ be a group acting on a set $X$. We define the \emph{transformation groupoid} as $\bG = G\times X$, with the structure of groupoid given by the operation $(g,x)(h,y) = (gh, y)$ if $x=h(y)$, but undefined otherwise.

\begin{definition}\label{defi:weyl-groupoid}
Consider $\bX$ the set of all ordered bases of $\zt$, and the canonical action of $GL(\zt)$ over $\bX$. The smallest subgroupoid of the transformation groupoid $GL(\zt) \times \bX$ that satisfies the following properties:
\begin{itemize}
    \item $(\id, E)\in W(\chi)$,
    \item if $(\id, F)\in W(\chi)$ and $\si$ is defined, then $(\si, F)\in
    W(\chi)$,
\end{itemize}
is called the \emph{Weyl groupoid} $W(\chi)$ of the bilinear form $\chi$.
\end{definition}

\bigbreak Let $\bP(\chi) = \{F: (\id, F)\in W(\chi)\}$ be the set of points of the groupoid
$W(\chi)$. The set
\begin{equation}\label{eqn:root-system}
\Delta(\chi) = \bigcup_{F\in \bP(\chi)} F.
\end{equation}
is called the \emph{generalized root system} associated to
$\chi$.

\bigbreak

We record for later use the following evident facts.

\begin{obs}\label{casouno}
Let $i \in \unon$ such that $s_{i,E}$ is defined. Let $F=
s_{i,E}(E)$ and $(\widetilde{q}_{ij})$ the braiding matrix with
respect to the basis $F$. Assume that
\begin{itemize}
    \item $q_{ii}=-1$ (and then, $m_{ik}=0$ if $q_{ik}q_{ki}=1$ or $m_{ik}=1$, for each $k \neq i$);
    \item there exists $j \neq i$ such that $q_{jj}q_{ji}q_{ij}=1$ (that is, $m_{ij}=m_{ji}=1$).
\end{itemize}
Then, $\widetilde{q}_{jj}=-1$.
\end{obs}
\bdem Simply,
$\widetilde{q}_{jj}=q_{ii}q_{ij}q_{ji}q_{jj}=q_{ii}=-1$. \edem

\begin{obs}\label{verticeCartan}
If the $m_{ij}$ satisfies $q_{ii}^{m_{ij}}q_{ij}q_{ji}=1$ for all $j \neq i$, then the
braiding of $V_i$ is twist equivalent with the corresponding to $V$.
\end{obs}

Let $\alpha: W(\chi) \rightarrow GL(\theta, \Z)$, $\alpha(s,F) = s$ if $(s, F)\in W(\chi)$, and denote by $\wo$ the subgroup generated by the image of $\alpha$.

\begin{definition}\label{standardbraiding}{\cite{AA}}
We say that $\chi$ is \emph{standard} if for any $F\in \bP(\chi)$,
the integers $m_{rj}$ are defined, for all $1\le r,j\le \theta$,
and the integers $m_{rj}$ for the bases $\si(F)$ coincide with
those for $F$ for all $i,r,j$. Clearly it is enough to assume this
for the canonical basis $E$.
\end{definition}

\emph{We assume now that $\chi$ is standard.} We set $C:=(a_{ij})
\in \Z^{\theta \times \theta}$, where $a_{ij}=-m_{ij}$: it is a
generalized Cartan matrix.

\begin{prop}{\cite{AA}}\label{prop:standard}
$\wo = \langle \sE: 1 \leq i \leq \theta \rangle$. Furthermore
$\wo$ acts freely and transitively on $\bP(\chi)$. \qed
\end{prop}

Hence, $\wo$ is a Coxeter group, and $\wo$ and $\bP(\chi)$ have
the same cardinal.

\begin{lema}{\cite{AA}}\label{lema:cartan-finite}
The following are equivalent:
\begin{enumerate}
    \item The groupoid $W(\chi)$ is finite.
    \item The set $\bP(\chi)$ is finite.
    \item The generalized root system $\Delta(\chi)$ is finite.
    \item The group $\wo$ is finite.
\end{enumerate}
If $C$ is symmetrizable, (1)-(4) are equivalent to
\begin{enumerate}
    \item[(5)] The Cartan matrix $C$ is of finite type. \qed
\end{enumerate}
\end{lema}

We shall prove in Theorem \ref{heckenberger}, that if
$\Delta(\chi)$ is finite, then the matrix $C$ is symmetrizable,
hence of finite type. Then, $\bB(V)$ is of finite dimension if and
only if the Cartan matrix $C$ is of finite type.

\subsection{Classification of standard braidings}\label{subsection:classification}

We now classify standard braidings such that the Cartan matrix is
of finite type. We begin by types $C_{\theta}, D_{\theta}, E_l
\quad (l=6,7,8)$ and $F_4$: these standard braidings are
necessarily of Cartan type.

\smallskip

\begin{prop}\label{solocartan}
Let $V$ be a braided vector space of standard type, $\theta= \dim
V$, and $C=(a_{ij})_{i,j \in \unon}$ the corresponding Cartan
matrix, of type $C_{\theta}, D_{\theta}, E_l \quad (l=6,7,8)$ or
$F_4$. Then $V$ is of Cartan type (associated to the corresponding
matrix of finite type).
\end{prop}
\bdem \textit{Let $V$ be standard of type $C_{\theta}$, $\theta
\geq 3$.}
\begin{equation}\label{diagramaCn}
\xymatrix{
{\circ}^1 \ar@{-}[r] & {\circ}^2\ar@{-}[r] & {\circ}^3  \ \ \cdots  & {\circ}^{\theta-2} \ar@{-}[r] & {\circ}^{\theta-1}  & \ar@{=>}[l]{\circ}^{\theta}  }
\end{equation}

Note that, if we suppose $q_{\theta-1,\theta-1}=-1$, as
$m_{\theta-1,\theta}=2$ and $q_{\theta-1,\theta-1}^3 \neq 1$, we
have
 $$1=q_{\theta-1,\theta-1}^2q_{\theta-1,\theta}q_{\theta,\theta-1}=q_{\theta-1,\theta}q_{\theta,\theta-1},$$
so $m_{\theta,\theta-1}=m_{\theta-1,\theta}=0$, but this is a
contradiction. Then $q_{\theta-1,\theta-1} \neq -1$, and
$m_{\theta-1,\theta-2}=1$, so $q_{\theta-1,\theta-1}
q_{\theta-1,\theta-2} q_{\theta-2,\theta-1}=1$. Using Remark
\ref{casouno} when $i=\theta-2,j=\theta-1$, as
$\widetilde{q}_{\theta-1,\theta-1} \neq -1$ when we transform by
$s_{\theta-2}$ (since the new braided vector space is also standard), we have $q_{\theta-2,\theta-2} \neq -1$, so
$$q_{\theta-2,\theta-2}q_{\theta-2,\theta-1}q_{\theta-1,\theta-2}=q_{\theta-2,\theta-2}q_{\theta-2,\theta-3}q_{\theta-3,\theta-2}=1,$$
and $q_{\theta-1,\theta-1}=q_{\theta-2,\theta-2}$. Inductively,
$$q_{kk}q_{k,k-1}q_{k-1,k}=q_{kk}q_{k,k+1}q_{k+1,k}=q_{11}q_{12}q_{21}=1, \quad k=2,..., \theta-1$$
and $q_{11}=q_{22}=\ldots=q_{\theta-1,\theta-1}$. So we look at
$q_{\theta\theta}$: as $m_{\theta,\theta-1}=1$, we have
$q_{\theta\theta}=-1$ or $q_{\theta\theta}q_{\theta,\theta-1}
q_{\theta-1,\theta}=1$. If $q_{\theta\theta}=-1$, transforming by
$s_\theta$, we have
$$\widetilde{q}_{\theta-1,\theta-1}=-q^{-1}, \quad
\widetilde{q}_{\theta-1,\theta}\widetilde{q}_{\theta,\theta-1}=q^{2},$$
and as $m_{\theta-1,\theta-2} =1$, we have $q^2=-1$. Then
$$q_{\theta\theta}q_{\theta,\theta-1} q_{\theta-1,\theta}=1, \quad q_{\theta\theta}=q^{2},$$
and the braiding is of Cartan type in both cases.
\medskip

\textit{Let $V$ be standard of type $D_{\theta}$, $\theta \geq
4$}.

We prove the statement by induction on $\theta$. Let $V$ be of
standard $D_4$ type, and suppose that $q_{22}=-1$. Let
$(\widetilde{q}_{ij})$ the braiding matrix with respect to
$F=s_{2,E}(E)$. We calculate for each pair $j \neq k \in
\left\{1,3,4\right\}$:
    $$\widetilde{q}_{jk}\widetilde{q}_{kj}= \left((-1)q_{2k}q_{j2}q_{jk}\right) \left((-1)q_{2j}q_{k2}q_{kj}\right)= \left(q_{2k}q_{k2}\right)\left(q_{2j}q_{j2}\right),$$
where we use that $q_{jk}q_{kj}=1$. As also
$\widetilde{q}_{jk}\widetilde{q}_{kj}=1$, we have
$q_{2k}q_{k2}=\left(q_{2j}q_{j2}\right)^{-1}, \quad j \neq k$, so
$q_{2k}q_{k2}=-1$, $k=1,3,4$, since $q_{2k}q_{k2} \neq 1$. In this
case, the braiding is of Cartan type, with $q=-1$. Suppose then
$q_{22} \neq -1$. From the fact that $m_{2j}=1$, we have
$$q_{22}q_{2j}q_{j2}=1, \quad j=1,3,4.$$
For each $j$, applying Remark \ref{casouno}, as
$\overline{q}_{22}\neq -1$, we have $q_{jj} \neq -1$, so
$q_{jj}q_{2j}q_{j2}=1, \quad j=1,3,4$, and the braiding is of
Cartan type.

\begin{equation}\label{diagramaDn}
\xymatrix{
{\circ}^1 \ar@{-}[r] & {\circ}^2\ar@{-}[r]  & {\circ}^3 \cdots  & {\circ}^{\theta-2}\ar@{-}[d] \ar@{-}[r]  & {\circ}^{\theta}  \\ & & & {\circ}_{\theta-1} }
\end{equation}

We now suppose the statement valid for $\theta$. Let $V$ be a
standard braided vector space of type $D_{\theta+1}$. The subspace
generated by $x_2, \ldots , x_{\theta+1}$ is a standard braided
vector space associated to the matrix $(q_{ij})_{i,j=2, \ldots,
\theta+1}$, of type $D_{\theta}$, so it is of Cartan type. To
finish, apply Remark \ref{casouno} when $i=1,j=2$, so we obtain
that $V$ is of Cartan type with $q=-1$, or if $q_{22}\neq -1$, we
have $q_{11} \neq -1$, and $q_{11}q_{12}q_{21}=1$, and in this
case it is of Cartan type too (because also $q_{1k}q_{k1}=1$ when
$k>2$).
\medskip

\textit{Let $V$ be standard of type $E_6$.} Note that $1,2,3,4,5$
determine a braided vector subspace, which is standard of type
$D_5$, so it is of Cartan type. Then to prove that
$q_{66}q_{65}q_{56}=1$, we use Remark \ref{casouno} as above.

\begin{equation}\label{diagramaE6}
\xymatrix{ {\circ}^1 \ar@{-}[r] & {\circ}^2\ar@{-}[r]  &
{\circ}^{3}\ar@{-}[d] \ar@{-}[r]  & {\circ}^5\ar@{-}[r]  &
{\circ}^6    \\ & & {\circ}_4 }
\end{equation}
\medskip

\textit{If $V$ is standard of type $E_7$ or $E_8$,} we proceed
similarly by reduction to $E_6$, respectively $E_7$.

\begin{eqnarray}\label{diagramaE7}
\xymatrix{ {\circ}^1 \ar@{-}[r] & {\circ}^2\ar@{-}[r]  &
{\circ}^{3} \ar@{-}[r]  & {\circ}^4 \ar@{-}[d] \ar@{-}[r] &
{\circ}^6\ar@{-}[r]  & {\circ}^7    \\ & & & {\circ}_5 }
\\ \label{diagramaE8} \xymatrix{ {\circ}^1 \ar@{-}[r] &
{\circ}^2\ar@{-}[r]  & {\circ}^{3} \ar@{-}[r]  & {\circ}^4
\ar@{-}[r] & {\circ}^5 \ar@{-}[d] \ar@{-}[r]  & {\circ}^7
\ar@{-}[r]  & {\circ}^8   \\ & & & & {\circ}_6 }
\end{eqnarray}
\medskip

\textit{Let $V$ be standard of type $F_4$.} The vertices $2,3,4$
determine a braided subspace, which is standard of type $C_3$, so
the $q_{ij}$ satisfy the corresponding relations. Let
$(\widetilde{q}_{ij})$ the braiding matrix with respect to
$F=s_{2,E}(E)$. As $\widetilde{q}_{13}\widetilde{q}_{31}=1$ and
$q_{22}q_{23}q_{32}=1$, we have $q_{22}q_{12}q_{21}=1$.

\begin{equation}\label{diagramaF4}
\xymatrix{
{\circ}^1 \ar@{-}[r] & {\circ}^2 \ar@{=>}[r] & {\circ}^3 \ar@{-}[r] & {\circ}^4  }
\end{equation}

Now, if we suppose $q_{11}=-1$, applying Remark \ref{casouno} we
have $q_{22}=-1=q_{21}q_{12}$, and then it is corresponding vector
space of Cartan $F_4$ type associated to $q \in \G_4$. If $q_{11}\neq -1$, then
$q_{11}q_{12}q_{21}=1$, and also it is of Cartan type.

\edem

To finish the classification of standard braidings, we describe
the standard braidings that are not of Cartan type. They are
associated to Cartan matrices of type $A_{\theta}, B_{\theta}$ or
$G_2$.

We use the same notation as in \cite{H3}; $C(\theta,q;i_1,\ldots
,i_j)$ corresponds to the generalized Dynkin diagram
\begin{equation}\label{cadena}
\xymatrix{
{\circ}^1 \ar@{-}[r] & {\circ}^2 \ar@{-}[r] & {\circ}^3 \ \cdots & {\circ}^{\theta-1} \ar@{-}[r] &  {\circ}^{\theta}  }
\end{equation}
where
\begin{itemize}
    \item $q=q_{\theta-1,\theta}q_{\theta,\theta-1} q_{\theta\theta}^2$
holds, $1 \leq i_1< \ldots < i_j \leq \theta$;
    \item equation $q_{i-1,i}q_{i,i-1}=q$, where $1\leq i\leq \theta$, is
valid if and only if $i\in \{i_1,i_2,\ldots ,i_j\}$, so each
$q_{i_t,i_t}=-1, \quad t=1, \ldots, j$;
    \item $q_{ii}=q^{\pm 1}$ if $i \neq i_1, \ldots , i_j$.
\end{itemize}

Then, the labels of vertices between $i_t$ and $i_{t+1}$ are all
equal, and they are labeled with the inverse of the scalar
associated to the vertices between $i_{t+1}$ and $i_{t+2}$; the
same is valid for the scalars that appear in the edges.
\medskip

\begin{prop}\label{casosAn}
Let $V$ be a braided vector space of diagonal type. Then $V$ is
standard of $A_{\theta}$ type if and only if its generalized
Dynkin diagram is of the form:
\begin{equation} \label{cadenaAn}
C(\theta,q;i_1,\ldots ,i_j).
\end{equation}
\end{prop}

Note that the previous braiding is of Cartan type if and only if
$j=0$, or $j=n$ with $q=-1$.

\bdem
Let $V$ be a braided vector space of standard $A_{\theta}$
type. For each vertex $1<i<\theta$ we have $q_{ii}=-1$ or
$q_{ii}q_{i,i-1}q_{i-1,i}=q_{ii}q_{i,i+1}q_{i+1,i}=1$, and the
corresponding formulas for $i=1,\theta$. So suppose that
$1<i<\theta$ and $q_{ii}=-1$. We transform by $s_i$ and obtain
$$ \widetilde{q}_{i-1,i+1}=-q_{i,i+1}q_{i-1,i}q_{i-1,i+1}, \quad  \widetilde{q}_{i+1,i-1}=-q_{i,i-1}q_{i+1,i}q_{i+1,i-1},$$
and using that $m_{i-1,i+1}=\widetilde{m}_{i-1,i+1}=0$, we have
$q_{i-1,i+1}q_{i+1,i-1}=1$ and
$\widetilde{q}_{i-1,i+1}\widetilde{q}_{i+1,i-1}=1$, so we deduce
that $q_{i,i+1}q_{i+1,i}=(q_{i,i-1}q_{i-1,i})^{-1}$. Then the
corresponding matrix $(q_{ij})$ is of the form \eqref{cadenaAn}.

Now, consider $V$ of the form \ref{cadenaAn}. Assume
$q_{ii}=q^{\pm1}$; if we transform by $s_i$, then the braided
vector space $V_i$ is twist equivalent with $V$ by Remark
\ref{verticeCartan}. Thus, $\overline{m}_{ij}=m_{ij}$.

Assume $q_{ii}=-1$. We transform by $s_i$ and calculate
\begin{eqnarray*}
    \widetilde{q}_{jj} &=& (-1)^{m_{ij}^2}(q_{ij}q_{ji})^{m_{ij}}q_{jj}
        \\ &=& \begin{cases}  q_{jj}, & \left| j-i \right| >1;   \\  (-1)q^{\mp1}q^{\pm1}=-1,  &  j=i \pm 1, \ q_{jj}=q^{\pm1};  \\   (-1)q^{\pm1}(-1)=q^{\pm1},  &   j=i \pm 1, \ q_{jj}=-1.   \end{cases}
\end{eqnarray*}
Also, $\widetilde{q}_{ij}\widetilde{q}_{ji}=q_{ij}^{-1}q_{ji}^{-1}$ if $\left| j-i \right|>1$, and
    \[\widetilde{q}_{kj}\widetilde{q}_{jk} = (q_{ik}q_{ki})^{m_{ij}} (q_{ij}q_{ji})^{m_{ik}} q_{kj} q_{jk}=  \left\{ \begin{array}{ll}  q_{kj}q_{jk} & \left| j-i \right| or \left| k-i \right| >1,   \\  1  &   j=i-1,k=i+1.   \end{array} \right. \]

Then $V_i$ has a braiding of the above form too, and $(-m_{ij})$
corresponds to the finite Cartan matrix of type $A_{\theta}$, so
it is a standard braiding of type $A_{\theta}$. Thus this is the
complete family of standard braidings of type $A_{\theta}$. \edem

\begin{prop}\label{casosBn}
Let $V$ a diagonal braided vector space. Then $V$ is standard of
type $B_{\theta}$ if and only if its generalized Dynkin diagram is
of one of the following forms:

\begin{enumerate}
    \item[(a)] \Dchaintwo{}{$\zeta$}{$q^{-1}$}{$q$},    $\quad \zeta \in \G_3, \quad q\in \G_N, N \geq 4 \quad (\theta = 2)$;

    \item[(b)] \begin{picture}(44,8) \put(15,3){\oval(30,6)} \put(0,0){\makebox(30,6){\scriptsize $C(\theta-1,q^2;i_1,\ldots ,i_j)$}} \put(30,3){\line(1,0){10}} \put(41,3){\circle{2}} \put(36,4){\makebox[0pt]{\scriptsize $q^{-2}$}} \put(41,5){\makebox[0pt]{\scriptsize $q$}}
\end{picture}, $\quad  q \neq 0 -1, \quad 0 \leq j\leq d-1$;

    \item[(c)] \begin{picture}(48,8)
\put(17,3){\oval(34,6)}
\put(0,0){\makebox(34,6){\scriptsize $C(\theta-1,-\zeta ^{-1};i_1,\ldots ,i_j)$}}
\put(34,3){\line(1,0){10}}
\put(45,3){\circle{2}}
\put(40,4){\makebox[0pt]{\scriptsize $-\zeta $}}
\put(45,5){\makebox[0pt]{\scriptsize $\zeta $}}
\end{picture}, $\quad \zeta \in \G_3, \quad 0 \leq j\leq d-1$.
\end{enumerate}
\end{prop}

Note that the previous braiding is of Cartan type if and only if
it is as in (b) and $j=0$.

\bdem First we analyze the case $\theta=2$. Let $V$ a standard
braided vector space of type $B_2$. There are several
possibilities:
\begin{itemize}
\item $q_{11}^2q_{12}q_{21}=q_{22}q_{21}q_{12}=1$: this braiding
is of Cartan type, with $q=q_{11}$. Note that $q \neq -1$. This
braiding has the form (b) with $\theta=2, j=0$.

    \item $q_{11}^2q_{12}q_{21}=1, \quad q_{22}=-1$. We transform by
    $s_2$, then
        $$\widetilde{q}_{11}=-q_{11}^{-1}, \quad \widetilde{q}_{12}\widetilde{q}_{21} =q_{12}^{-1}q_{21}^{-1}. $$
Thus $\widetilde{q}_{11}^2 \widetilde{q}_{12}
\widetilde{q}_{21}=1$. It has the form (b) with $j=1$.

    \item $q_{11}^3=1, \quad q_{22}q_{21}q_{12}=1$. We transform by $s_1$,
        $$\widetilde{q}_{22}=q_{11}q_{12}q_{21}, \quad \widetilde{q}_{12}\widetilde{q}_{21} =q_{11}^2q_{12}^{-1}q_{21}^{-1}. $$
So $\widetilde{q}_{22}\widetilde{q}_{21}\widetilde{q}_{12}=1$,
which is the case (a).

    \item $q_{11}^3=1, \quad q_{22}=-1$: we transform by $s_1$,
        $$\widetilde{q}_{22}=-q_{12}^2q_{21}^2q_{11}, \quad \widetilde{q}_{12}\widetilde{q}_{21} =q_{11}^2q_{12}^{-1}q_{21}^{-1}. $$
If we transform by $s_2$,
        $$\widetilde{q}_{11}=-q_{12}q_{21}q_{11}, \quad \widetilde{q}_{12}\widetilde{q}_{21} =q_{12}^{-1}q_{21}^{-1}. $$
So $q_{12}q_{21}= \pm q_{11}$, and we discard the case
$q_{12}q_{21}= q_{11}$ because it was considered before. Then it
has the form in case (c) with $j=0$, and is standard.
\end{itemize}

Conversely, all braidings (a), (b) and (c) are standard of type
$B_2$.

Let now $V$ of type $B_{\theta}$, with $\theta \geq 3$. Note that
the first $\theta-1$ vertices determine a braiding of standard
$A_{\theta-1}$ type, and the last two determine a braiding of
standard $B_2$ type; so we have to 'glue' the possible such
braidings. The possible cases are the two presented in the
Proposition, and
\begin{center}
\rule[-4\unitlength]{0pt}{5\unitlength}
\begin{picture}(50,5)(0,3)
\put(13,3){\oval(26,6)}
\put(0,0){\makebox(26,6){\tiny $C(\theta-2,q;i_1, \ldots,i_j)$}}
\put(26,3){\line(1,0){10}}
\put(37,3){\circle{2}}
\put(38,3){\line(1,0){10}}
\put(49,3){\circle{2}}
\put(31,4){\makebox[0pt]{\scriptsize $q^{-1}$}}
\put(37,5){\makebox[0pt]{\scriptsize $q$}}
\put(43,4){\makebox[0pt]{\scriptsize $q^{-1}$}}
\put(49,5){\makebox[0pt]{\scriptsize $\zeta$}}
\end{picture} \quad .
\end{center}
But if we transform by $s_{\theta}$, we obtain
        $$\widetilde{q}_{\theta-1,\theta-1}=\zeta q^{-1}, \quad \widetilde{q}_{\theta-1, \theta-2}\widetilde{q}_{\theta-2, \theta-1} =q^{-1}, $$
so $1=\widetilde{q}_{\theta-1,\theta-1}\widetilde{q}_{\theta-1,
\theta-2}\widetilde{q}_{\theta-2, \theta-1}$ and we obtain $q= \pm
\zeta$, or $\widetilde{q}_{\theta-1,\theta-1}=-1$. Then,
$q=-\zeta$ or $q=-1$, so it is of some of the above forms.

To prove that (b), (c) are standard braidings, we use the
following fact: if $m_{ij}=0$ (that is,  $q_{ij}q_{ji}=1$) and we
transform by $s_i$, then
$$\widetilde{q}_{jj}=q_{jj}, \quad \widetilde{q}_{jk}\widetilde{q}_{jk}=q_{jk}q_{kj} \ (k \neq i).$$
In this case, if $\left|i-j\right| >1$, then $m_{ij}=0$; if $j=i
\pm 1$ we use the fact that the subdiagram determined by these two
vertices is standard of type $B_2$ or type $A_2$. So this is the
complete family of all twist equivalence classes of standard
braidings of type $B_{\theta}$. \edem

\begin{prop}\label{casosG2}
Let $V$ a braided vector space of diagonal type. Then $V$ is
standard of type $G_2$ if and only if its generalized Dynkin
diagram is one of the following:

\begin{enumerate}
    \item[(a)] \Dchaintwo{}{$q$}{$q^{-3}$}{$q^3$}, $\quad \ord \ q  \geq 4$;

    \item[(b)] There exists $\zeta \in \G_8$ such that

        \emph{(i)} \Dchaintwo{}{$\zeta ^2$}{$\zeta
        $}{$\zeta^{-1}$}, or

        \emph{(ii)} \Dchaintwo{}{$\zeta ^2$}{$\zeta ^3$}{$-1$}, or

        \emph{(iii)} \Dchaintwo{}{$\zeta$}{$\zeta^5$}{$-1$}.
\end{enumerate}
\end{prop}

Note that the previous braiding is of Cartan type iff it is as in
(a).

\bdem Let $V$ be a standard braiding of $G_2$ type. There are four
possible cases:
\begin{itemize}
\item $q_{11}^3q_{12}q_{21}=1, \quad q_{22}q_{21}q_{12}=1$: this
braiding is of Cartan type, as in (a), with $q=q_{11}$. Note that
if $q$ is a root of 1, then $\ord \ q \geq 4$ because $m_{12}=3$.

\item $q_{11}^3q_{12}q_{21}=1, \quad q_{22}=-1$: we transform by
$s_2$,
        $$\widetilde{q}_{11}=-q_{11}^{-2}, \quad \widetilde{q}_{12}\widetilde{q}_{21} =q_{12}^{-1}q_{21}^{-1}. $$
If $1=\widetilde{q}_{11}^3 \widetilde{q}_{12} \widetilde{q}_{21}
=-q_{11}^{-3}$, then $q_{12}q_{21}=-1$, and the braiding is of
Cartan type with $q_{11} \in \G_6$. If not,
$1=\widetilde{q}_{11}^{4}=q_{11}^{-8}$ and $\ord
\widetilde{q}_{11} =4$, so $\ord q_{11}=8$. Then we can express
the braiding in the form (b)-(iii).

\item $q_{11}^4=1, \quad q_{22}q_{21}q_{12}=1$: we transform by
$s_1$,
        $$\widetilde{q}_{22}=q_{11}q_{12}^2q_{21}^2, \quad \widetilde{q}_{12}\widetilde{q}_{21} =-q_{12}^{-1}q_{21}^{-1}. $$
If $1=\widetilde{q}_{22}\widetilde{q}_{21}\widetilde{q}_{12}=
-q_{11}q_{12}q_{21}$, we have $q_{11}^3q_{12}q_{21}=1$ because
$q_{11}^2=-1$, and this is a braiding of Cartan type. So we
consider now the case
$-1=\widetilde{q}_{22}=q_{11}q_{12}^2q_{21}^2$, so we obtain
$q_{22}^2=q_{11}^{-1}$ and $q_{22} \in \G_8$. Then we obtain a
braiding of the form (b)-(i).

\item $q_{11}^4=1, \quad q_{22}=-1$: we transform by $s_2$,
        $$\widetilde{q}_{11}=-q_{12}q_{21}q_{11}, \quad \widetilde{q}_{12}\widetilde{q}_{21} =q_{12}^{-1}q_{21}^{-1}. $$
If $\widetilde{q}_{11} \in \G_4$, then $(q_{12}q_{21})^4=1$.
$q_{12}q_{21} \neq 1$ and $q_{12}q_{21} \neq q_{11}^{-1}$ because
$m_{12} =3$. So, $q_{12}q_{21}=-1$ or $q_{12}q_{21}
=q_{11}=q_{11}^{-3}$, but these cases already were considered. So
we analyze the case
$$1= \widetilde{q}_{11}^3\widetilde{q}_{12}\widetilde{q}_{21}=q_{11}q_{12}^2q_{21}^2,$$
so we can express it in the form (b)-(ii) for some $\zeta \in
\G_8$.
\end{itemize}
A simple calculation proves that this braidings are of standard type, so they are all the standard braidings of $G_2$ type.
\edem

\section{Nichols algebras of standard braided vector spaces}\label{section:nichols}

In this section we study Nichols algebras associated to standard
braidings. We assume that the Dynkin diagram is connected, as in
Section \ref{section:standard}. In subsection
\ref{subsection:PBWbases} we prove that the set $\Delta^{+}
(\bB(V))$ is in bijection with $\de_C^{+}$, the set of positive
roots associated with the finite Cartan matrix $C$.

We describe an explicit set of generators in subsection
\ref{subsection:generators}, following \cite{LR}. We adapt their
proof since they work on enveloping algebras of simple Lie
algebras. In subsection \ref{subsection:dimension}, we calculate
the dimension of Nichols algebra associated to a standard braided
vector space, type by type.

\subsection{PBW bases of Nichols algebras}\label{subsection:PBWbases}
\ The next result is the analogous to \cite[Theorem 1]{H2} but for
braidings of standard type.

\begin{theorem}\label{heckenberger}
Let $V$ be a braided vector space of standard type with Cartan
matrix $C$. Then the following are equivalent:
\begin{enumerate}
    \item The set $\Delta(\bB(V))$ is finite.
    \item The Cartan matrix $C$ is of finite type.
\end{enumerate}
\end{theorem}

\bdem (1) $\Rightarrow$ (2) If $\de \left( \bB \left( V \right)
\right)$ is finite, then $\dc \subseteq \de(\bB(V))$ is also
finite.

\smallskip
$\blacktriangleright$ If $C$ is symmetrizable, (2) holds by Lemma
\ref{lema:cartan-finite}.

\smallskip
$\blacktriangleright$ Let $C$ be non symmetrizable. We prove that
either the corresponding set $\dc$ is not finite, or else there
does not exist any standard braided vector space associated with
this matrix $C$. The proof follows the same steps as in \cite{H2}
for the corresponding result about braided vector spaces of Cartan
type. The unique step where he uses the Cartan type condition is
the following, that we adapt to the standard type case. We
restrict the proof to the case where the generalized Cartan matrix
is not symmetrizable, and the corresponding Dynkin diagram is not
simply laced cycle, such that after removing an arbitrary vertex
the resulting diagram is of finite type. At this stage, as in loc.
cit., we reduce to the following cases:

\begin{itemize}
\item For $\theta=5$, there is only one multiple edge, because
there are no two multiple edges at distance one (we remove a
vertex which is not an extreme of these multiple edges and obtain
a diagram of non finite type). Then we have an unique double edge,
$$C_0=\left( \begin{array}{ccccc}   2 &   -2 &   0 &   0 &   0 \\   -1 &   2 &   -1 &   0 &   0 \\   0 &   -1 &   2 &   -1 &   0 \\   0 &   0 &   -1 &   2 &   -1 \\   -1 &   0 &   0 &   -1 &   2 \end{array} \right) . $$
Note that, if we suppose $q_{11}=-1$, from $m_{12}=2$ and
$q_{11}^3 \neq 1$, we have $1=q_{11}^2q_{12}q_{21}=q_{12}q_{21}$,
and then $m_{12}=m_{21}=0$, which is not possible. Then $q_{11}
\neq -1$, and from $m_{15}=1$, we have $q_{11}q_{15}q_{51}=1$.
Following Remark \ref{casouno} for $i=5, j=1$, we have $q_{55}
\neq -1$, because in other case, $\widetilde{q}_{11} = -1$ when we
transform by $s_5$, which is not possible (the new braiding is
also standard). Then $q_{55}q_{51}q_{15}=q_{55}q_{54}q_{45}=1$,
and $q_{11}=q_{55}$. Now, also from Remark \ref{casouno} but for
$i=4, j=5$, as $\widetilde{q}_{55} \neq -1$, it follows that
$q_{44} \neq -1$, and then
$$q_{44}q_{45}q_{54}=q_{44}q_{43}q_{34}=1, \quad q_{44}=q_{55}.$$
Following,
$$q_{33}q_{34}q_{43}=q_{33}q_{32}q_{23}=q_{22}q_{23}q_{32}=q_{22}q_{21}q_{12}=1,$$
and $q_{22}=q_{33}=q_{44}=q_{55}=q_{11}$. But then
$q_{11}q_{12}q_{21}=1$, and $m_{12}=1$, a contradiction. Then
there are no standard braidings with Cartan matrix $C_0$.
\smallskip

\item For $\theta=4$, we consider the matrix $$C=
\left(\begin{array}{cccc}   2 &   -2 &
  0 &   -b \\   -1 &
2 &   -c &   0 \\   0 &
  -f &   2 &   -d \\
-e &   0 &   -g &   2 \end{array} \right),$$ where $be,cf, dg =
1,2$, because there are no triple edges. The proof is the same as
in \cite{H2}, and we obtain that $\dc$ is infinite in this case.
\smallskip

\item For $\theta=3$, we consider the matrices $t_i$ corresponding to $s_i, \ i=1,2,3$:
\begin{eqnarray*}
t_1 &=& \left( \begin{array}{ccc} -1 & -a_{12} & -a_{13} \\ 0 & 1 & 0 \\ 0 & 0 & 1 \end{array} \right),
    \quad t_2 = \left( \begin{array}{ccc} 1 & 0 & 0 \\ -a_{21} & -1 & -a_{31} \\ 0 & 0 & 1 \end{array} \right),
    \\ t_3 &=& \left( \begin{array}{ccc} 1 & 0 & 0 \\ 0 & 1 & 0 \\ -a_{31} & -a_{32} & -1 \end{array} \right).
\end{eqnarray*}
\end{itemize}

The proof is as in \cite{H2}, and $\dc$ is infinite in this case
too.

\medskip
(2) $\Rightarrow$ (1) Let $V$ be a standard braided vector space
with Cartan matrix of finite type. Then the matrix $C$ is
symmetrizable. We fix $\pi = \{\alpha_1,\ldots, \alpha_{\theta}
\}$ a set of simple roots corresponding to the root system $\de_C$
of $C$. We define the $\Z$-linear map $$\phi: \Z \pi \rightarrow
\zt, \quad \phi(\alpha_i):=\e_i.$$

Consider the action of $s_i$ over $\de_C$ as the reflection
corresponding to the simple root $\alpha_i$; then $\phi$ is a
$W$-module morphism. Each $\beta \in \de_C$ is of the form
$\beta=w(\alpha_i)$, and $w=s_{i_1} \cdots s_{i_k}$, for some
$i_1, \ldots, i_k \in \unon$. Then
    \[ \phi(\beta)= \phi \left( s_{i_1} \cdots s_{i_k}(\alpha_i) \right) = s_{i_1} \cdots s_{i_k} \phi(\alpha_i)=s_{i_1} \cdots s_{i_k}(\e_i), \]
thus $\phi (\de_C) \subseteq \de \left( \bB(V) \right)$.

Suppose that $\de \left( \bB(V) \right) \supsetneq \phi (\de_C)$.
In this case, we give a different proof to the one in \cite{H2},
based in the fact that $C$ is positive definite. Let $\alpha$ be a
root of minimum height in the non empty set $\de \left( \bB(V)
\right) - \phi (\de_C) $. First, $\alpha \neq m\alpha_i$, for all
$m \in \N$, and $i=1,\ldots, \theta$, because $m\alpha_i \in \de
\left( \bB(V) \right) \Leftrightarrow m= \pm 1$, but $\pm \alpha_i
\in \phi (\de_C)$. Therefore, for each $s_i$, as $\alpha$ is not a
multiple of $\alpha_i$, we have $s_i(\alpha) \in \de \left( \bB(V)
\right) - \phi (\de_C) $, and then $\gr(s_i(\alpha))-\gr(\alpha)
\geq 0$. As $\alpha=\sum_{i=1}^{\theta} b_i\alpha_i$, we have
$\sum_{i=1}^{\theta} b_ia_{ij} \leq 0$, and as $b_j \geq 0$, we
have $\sum_{i,j=1}^{\theta} b_ia_{ij}b_j \leq 0$. This contradicts
the fact that $(a_{ij})$ is definite positive, and $(b_i) \geq 0,
(b_i) \neq 0$. Then $\phi (\de_C) = \de \left( \bB(V) \right)$.
\edem

\bigskip

\begin{cor}\label{corollary:genPBW}
Let $V$ be a braided vector space of standard type, $\theta= \dim
V$, and $C=(a_{ij})_{i,j \in \unon}$ the corresponding generalized
Cartan matrix of finite type. Then
\begin{enumerate}
    \item[(a)] $\phi (\de_C) = \de \left( \bB(V) \right)$, where
    as before $\phi: \Z \pi \rightarrow \zt$ is the $\Z$-linear
    map determined by $\phi(\alpha_i):=e_i$.
    \item[(b)] The multiplicity of each root in $\de$ is one.
\end{enumerate}
\end{cor}

\bdem (a) follows from the proof of (2) $\Rightarrow$ (1) of the
preceding Theorem.

Using this condition, as each root is of the form
$\beta=w(\alpha_i), w \in W, i \in \unon$, doing a certain
sequence of transformations $s_i$'s, this is the degree
corresponding to a generator of the corresponding Nichols algebra,
so the multiplicity (invariant by these transformations) is 1.
\edem

\subsection{Explicit generators for a PBW basis}\label{subsection:generators}
\

From Corollary \ref{corollary:genPBW}, we restrict our attention
to find one Lyndon word for each positive root of the root system
associated with the corresponding finite Cartan matrix.

\begin{prop}{\cite[Proposition 2.9]{LR}}
Let $l$ be an element of $S_I$. Then $l$ is of the form $l=l_1
\ldots l_k a$, where
\begin{itemize}
    \item $l_i \in S_I$, for each $i =1, \ldots , k$;
    \item $l_i$ is a beginning of $l_{i-1}$, for each $i>1$;
    \item $a$ is a letter.
\end{itemize}
Also, if $l=uv$ is the Shirshov decomposition, then $u,v \in
S_{I}$. \qed
\end{prop}

In what follows, we describe a set of Lyndon words for each Cartan
matrix of finite type $C$.

Consider $\alpha = \sum^{\theta}_{j=1} a_j \alpha_j \in \de^{+}$,
let $l_{\alpha} \in S_{I}$ be such that $\deg{l_{\alpha}}=\alpha$.
Let $l_{\alpha}=l_{\beta_1} \ldots l_{\beta_k} x_s$ be a
decomposition as above, where $s \in \{1, \ldots , \theta \}$ and
$\deg{l_{\beta_j}}=\beta_j$. Note that, as each $l_{\beta_j}$ is a
beginning of $l_{\beta_{j-1}}$, all the words begin with the same
letter $x'$, and as $l$ is a Lyndon word, $x'<x_s$. Therefore,
$x'$ is the least letter of $l$, so
$$x'=x_i, \quad i=\min \{j: a_j \neq 0 \} \qquad \Rightarrow \quad
\alpha = \sum^{\theta}_{j=i} a_j \alpha_j .$$ Then $k \leq a_i
\leq 3$ -- for the order given in \eqref{cadena},
\eqref{diagramaDn}, \eqref{diagramaE6}, \eqref{diagramaE7},
\eqref{diagramaE8}, \eqref{diagramaF4} ($a_i=3$ appears only when
$C$ is of type $G_2$).

Now, each $l_{\beta_j} \in S_I$, so $\beta_j \in \de^{+}$; i.e.,
it corresponds with a term of the PBW basis. Also $\sum^{k}_{j=1}
\beta_j + \alpha_s=\alpha$. If $k=2$, we have $\beta_1-\beta_2=
\sum^{\theta}_{j=1} b_j \alpha_j, \quad b_j \geq 0$, because
$\beta_2$ is a beginning of $\beta_1$ (the analogous claim is
valid when the matrix is of type $G_2$, and $k=3$). With these
rules we define inductively Lyndon words for a PBW basis
corresponding with a standard braiding for a fixed order on the
letters as in \cite{LR}, but taking care that in their work they
use the Serre relations. Now we have Serre quantum relations and
some quantum binomial coefficients maybe are zero. \bigbreak

\textbf{Type $A_{\theta}$:} In this case, the roots are of the
form
    $$\ub_{i,j}:=\sum^{j}_{k=i} \alpha_k, \quad 1 \leq i \leq j \leq \theta.$$
    By induction on $s=j-i$, we have
    $$l_{\ub_{i,j}}=x_i x_{i+1} \ldots x_j.$$
    This is because when $s=0$ we have $i=j$, and the unique possibility is
    $l_{\ub_{i,i}}=x_i$. Then if we remove the last letter (when $j-i>0$),
    we must obtain a Lyndon word, so the last letter must be $x_j$.

\bigbreak \textbf{Type $B_{\theta}$:} For convenience, we use the
following enumeration of vertices:

\begin{equation}\label{diagramaBn}
\xymatrix{ {\circ}^1 \ar@{<=}[r] & {\circ}^2 \ar@{-}[r] &
{\circ}^3 \cdots & {\circ}^{\theta-1} \ar@{-}[r] &
{\circ}^{\theta} }.
\end{equation}

The roots are of the form $\ub_{i,j}:=\sum^{j}_{k=i} \alpha_k$, or
$$\vb_{i,j}:= 2 \sum^{i}_{k=1} \alpha_k +
\sum^{j}_{k=i+1} \alpha_k.$$ In the first case, as above we have
$l_{\ub_{i,j}}=x_i x_{i+1} \ldots x_j$. In the second case,
 note that if $j=i+1$, we must have $x_{i+1}$ as the last letter
 to obtain a decomposition in two words $x_1 \cdots x_i$; if
 $j>i+1$, then the last letter must be $x_j$, so we obtain that
$$l_{\vb_{i,j}}=x_1 x_2 \ldots x_{i} x_1 x_2 \ldots x_j.$$

\bigbreak \textbf{Type $C_{\theta}$:} The roots are of the form
$\ub_{i,j}:=\sum^{j}_{k=i} \alpha_k$, or
    $$\wb_{i,j}:=\sum^{j-1}_{k=i} \alpha_k+ 2 \sum^{\theta-1}_{k=j} \alpha_k + \alpha_{\theta}, \quad i \leq j < \theta . $$
As before, $l_{\ub_{i,j}}=x_i x_{i+1} \ldots x_j$. Now, if $i<j$,
the least letter $x_i$ has degree 1, so if we remove the last
letter, we obtain a Lyndon word; i. e., $\wb_{i,j}-x_s$ is a root,
and then $x_s=x_j$, so
    $$l_{\wb_{i,j}}=x_ix_{i+1} \ldots x_{\theta-1} x_{\theta} x_{\theta-1} \ldots x_j .$$
When $i=j$, $a_i=2$, so there are one or two Lyndon words
$\beta_j$ as before. As $\wb -x_s$ is not a root, for
$s=i+1,...,\theta$, and $i<s$, there are two Lyndon words $\beta_1
\geq \beta_2$, and $\beta_1+\beta_2= 2 \sum^{\theta-1}_{k=i}
\alpha_k$. The unique possibility is $\beta_1=\beta_2=x_i x_{i+1}
\ldots x_{\theta-1}$; i. e.,
$$l_{\wb_{i,i}}=x_ix_{i+1} \ldots x_{\theta-1} x_ix_{i+1} \ldots x_{\theta-1} x_{\theta}.$$

\bigbreak   \textbf{Type $D_{\theta}$:} the roots are of the form
$\ub_{i,j} :=\sum^{j}_{k=i} \alpha_k, \quad 1 \leq i \leq j \leq
\theta$, or
\begin{eqnarray*}
    \zb_{i,j} &:=& \sum^{j-1}_{k=i} \alpha_k+ 2 \sum^{\theta-2}_{k=j} \alpha_k + \alpha_{\theta-1}+ \alpha_{\theta}, \quad i < j \leq \theta-2,
    \\ \bar{\zb}_i &:=& \sum^{\theta-2}_{k=i} \alpha_k + \alpha_{\theta}, \quad 1 \leq i \leq \theta -2.
\end{eqnarray*}
As above, $l_{\ub_{i,j}}=x_i x_{i+1} \ldots x_j$ if $j \leq n-1$.
When the roots are of type $\bar{\zb}_i$, as $\bar{\zb}_i-x_s$
must be a root (if $x_s$ is the last letter), we have $s=\theta$,
and then $l_{\bar{\zb}_i}=x_i x_{i+1} \ldots x_{\theta-2}
x_{\theta}$ is the unique possibility.

Now, when $\alpha= \ub_{i,\theta}$, the last letter is
$x_{\theta-1}$ or $x_{\theta}$: if it is $x_{\theta}$, we have
$l_{\ub_{i,\theta}}=x_i x_{i+1} \ldots x_{\theta-1} x_{\theta}$.
As $m_{\theta-1, \theta}=0$, we have
$x_{\theta-1}x_{\theta}=q_{\theta-1,
\theta}x_{\theta}x_{\theta-1}$, so $$ x_i x_{i+1} \ldots
x_{\theta-1} x_{\theta} \equiv x_i x_{i+1} \ldots
x_{\theta-2}x_{\theta}x_{\theta-1} \quad \mod I,$$ and then $x_i
x_{i+1} \ldots x_{\theta-1} x_{\theta} \notin S_I$. So, $l_{\ub_{i,\theta}}=x_i \ldots x_{\theta-2}x_{\theta}x_{\theta-1}$.

In the last case, note that if $j=n-2$, the unique possibility is
$\beta_t$ as before, because the least letter $x_i$ has degree 1
and as $\alpha-\alpha_s$ is a root, $x_s=x_{\theta-2}$. Then
$l_{\zb_{i,\theta-2}}= x_i \ldots
x_{\theta-2}x_{\theta}x_{\theta-1}x_{\theta-2}$, and inductively,
$$l_{\zb_{i,j}}= x_i \ldots
x_{\theta-2}x_{\theta}x_{\theta-1}x_{\theta-2} \ldots x_j. $$

\bigbreak \textbf{Type $E_6$:} Note that if $\alpha= \sum_{j=1}^6
a_j\alpha_j$ and $a_6=0$, then it corresponds with the Dynkin
subdiagram of type $D_5$ determined by $1,2,3,4,5$, and we obtain
$l_{\alpha}$ as above. If $a_1=0$ it corresponds with the Dynkin
subdiagram of type $D_5$ determined by $2,3,4,5,6$ -- the
numeration is different of the one given in \ref{diagramaDn}.
Anyway, the roots are defined in a similar way, and we obtain the
same list as in \cite[Fig.1]{LR}. If $a_4=0$, then $\alpha$
corresponds with the Dynkin subdiagram of type $A_5$ determined by
$1,2,3,5,6$.

So we restrict our attention to the case $a_i \neq 0, \quad
i=1,2,3,4,5,6$. We consider each case:
\begin{itemize}
    \item $\alpha=\alpha_1+\alpha_2+\alpha_3+\alpha_4+\alpha_5+\alpha_6$: as $a_1=1$,
    $\alpha-\alpha_s= \beta_1$ is a root, where $\alpha_s$ is the last letter. Then
    $s=2$ or $s=6$. In the second case, $l_{\beta_1}=x_1x_2x_3x_4x_5$, but using that
    $x_2x_3=q_{23}x_3x_2$, we have that $x_1x_2x_3x_4x_5 \notin S_I$. So $s=2$, and
    $l_{\alpha}=x_1x_3x_4x_5x_6x_2$.

    \item $\alpha=\alpha_1+\alpha_2+\alpha_3+2\alpha_4+\alpha_5+\alpha_6$: from $a_1=1$,
    we note that $\alpha-\alpha_s= \beta_1$ is a root. Then $s=4$, and
    $l_{\alpha}=x_1x_3x_4x_5x_6x_2x_4$.

    \item $\alpha=\alpha_1+\alpha_2+2\alpha_3+2\alpha_4+\alpha_5+\alpha_6$: from $a_1=1$,
    $\alpha-\alpha_s= \beta_1$ is a root. So $s=3$, and
    $l_{\alpha}=x_1x_3x_4x_5x_6x_2x_4x_3$.

    \item $\alpha=\alpha_1+\alpha_2+\alpha_3+2\alpha_4+2\alpha_5+\alpha_6$: from $a_1=1$,
    $\alpha-\alpha_s= \beta_1$ is a root. The unique possibility is $s=5$, and
    $l_{\alpha}=x_1x_3x_4x_5x_6x_2x_4x_5$.

    \item $\alpha=\alpha_1+\alpha_2+2\alpha_3+2\alpha_4+2\alpha_5+\alpha_6$: as above
    $a_1=1$, and $\alpha-\alpha_s= \beta_1$ is a root. So $s=3$, and
    $l_{\alpha}=x_1x_3x_4x_5x_6x_2x_4x_5x_3$.

    \item $\alpha=\alpha_1+\alpha_2+2\alpha_3+3\alpha_4+2\alpha_5+\alpha_6$: from $a_1=1$,
    $\alpha-\alpha_s= \beta_1$ is a root. Then $s=4$ and $l_{\alpha}=x_1x_3x_4x_5x_6x_2x_4x_5x_3x_4$.

    \item $\alpha=\alpha_1+2\alpha_2+2\alpha_3+3\alpha_4+2\alpha_5+\alpha_6$: from $a_1=1$,
    $\alpha-\alpha_s= \beta_1$ is a root. So $s=2$, and $l_{\alpha}=x_1x_3x_4x_5x_6x_2x_4x_5x_3x_4$.
\end{itemize}

\bigbreak \textbf{Type $E_7$:} If $\alpha= \sum_{j=1}^7
a_j\alpha_j$ and $a_7=0$, the root corresponds
    to the subdiagram of type $D_6$ determined by $1,2,3,4,5,6$, and we obtain $l_{\alpha}$
    as above. If $a_1=0$, it corresponds to the subdiagram of type $E_6$ determined by
    $2,3,4,5,6,7$. If $a_5=0$, then $\alpha$ corresponds to the subdiagram of type $A_6$
    determined by $1,2,3,4,6,7$.

As above, consider each case where $a_i \neq 0, \quad
i=1,2,3,4,5,6,7$:
\begin{itemize}
        \item $\alpha=\alpha_1+\alpha_2+\alpha_3+\alpha_4+\alpha_5+\alpha_6+ \alpha_7$:
        as $a_1=1$, $\alpha-\alpha_s= \beta_1$ is a root, if $\alpha_s$ is the last letter.
        Then $s=2$ o $s=7$. In the second case, $l_{\beta_1}=x_1x_2x_3x_4x_5x_6$, but
        from $x_2x_3=q_{23}x_3x_2$, we have $x_1x_2x_3x_4x_5x_6x_7 \notin S_I$. So $s=2$,
        and $l_{\alpha}=x_1x_3x_4x_5x_6x_7x_2$.

    \item $\alpha=\alpha_1+\alpha_2+\alpha_3+2\alpha_4+\alpha_5+\alpha_6+ \alpha_7$: now,
    $s=4,7$. We discard the case $s=7$ using that $m_{47}=0$, and then $s=4$:
    $l_{\alpha}=x_1x_3x_4x_5x_6x_7x_2x_4$.

    \item $\alpha=\alpha_1+\alpha_2+2\alpha_3+2\alpha_4+\alpha_5+\alpha_6+ \alpha_7$: as
    above, $s=3,7$, but we discard $s=7$ using that $m_{37}=0$, so
    $l_{\alpha}=x_1x_3x_4x_5x_6x_7x_2x_4x_3$.

    \item $\alpha=\alpha_1+\alpha_2+\alpha_3+2\alpha_4+2\alpha_5+\alpha_6+ \alpha_7$: now,
    $s=5,7$, and discard the case $s=7$ because $m_{57}=0$, and $l_{\alpha}=x_1x_3x_4x_5x_6x_7x_2x_4x_5$.

    \item $\alpha=\alpha_1+\alpha_2+2\alpha_3+2\alpha_4+2\alpha_5+\alpha_6+ \alpha_7$:
    $s=3,7$, and as above we discard the case $s=7$, so
    $l_{\alpha}=x_1x_3x_4x_5x_6x_7x_2x_4x_5x_3$.

    \item $\alpha=\alpha_1+\alpha_2+2\alpha_3+3\alpha_4+2\alpha_5+\alpha_6+ \alpha_7$:
    $s=4$, and $l_{\alpha}=x_1x_3x_4x_5x_6x_7x_2x_4x_5x_3x_4$.

    \item $\alpha=\alpha_1+2\alpha_2+2\alpha_3+3\alpha_4+2\alpha_5+\alpha_6+ \alpha_7$:
    $s=2$, as above, and $l_{\alpha}=x_1x_3x_4x_5x_6x_7x_2x_4x_5x_3x_4x_2$.

    \item $\alpha=\alpha_1+\alpha_2+\alpha_3+2\alpha_4+2\alpha_5+2\alpha_6+
    \alpha_7$:as above, the unique possibility is $s=6$, so
    $l_{\alpha}=x_1x_3x_4x_5x_6x_7x_2x_4x_5x_6$.

    \item $\alpha=\alpha_1+\alpha_2+2\alpha_3+2\alpha_4+2\alpha_5+2\alpha_6+ \alpha_7$:
    $s=3$, and $l_{\alpha}=x_1x_3x_4x_5x_6x_7x_2x_4x_5x_6x_3$.

    \item $\alpha=\alpha_1+\alpha_2+2\alpha_3+3\alpha_4+2\alpha_5+2\alpha_6+ \alpha_7$:
    $s=4$, and $l_{\alpha}=x_1x_3x_4x_5x_6x_7x_2x_4x_5x_6x_3x_4$.

    \item $\alpha=\alpha_1+2\alpha_2+2\alpha_3+3\alpha_4+2\alpha_5+2\alpha_6+ \alpha_7$:
    $s=2$, and $l_{\alpha}=x_1x_3x_4x_5x_6x_7x_2x_4x_5x_6x_3x_4x_2$.

    \item $\alpha=\alpha_1+\alpha_2+2\alpha_3+3\alpha_4+3\alpha_5+2\alpha_6+ \alpha_7$:
    $s=5$, and $l_{\alpha}=x_1x_3x_4x_5x_6x_7x_2x_4x_5x_6x_3x_4x_5$.

    \item $\alpha=\alpha_1+2\alpha_2+2\alpha_3+3\alpha_4+3\alpha_5+2\alpha_6+ \alpha_7$: as above,
    $s=2$, and $l_{\alpha}=x_1x_3x_4x_5x_6x_7x_2x_4x_5x_6x_3x_4x_5x_2$.

    \item $\alpha=\alpha_1+2\alpha_2+2\alpha_3+4\alpha_4+3\alpha_5+2\alpha_6+ \alpha_7$:
    $s=4$, and $l_{\alpha}=x_1x_3x_4x_5x_6x_7x_2x_4x_5x_6x_3x_4x_5x_2x_4$.

    \item $\alpha=\alpha_1+2\alpha_2+3\alpha_3+4\alpha_4+3\alpha_5+2\alpha_6+ \alpha_7$:
    $s=3$, and $l_{\alpha}=x_1x_3x_4x_5x_6x_7x_2x_4x_5x_6x_3x_4x_5x_2x_4x_3$.

    \item $\alpha=2\alpha_1+2\alpha_2+3\alpha_3+4\alpha_4+3\alpha_5+2\alpha_6+ \alpha_7$:
    now, there are one or two words $\beta_j$. As $\alpha -\alpha_s \in \de^{+}$ iff $s=1$
    and $x_1$ is not the last letter (because it is the least letter), there are two words
    $\beta_{j}$. So looking at the roots we obtain $s=7$, and
    $$l_{\alpha}= (x_1x_3x_4x_5x_6x_2x_4x_5x_3x_4x_2)(x_1x_3x_4x_5x_6)x_7$$
\end{itemize}

\bigbreak \textbf{Type $E_8$:} Consider $\alpha= \sum_{j=1}^8
a_j\alpha_j$; if $a_8=0$, the root corresponds to the subdiagram
of type $D_7$ determined by $1,2,3,4,5,6,7$, and we obtain
$l_{\alpha}$ as in that case. If $a_1=0$, it corresponds to the
subdiagram of type $E_7$ determined by $2,3,4,5,6,7,8$. If
$a_6=0$, then $\alpha$ corresponds to a subdiagram of type $A_7$
determined by $1,2,3,4,5,7,8$.

So, we consider the case $a_i \neq 0, \quad i=1,2,3,4,5,6,7,8$,
and solve it case by case in a similar way as for $E_7$, by
induction on the height.

\bigbreak \textbf{Type $F_4$:} Now, $\alpha= \sum_{j=1}^4
a_j\alpha_j$. If $a_4=0$, then it corresponds to the subdiagram of
type $B_3$ determined by $1,2,3$, so we obtain $l_{\alpha}$ as
before. If $a_1=0$, $\alpha$ corresponds to the subdiagram of
type $C_3$ determined by $2,3,4$.

So consider the case $a_i \neq 0, \quad i=1,2,3,4$:
\begin{itemize}
    \item $\alpha=\alpha_1+\alpha_2+\alpha_3+\alpha_4$: $a_1=1$, so
    $\alpha-\alpha_s= \beta_1$ is a root, where $\alpha_s$ is the last letter. Then
    $s=4$, and $l_{\alpha}=x_1x_2x_3x_4$.

    \item $\alpha=\alpha_1+\alpha_2+2\alpha_3+\alpha_4$: $a_1=1$, so $\alpha-
    \alpha_s= \beta_1$ is a root. Now, $s=3$ or $s=4$. If $s=4$, then
    $l_{\alpha}=x_1x_2x_3^2x_4$. But $m_{34}=2$, so
    $$x_3^2x_4 \equiv q_{34}(1+q_{33})x_3x_4x_3-q_{33}q_{34}x_4x_3^2 \quad \mod I,$$
    and $x_1x_2x_3^2x_4 \notin S_I$, a contradiction. So $s=3$, and we have that
    $l_{\alpha}=x_1x_2x_3x_4x_3$.

    \item $\alpha=\alpha_1+2\alpha_2+2\alpha_3+\alpha_4$: $a_1=1$, and as above, $s=2$
    or $s=4$: if $s=4$, then $l_{\alpha}=x_1x_2x_3^2x_2x_4$, but it is not an element of
    $S_I$, because $x_2x_4 \equiv q_{24}x_2x_4 \quad \mod I$. Then $s=2$, and
    $l_{\alpha}=x_1x_2x_3x_4x_3x_2$.

    \item $\alpha=\alpha_1+2\alpha_2+3\alpha_3+\alpha_4$: $a_1=1$, so $s=3$, and we have that
    $l_{\alpha}=x_1x_2x_3x_4x_3x_2x_3$.

    \item $\alpha=\alpha_1+\alpha_2+2\alpha_3+2\alpha_4$: $a_1=1$, so $s=4$, and
    $l_{\alpha}=x_1x_2x_3x_4x_3x_4$.

    \item $\alpha=\alpha_1+2\alpha_2+2\alpha_3+2\alpha_4$: $a_1=1$, so $s=2$ or $s=4$,
    but we discard the case $s=4$ since $x_2x_4 \equiv q_{24}x_2x_4 \quad \mod I$. So,
    $l_{\alpha}=x_1x_2x_3x_4x_3x_4x_2$.

    \item $\alpha=\alpha_1+2\alpha_2+3\alpha_3+2\alpha_4$: $a_1=1$, so $s=3$, and
    $$l_{\alpha}=x_1x_2x_3x_4x_3x_4x_2x_3.$$

    \item $\alpha=\alpha_1+2\alpha_2+4\alpha_3+2\alpha_4$: $a_1=1$, so $s=3$, and
    $$l_{\alpha}=x_1x_2x_3x_4x_3x_4x_2x_3^2.$$

    \item $\alpha=\alpha_1+3\alpha_2+4\alpha_3+2\alpha_4$: $a_1=1$, so $s=2$, and
    $$l_{\alpha}=x_1x_2x_3x_4x_3x_4x_2x_3^2x_2.$$

    \item $\alpha=2\alpha_1+3\alpha_2+4\alpha_3+2\alpha_4$: $a_1=2$, then there are one
    or two Lyndon words $\beta_j$. If there is only one, $\beta_1=\alpha-\alpha_s
    \in \de^{+}$. The unique possibility is $s=1$, but it contradicts that
    $l_{\alpha}$ is a Lyndon word. Then there exist $\beta_1, \beta_2 \in \de^{+}$ such
    that $\beta_1+\beta_2=\alpha-\alpha_s$, and $\beta_2$ is a beginning of $\beta_1$.
    So $s=2$ and $\beta_1=\beta_2= \alpha_1 +\alpha_2 +2\alpha_3 +\alpha_4$,
    i.e., $ \quad l_{\alpha}=x_1x_2x_3x_4x_3x_1x_2x_3x_4x_3x_2$.
\end{itemize}

\bigbreak  \textbf{Type $G_2$:} the roots are $\alpha_1,
\alpha_2,\alpha_1+\alpha_2,
    2\alpha_1+\alpha_2, 3\alpha_1+\alpha_2,3\alpha_1+2\alpha_2$:
    $$l_{\alpha_1}=x_1, \qquad l_{\alpha_2}=x_2, \qquad
    l_{m\alpha_1+\alpha_2}=x_1^m x_2, \quad m=1,2,3.$$
If $\alpha= 3\alpha_1+2\alpha_2$, the last letter is $x_2$. If we
suppose $\beta_1=3\alpha_1+\alpha_2$, then
$l_{\alpha}=x_1^3x_2^2$, but
$$(ad x_2)^2 x_1 = x_2^2x_1-q_{21}(1+q_{22})x_2x_1x_2+q_{22}q_{21}x_1x_2^2
\equiv 0 \quad \mod I,$$ so we have
$$x_1^3x_2^2 \equiv (q_{22}^{-1}+1)x_1^2x_2x_1x_2 - q_{22}^{-1}q_{21}^{-1}
x_1^2x_2^2x_1 \quad \mod I,$$ and then $l_{\alpha}=x_1^3x_2^2
\notin S_I$ because $q_{22}^{-1}q_{21}^{-1} \neq 0$, so there are
at least two words $\beta_{j}$. Analogously, if we suppose that
there are three words $\beta_j$, as $\beta_1 \geq \beta_2 \geq
\beta_3$ and $\beta_1 + \beta_2 + \beta_3= 3\alpha_1+\alpha_2$, we
have $l_{\beta_1}=l_{\beta_2}=x_1 > l_{\beta_3}=x_1x_2$, and also
$l_{\alpha}=x_1^3x_2^2 \notin S_I$. So there are two Lyndon words
of degree $\beta_1 \geq \beta_2$, so the unique possibility is
$\beta_1=2\alpha_1+\alpha_2, \quad \beta_2=\alpha_1$; i. e.,
$l_{\alpha}=x_1^2x_2x_1x_2$.

\subsection{Dimension of Nichols algebras of standard braidings}\label{subsection:dimension}
\

We begin by standard braidings of types $C_{\theta}, D_{\theta},
E_6, E_7, E_8, F_4$, which are just of Cartan type.

\begin{prop}
Let $V$ a braided vector space of Cartan type, where $q_{44} \in
\G_{N}$ if $V$ is of type $F_4$, or $q_{11} \in \G_N$ otherwise,
for some $N \in \N$. Then, for the associated Nichols algebra
$\bB(V)$, we give
\begin{description}
    \item[Type $C_{\theta}$] $\quad \dim \bB(V)= \begin{cases}N^{\theta^2} & N \mbox{ odd,} \\ N^{\theta^2}/2^{\theta} & N \mbox{ even;}  \end{cases} $

    \item[Type $F_{4}$] $\quad \dim \bB(V)= \begin{cases}N^{24} & N \mbox{ odd,} \\ N^{24}/2^{12} & N \mbox{ even;}  \end{cases} $

    \item[Types $D_{\theta}, E_6,E_7,E_8$] $\quad \dim \bB(V)= N^{\left| \de^{+} \right|}$.
\end{description}
\end{prop}
Note that the last case corresponds to simply-laced Dynkin
diagrams.

\bdem Note that if $N$ is odd, then $\ord q^2 = \ord q=N$, but if $N$ is even, we have $\ord q^2 = N/2$. Also, as the braiding is of Cartan type,
$$ q_{s_i(\alpha)}= \chi \left( s_i(\alpha) , s_i(\alpha) \right)=  \widetilde{\chi} (\alpha, \alpha ) = \chi (\alpha , \alpha ) = q_{\alpha}. $$
Using this, we just have to determine how many roots there are in
the orbit of each simply root.
\smallskip

\textit{When $V$ is of type $C_{\theta}$}, $q_{ii}=q$, except for $q_{\theta\theta} = q^2$. The roots in the orbit of $\alpha_{\theta}$ by the action of the Weyl group are $q_{\wb_{ii}}$ for $1 \leq i < \theta$, and the others are in the orbit of $\alpha_j$, for some $j < \theta$. Then, there are $\theta$ roots such that $q_{\alpha}=q^2$, and $q_{\alpha}=q$ for the rest.
\smallskip

\textit{When $V$ is of type $F_4$}, we have $q_{11}=q_{22}=q^2$, and $q_{33}=q_{44}=q$. There are exactly 12 roots in the union of orbits corresponding to $\alpha_1$ and $\alpha_2$, and the other 12 in the union of orbits corresponding to $\alpha_3$ and $\alpha_4$. So
$$ \left| \left\{ \alpha \in \de^+ :  q_{\alpha}=q \right\} \right| = \left| \left\{ \alpha \in \de^+ :  q_{\alpha}=q^2 \right\} \right| = 12. $$
\smallskip

\textit{When $V$ is of type $D$ or $E$}, all $q_{\alpha}=q$ because $q_{ii}=q$, for all $1 \leq i \leq \theta$.
\smallskip

The formula for the dimension is a consequence of the theory of
PBW bases above and Corollary \ref{corollary:genPBW}. \edem

Now we treat the types $A_{\theta}, B_{\theta}$ and $G_2$.

\begin{prop}
Let $V$ be a standard braided vector space of type $A_{\theta}$ as
in Proposition \ref{casosAn}. Then the associated Nichols algebra
$\bB(V)$ is of finite dimension if and only if $q$ is a root of
unit of order $N \geq 2$. In such case,
\begin{equation}\label{dimensionAn} \dim \bB(V)= 2^{ \binom{\theta +1 }{ 2
} -  \binom{t}{2} - \binom{\theta +1-t}{2}}N^{\binom{t}{2} +
\binom{\theta +1-t}{2}},
\end{equation}
where $ t= \theta- \sum ^{j}_{k=1}(-1)^{j-k}i_k $.
\end{prop}

\bdem $q$ is a root of unit of order $N \geq 2$ because the height
of each PBW generator is finite. To calculate the dimension,
recall that from Corollary \ref{corollary:genPBW}, we have to
determine $q_{\alpha}$ for $\alpha \in \de_C$. As before,
$\ub_{ij} = \sum^{j}_{k=i} \e_k, \quad i \leq j$, and we have
$$\de(\bB(V))= \{ \ub_{ij}: 1\leq i \leq j \leq \theta \}.$$
If $1\leq i \leq j \leq \theta$, we define $\quad \kappa_{ij}:= \sharp \{ k \in \{i, \ldots , j\}: q_{kk}=-1 \}$.

We prove by induction on $j-i$ that
\begin{itemize}
    \item if $\kappa_{ij}$ is odd, then $q_{\ub_{ij}}=-1$;
    \item if $\kappa_{ij}$ is even, then $q_{\ub_{ij}}=q_{i,i+1}^{-1}q_{i+1,i}^{-1}$.
\end{itemize}
If $j-i=0$, then $q_{\ub_{ii}}=q_{ii}$; in this case,
$\kappa_{ii}=1$ if $q_{ii}=-1$ or $\kappa_{ii}=0$ if
$q_{ii}=(q_{i,i+1}q_{i+1,i})^{-1} \neq -1$. Now, assume this is
valid for certain $j$, and calculate it for $j+1$:
\begin{align*}
    q_{\ub_{i,j+1}} =& \chi(\ub_{ij}+\e_{j+1},\ub_{ij}+\e_{j+1})=q_{\ub_{ij}} \chi(\ub_{ij},\e_{j+1})\chi( \e_{j+1},\ub_{ij})q_{j+1,j+1}
    \\ =& \ q_{\ub_{ij}}q_{j,j+1}q_{j+1,j}q_{j+1,j+1}
    \\ =& \left\{ \begin{array}{lc}  q_{\ub_{ij}} & q_{j+1,j+1} \neq -1 \  (\kappa_{i,j+1}=\kappa_{ij}),   \\  (-1)qq^{-1}=-1  & q_{j+1,j+1}=-1, \ \kappa_{ij} \mbox{ even},   \\   (-1)q(-1)=q  &   q_{j+1,j+1}=-1, \ \kappa_{ij} \mbox{ odd}.   \end{array} \right.
\end{align*}
So this proves the inductive step, and to calculate the dimension
of $\bB(V)$ we have to calculate the number of $\ub_{ij}$ such that
$$q_{\ub_{ij}}=q_{i,i+1}^{-1}q_{i+1,i}^{-1}=q^{\pm1},$$ this is,
$\sharp \left\{ \kappa_{ij}: i \leq j, \kappa_{ij} \mbox{ even}
\right\}$.

We consider an $1\times (\theta+1)$ board, numbered from 1 to
$\theta+1$, and paint its squares of white or black: the square
$\theta +1$ is white, and then the $i$-th square is the same color
of the $i+1$-th square if $q_{ii} \neq -1$, or different color if
$q_{ii}=-1$. All the possible colorations of this board are in
bijective correspondence with the choices of $1 \leq i_1< \ldots <
i_j \leq \theta$ for all $j$ (the positions where we put a $-1$ in
the corresponding $q_{ii}$ of the braiding), and the number of
white squares is
    $$ t= (\theta-i_j)+ (i_{j-1}-i_{j-2})+ \ldots = \theta - \sum ^{j}_{k=1}(-1)^{j-k}i_k $$
So $\sharp \left\{ \kappa_{ij}: i \leq j, \kappa_{ij} \mbox{ even}
\right\}$ is the number of pairs $(a,b), \quad 1 \leq a<b \leq
\theta+1$ ($a=i$ and $b=j+1$) such that the squares in positions
$a$ and $b$ are of the same color, that is, $$\binom{t}{2} +
\binom{\theta +1-t}{2}.$$ Then we obtain the formula
\eqref{dimensionAn} for the dimension of $\bB(V)$. \edem

\begin{prop}
Let $V$ be a standard braided vector space of type $B_{\theta}$ as
in Proposition \ref{casosBn}. Then the associated Nichols algebra
$\bB(V)$ is of finite dimension if and only if $q$ is a root of
unit of order $N \geq 2$. In such cases,
\begin{itemize}
    \item if the braiding is as in \emph{(a)} of Proposition \ref{casosBn},
    \begin{align}
    &\dim \bB(V)=3^2N^2 & \mbox{  if 3 divides }N, \\
    &\dim \bB(V)=3^3N^2 & \mbox{  if 3 does not divide }N;
    \end{align}
    \item if the braiding is as in \emph{(b)}, then $0\leq j\leq d-1$, and
    \begin{align}
    & \dim \bB(V)=2^{2t(\theta-t)+\theta}k^{\theta^2-2t\theta+2t^2} & \mbox{if }N=2k,
    \\ &\dim \bB(V)= 2^{(2t+1)(\theta-t)+1}N^{\theta^2-2t\theta+2t^2} & \mbox{if }N \mbox{ is
    odd;}
    \end{align}
    \item if the braiding is as in \emph{(c)}, then
    \begin{equation}
    \dim \bB(V)=2^{\theta(\theta-1)}3^{\theta^2-2t\theta+2t^2}.
    \end{equation}
\end{itemize}
Here, $ t= \theta - \sum ^{j}_{k=1}(-1)^{j-k}i_k $.
\end{prop}

\bdem It is clear that $q$ should be a root of 1.

Now, we proceed to calculate $\dim \bB(V)$. From Corollary
\ref{corollary:genPBW}, we have to determined $q_{\alpha}$ for
$\alpha \in \de_C$, and multiply their orders. As before,
$\ub_{ij}= \sum^{j}_{k=i} \e_k, \quad 1 \leq i \leq j \leq \theta$
and $\vb_{ij}= 2\sum^{i}_{k=1} e_k +\sum^{j}_{k=i+1} e_k=
2e_{1,i}+e_{i+1, j}, \quad 1 \leq i < j$, so
$$\de(\bB(V))= \{ \ub_{ij}: 1\leq i \leq j \leq \theta \} \cup \{ \vb_{ij}: 1\leq i < j \leq \theta \}.$$
We calculate $q_{\ub_{ij}}, \ 1 < i \leq j \leq \theta$ as above,
because they correspond with a braiding of standard $A_{\theta-1}$
type, and
\begin{eqnarray*}
    q_{\vb_{ij}} &=& \chi (\vb_{ij}, \vb_{ij})= \chi(\ub_{1i},\ub_{1i})^4 \chi(\ub_{1i}, \ub_{i+1,j}) ^2\chi(\ub_{i+1,j},\ub_{1i})^2 q_{\ub_{i+1,j}}
    \\ &=& q_{11}^4q_{12}^2q_{21}^2 \left( \prod^{i}_{k=2} q_{kk}^2q_{k-1,k} q_{k-1,k} q_{k+1,k} q_{k+1,k} \right)^2 q_{\ub_{i+1,j}} = q_{\ub_{i+1,j}},
\end{eqnarray*}
where we use that
\begin{itemize}
    \item $q_{ij}q_{ji}=1$ if $\left|i-j\right| >1$,
    \item $q_{11}^4q_{12}^2q_{21}^2=1$, and
    \item $q_{kk}^2q_{k-1,k}q_{k-1,k}q_{k+1,k}q_{k+1,k}=1$ if $2 \leq k \leq \theta-1$.
\end{itemize}
To calculate the other $q_{\alpha}$'s, we analyze each case:
\smallskip

(a) Note that $ q_{\e_1}=\zeta, \quad q_{\e_1+\e_2}=\zeta, \quad
q_{2\e_1+\e_2}=\zeta q^{-1}, \quad q_{\e_2}=q,$ so there are two
possibilities: $\dim \bB(V)=3^2N^2$ if 3 divides $N$, and $\dim
\bB(V)=3^3N^2$ if 3 does not divide $N$.
\medskip

(b) We have that
\begin{eqnarray*}
    q_{\ub_{1k}} &=& q^{-1} q_{\ub_{2k}}=\left\{ \begin{array}{lc}  q^2q^{-1}=q & \kappa_{2k} \mbox{ even},   \\   -q^{-1}  &  \kappa_{2k} \mbox{ odd};   \end{array} \right.
\end{eqnarray*}
and also $q_{11}=q$. We have that $\kappa_{2k}$ is even iff $j \in
\left\{i_j+1 , \theta \right\}$, or $i \in \left\{i_{j-2}+1,
i_{j-1}\right\}$, and so on. Then there are
    $$ t= (\theta-i_j)+ (i_{j-1}-i_{j-2})+ \ldots = \theta- \sum ^{j}_{k=1}(-1)^{j-k}i_k $$
numbers (the corresponding with the number in the above
Proposition) such that $\kappa_{i,\theta-1}$ is even. There are $2
\left( \binom{t}{2} + \binom{\theta -t}{2} \right)$ roots such
that $q_{\alpha}=q^{2}$, $2 \left( \binom{\theta}{2}-\binom{t}{2}
- \binom{\theta -t}{2} \right)$ roots such that $q_{\alpha}=-1$,
$t+1$ roots such that $q_{\alpha}=q$ and $\theta-1-t$ roots such
that $q_{\alpha}=-q^{-1}$.

Note that if $N=2k$, then $\ord(-q^{-1})=2k$ and $\ord(q^2)=k$, so
\begin{eqnarray*}
    \dim \bB(V)&=& 2^{(\theta -1)\theta-t(t-1)- (\theta -t)(\theta-t-1)}k^{t(t-1)+(\theta -t) (\theta-t-1) } (2k)^{\theta}
    \\ &=&  2^{2t(\theta-t)+\theta}k^{\theta^2-2t\theta+2t^2};
\end{eqnarray*}
if $N$ is odd, then $\ord(-q^{-1})=2N$ and $\ord(q^2)=N$, so
\begin{eqnarray*}
    \dim \bB(V)&=& 2^{\theta(\theta -1)-t(t-1)- (\theta-t)(\theta -1-t)}N^{t(t-1)+ (\theta-t) (\theta -1-t)+t+1}
    \\ &&  (2N)^{\theta-1-t} = 2^{(2t+1)(\theta-t)+1}N^{\theta^2-2t\theta+2t^2}.
\end{eqnarray*}

(c) In a similar way,
\begin{eqnarray*}
    q_{\ub_{1i}} &=& (-\zeta^2) q_{\ub_{2i}}=\left\{ \begin{array}{lc} (-\zeta^2)^2=\zeta & \kappa_{2i} \mbox{ even},   \\   (-1)(-\zeta^2)=\zeta^2  &  \kappa_{2i} \mbox{ odd};   \end{array} \right.
\end{eqnarray*}
and also $q_{11}=\zeta$. There are $2 \left( \binom{t}{2} +
\binom{\theta -t}{2}\right)$ roots such that
$q_{\alpha}=-\zeta^{2}$, $2\left( \binom{\theta}{2} -\binom{t}{2}
- \binom{\theta -t}{2} \right)$ roots such that $q_{\alpha}=-1$,
$t+1$ roots such that $q_{\alpha}=\zeta$ and $\theta-1-t$ roots
such that $q_{\alpha}=\zeta^2$. As $\ord \zeta= \ord \zeta^2=3$
and $\ord (-\zeta^2)=6$, we have
\begin{eqnarray*}
    \dim \bB(V)&=& 2^{\theta(\theta -1)-t(t-1)- (\theta-t)(\theta -1-t)}6^{t(t-1)+ (\theta-t) (\theta -1-t)} 3^{\theta}
    \\ &=& 2^{\theta(\theta-1)}3^{\theta^2-2t\theta+2t^2}.
\end{eqnarray*}
So, the proof is completed.
\edem

\begin{prop}\label{dimensionG2}
Let $V$ be a standard braided vector space of type $G_2$ as in
Proposition \ref{casosG2}. Then the associated Nichols algebra
$\bB(V)$ is of finite dimension if and only if $q$ is a root of
unit of order $N \geq 4$. Then
\begin{itemize}
    \item in case \emph{(a)} of Proposition \ref{casosG2},
    \begin{align*}
    & \dim \bB(V)=N^6 & \mbox{ if 3 does not divide }N,
    \\ & \dim \bB(V)=27k^6 & \mbox{ if }N=3k;
    \end{align*}
    \item in case \emph{(b)}, $\dim \bB(V)= 2^{12}$.
\end{itemize}
\end{prop}

\bdem For (a) note that $q$ is a root of 1, because $x_1$ has
finite height, and
\begin{itemize}
    \item $q_{\alpha}=q$ if $\alpha \in \left\{ e_1, e_1+e_2, 2e_1+e_2 \right\}$,
    \item $q_{\alpha}=q^3$ if $\alpha \in \left\{ e_2, 3e_1+e_2, 3e_1+2e_2 \right\}$,
\end{itemize}
so the dimension is $\dim \bB(V)=N^6$ if 3 does not divide $N$,
and $\dim \bB(V)=27k^6$ if $N=3k$. For (b) we calculate
\medskip

\begin{tabular}{l|l|l|l|l|l|l|l}
 type & $q_{x_2}$ & $q_{x_1x_2}$  & $q_{x_1^3x_2^2}$ & $q_{x_1^2x_2}$  & $q_{x_1^3x_2}$ & $q_{x_1}$ & $\dim \bB(V)$   \\
\hline \hline
       \Dchaintwo{}{$\zeta ^2$}{$\zeta $}{$\zeta ^{-1}$} & 8 & 4 & 2 & 8 & 2 & 4 & $2^{12}$ \\
       \Dchaintwo{}{$\zeta ^2$}{$\zeta ^3$}{$-1$} & 2 & 8 & 2 & 4 & 8 & 4 & $2^{12}$ \\
       \Dchaintwo{}{$\zeta $}{$\zeta^5$}{$-1$} & 2 & 4 & 8 & 4 & 2 & 8 & $2^{12}$
       \\
\end{tabular}

so the proof is complete.

\edem

\section{Presentation by generators and relations of Nichols algebras of
standard braided vector spaces}\label{section:genrel} \

In this section we give a presentation by generators and relations
of Nichols algebras of standard braided vector spaces. To do this,
 we give some technical results about relations and PBW-bases in Subsection \ref{subsection:relations};
also we calculate the coproduct of some hyperwords in $T(V)$. In
Subsections \ref{subsection:An}, \ref{subsection:Bn} and
\ref{subsection:presentationG2} we express the braided commutator
of two PBW-generators as combination of elements of the PBW-basis
under some assumptions. Then, we obtain the desired presentation
with a proof similar to the ones in \cite{AD} and \cite{AS5}. In
Subsection \ref{subsection:presentation} we solve the problem when
the braiding is of Cartan type using the transformation in
Subsection \ref{subsection:transformation}.

There is a procedure to describe a (non-minimal) set of relations
for Nichols algebras of rank 2 in \cite[Th. 4]{H4}.

\subsection{Some general relations}\label{subsection:relations} \

Let $V$ be a standard braided vector space with connected Dynkin
diagram. Let $x_1, \ldots, x_{\theta}$ be an ordered basis of $V$,
and $\left\{ x_{\alpha}: \alpha \in \de^+(\bB(V)) \right\}$ a set
of PBW generators. Here, $x_{\alpha} \in \bB(V)$ is, by abuse of
notation, the image by the canonical projection of $x_{\alpha} \in
T(V)$, the hyperword corresponding to a Lyndon word $l_{\alpha}$.
We denote
$$ q_{\alpha}:= \chi(\alpha, \alpha), \qquad N_{\alpha}:= \ord q_{\alpha}, \quad \alpha \in \de^+(\bB(V)). $$

Note that each $x_{\alpha}$ is homogeneous and has the same degree as $l_{\alpha}$. Also,
\begin{equation}\label{bihomog}
x_{\alpha}\in T(V)^{\chi_{\alpha}}_{g_{\alpha}},
\end{equation}
where if $\alpha = b_1\e_1 + \dots + b_{\theta}\e_{\theta}$, then
$g_{\alpha} = g_{1}^{b_1} \dots g_{\theta}^{b_{\theta}}$,
$\chi_{\alpha} = \chi_{1}^{b_1} \dots \chi_{\theta}^{b_{\theta}}$.

\begin{prop}\label{Prop:formabilineal}
If the matrix of the braiding is symmetric, then the PBW basis is
orthogonal with respect to the bilinear form in Proposition
\ref{formabilineal}.
\end{prop}

\bdem We prove by induction on $k:= \max \{ \ell (u), \ell (v) \}$
that $(u | v)=0$, where $u \neq v$ are products of PBW generators
(we also allow powers greater than the corresponding heights). If
$k=1$, then $u=1$ or $x_i$, $v=x_j$, for some $i,j \in \unon$, and
$(x_i | x_j)= \delta_{ij}$.

Suppose it is valid when the length of both words is least than
$k$, and let $u, v \in B_{I(V)}, u \neq v$ be hyperwords such that
one (or both) has length $k$. If both are hyperletters, they have
different degrees $\alpha \neq  \beta \in \zt$, so $u=x_{\alpha}$,
$v=x_{\beta}$, and $(x_{\alpha} | x_{\beta})=0$, since the
homogeneous components are orthogonal for $( | )$.

Suppose that $u=x_{\alpha}$ and $v=x_{\beta_1}^{h_1} \ldots
x_{\beta_m}^{h_m}$, for some $x_{\beta_1} > \ldots > x_{\beta_m}$.
If they have different $\zt$-degree, they are orthogonal. Then, we
assume that $\alpha= \sum_{j=1}^{m}h_m \beta_m$. By \cite[Ch. VI,
Prop. 19]{B}, we can reorder the $\beta_i$'s, using $h_i$ copies
of $\beta_i$, in such form that each partial sum is a root. Using
\cite[Prop. 21]{R2}, the order induced by the Lyndon words
$l_{\alpha}$ is convex, so $\beta_n < \alpha$. Using Lemma
\ref{copro} and \eqref{bilinearprop2},
\begin{align*}
(u | v) =& ( x_{\alpha} | w) (1 | x_{\beta_m}) + (1 |w)
(x_{\alpha_n} | x_{\beta_m}) \\ &+ \sum_{ l_1\geq \dots \geq l_p
>\alpha, l_i \in
        L}  ( x_{l_1,\dots ,l_p} | w) ( [l_1]_c \cdots [l_p]_c | x_{\beta_m} )
\end{align*}
where $v=wx_{\beta_n}$. Note that $(1 | x_{\beta_m})= (1 |w) =0$.
Also, $[l_1]_c \cdots [l_p]_c$ is a linear combination of greater
hyperwords of the same degree and an element of $I(V)$. By
inductive hypothesis and the fact that $I(V)$ is the radical of
the bilinear form, $( [l_1]_c \cdots [l_p]_c | x_{\beta_m} ) =0$.

Consider now
$$u=x_{\alpha_1}^{j_1} \ldots x_{\alpha_n}^{j_n}, \ x_{\alpha_1} >  \ldots > x_{\alpha_n}, \quad v=x_{\beta_1}^{h_1} \ldots x_{\beta_m}^{h_m}, \ x_{\beta_1} > \ldots > x_{\beta_m}, $$
and suppose that $x_{\alpha_n} \leq x_{\beta_m}$ (if not, use that
the bilinear form is symmetric). Using Lemma \ref{copro} and
\eqref{bilinearprop2},
\begin{align*}
    (u | v) &= ( w | 1) (x_{\alpha_n} | v) + \sum ^{h_m}_{i=0} \binom{ h_m }{ i } _{q_{\beta_m}} (w | x_{\beta_1}^{h_1} \ldots x_{\beta_{m-1}}^{h_{m-1}}x_{\beta_m}^{i} ) (x_{\alpha_n} | x_{\beta_m}^{h_m-i})
        \\ & \quad + \sum_{ \substack{ l_1\geq \dots  \geq l_p >l, l_i \in L \\ 0\leq j \leq m } } (w | x_{l_1,\dots ,l_p}^{(j)}) (x_{\alpha_n} | \left[l_1\right]_c \dots      \left[l_p\right]_c \left[x_{\beta_m}\right]_c^j )
\end{align*}
where $w= x_{\alpha_1}^{h_1} \ldots x_{\alpha_m}^{h_m-1}$. Note
that for the first summand, $( w | 1)=0$. In the last sum,
$(x_{\alpha_n} | \left[l_1\right]_c \dots      \left[l_p\right]_c
\left[x_{\beta_m}\right]_c^j )=0$, because by the previous
results, $\left[l_1\right]_c \dots \left[l_p\right]_c
\left[x_{\beta_m}\right]_c^j$ is a combination of hyperwords of
the PBW basis greater or equal than it and an element of $I(V)$,
then we use induction hypothesis and the fact that $I(V)$ is the
radical of this bilinear form. As also $x_{\alpha_n},
x_{\beta_m}^{h_m-i} $ are different elements of the PBW basis for
$h_m-i \neq 1$, we have that
    \[ (u | v)= (h_m )_{q_{\beta_m}} (w | x_{\beta_1}^{h_1} \ldots x_{\beta_{m-1}}^{h_{m-1}} x_{\beta_m}^{h_m-1} ) (x_{\alpha_n} | x_{\beta_m}) . \]
Then it is zero if $\alpha_n \neq \beta_m$, but also if $\alpha_n=
\beta_m$, because in that case $w$, $x_{\beta_1}^{h_1} \ldots
x_{\beta_{m-1}}^{h_{m-1}} x_{\beta_m}^{h_m-1}$ are different
products of PBW generators, and we use induction hypothesis. \edem

\begin{cor}\label{heigthgenerators}
If $\alpha \in \de^+(\bB(V))$, then
\begin{equation}\label{powerrootvector}
x_{\alpha}^{N_{\alpha}}=0.
\end{equation}
\end{cor}
\bdem Let $(q_{ij})$ be symmetric. If $u=x_{\alpha_1}^{j_1} \ldots
x_{\alpha_n}^{j_n}, \quad x_{\alpha_1} >  \ldots > x_{\alpha_n}$,
then
\begin{equation} \label{bilformword}
(u | u )  =  \prod ^{n}_{i=1} (j_i)_{q_{\alpha_i}} ! (x_{\alpha_i} | x_{\alpha_i})^{j_i} ,
\end{equation}
where $(x_{\alpha} | x_{\alpha} ) \neq 0$ for all $\alpha \in
\de^+(\bB(V))$.

If we consider $u=x_{\alpha}^{N_{\alpha}}$, we have that $(u | v
)=0$, for each element $v$ of the PBW basis, because they are
ordered products of $x_{\alpha}$ different of $u$, and $(u|u)=0$
since $q_{\alpha} \in \G_{N_{\alpha}}$. Also, $(I(V) |
x_{\alpha}^{N_{\alpha}})=0$, because it is the radical of this
bilinear form, so $(T(V) | x_{\alpha}^{N_{\alpha}})=0$, and then
$x_{\alpha}^{N_{\alpha}} \in I(V)$. That is, we have
$x_{\alpha}^{N_{\alpha}}=0$  in $\bB(V)$.
\medskip

For the general case, we recall that a diagonal braiding is twist
equivalent to a braiding with a symmetric matrix, see
\cite[Theorem 4.5]{AS3}. Also, there exists a linear isomorphism
between the corresponding Nichols algebras. The corresponding
$x_{\alpha}$ are related by a non-zero scalar, because they are
an iteration of braided commutators between the hyperwords.
\edem
\medskip

We shall need some technical results about graded algebras
intermediated between $T(V)$ and $\bB(V)$.

\begin{lema}\label{buscarmij}
Let $i, \neq j \in \unon$. Let $\bB$ be a graded algebra provided
with an inclusion of braided vector spaces $V \hookrightarrow
\bB^1$. Assume that:
\begin{itemize}
    \item there exist skew derivations $D_i$ of $\bB$ as in Proposition \ref{derivations};
    \item $x_i^N=0$ if $N:= \ord q_{ii}< \min \{ n \in \N: q_{ii}^nq_{ij}q_{ji}=1 \}-1$.
\end{itemize}
For each $m \in \N$, $x_i^mx_j$ is a linear combination of greater hyperwords (for a fixed order such that $x_i<x_j$) if and only if
\begin{equation}\label{qserre}
    (\ad x_i)^{m_{ij}+1}x_j=0, \quad i \neq j.
\end{equation}
\end{lema}
\bdem
If $(ad_c x_i)^{m}x_j=0$, there exist $a_k \in \kk$ such that
$$0= \left[ x_i^mx_j \right]_c= (\ad_cx_i)^m x_j=x_i^mx_j +\sum ^{m-1}_{k=0}
a_k x_i^k x_jx_i^{m-k}. $$

Conversely, suppose that there exist $m \in \N$ such that $x_i^mx_j$ is a linear combination of greater hyperwords. Let $$n= \min \{ m \in \N: x_i^mx_j \mbox{ is a linear combination of greater hyperwords} \}.$$

If $x_i^n=0$, then $q_{ii}$ is a root of 1, because of the derivations. In this case, if $N$ is the order of $q_{ii}$, then $x_i^N=0$ and $x_i^{N-1} \neq 0$. Also, $(\ad_c x_i)^Nx_j=0$. Hence, we can assume $x_i^n \neq 0$ and $(n)_{q_{ii}}! \neq 0$.

Note that $ [x_i^{n-k}x_jx_i^k ]_c = [x_i^{n-k}x_j]_c x_i^k$. As $\bB$ is graded, $x_i^nx_j$ is a linear combination of $x_i^{n-k}x_jx_i^k$, $0\leq k <n$   Hence, there exist $\alpha_k \in \kk$ such that
    \[\left[ x_i^n x_j\right]_c= \sum^{n}_{k=1} \alpha_k
    \left[x_i^{n-k}x_j\right]_cx_i^k.\]
Applying $D_i$ we obtain
$$  0 = D_i(\left[ x_i^n x_j\right]_c)= \sum^{n}_{k=1} \alpha_k D_i\left( \left[x_i^{n-k}x_j\right]_cx_i^k
  \right) = \sum^{n}_{k=1} \alpha_k (k)_{q_{ii}} \left[x_i^{n-k}x_j\right]_cx_i^k. $$
By the hypothesis about $n$, $\alpha_1=0$. As $(n)_{q_{ii}}! \neq 0$, applying $D_i$ several times we conclude that $\alpha_k=0$ for $k=2,\ldots,n$. Then, $\left[ x_i^n x_j\right]_c=0$.
\edem

Recall that \ref{qserre} holds in $\bB(V)$, for $1 \leq i \neq j \leq \theta$

\bigbreak

The second lemma is related to Dynkin diagrams of a standard
braiding which have two consecutive simple edges.

\begin{lema}\label{lemaAn}
Let $\bB$ be a graded algebra provided with an inclusion of
braided vector spaces $V \hookrightarrow \bB^1$. Assume that:
\begin{itemize}
    \item there exist skew derivations $D_i$ in $\bB$ as in Proposition \ref{derivations};
    \item there exist different $j,k,l \in \{ 1, \ldots, \theta \}$ such that $m_{kj}=m_{kl}=1$, $m_{jl}=0$;
    \item $(\ad x_k)^2x_j=(\ad x_k)^2x_l= (\ad x_j) x_l=0$ hold in $\bB$;
    \item $x_k^2=0$ if $q_{kk}q_{kj}q_{jk} \neq 1$ or $q_{kk}q_{kl}q_{lk} \neq 1$.
\end{itemize}

(1) If we order the letters $x_1, \ldots, x_{\theta}$ such that $x_j< x_k < x_l$, then $x_jx_kx_lx_k$ is a linear combination of greater words if and only if
\begin{equation}\label{relA}
    \left[ (\ad x_j)(\ad x_k)x_l, x_k \right]_c=0.
\end{equation}

(2) If $V$ is standard and $q_{kk} \neq -1$, then \eqref{relA} holds
     in $\bB$.
\smallskip

(3) If $V$ is standard and $\dim \bB(V)< \infty$, then \eqref{relA} holds in $\bB=\bB(V)$.
\end{lema}

\bdem (1) ($\Leftarrow$) If \eqref{relA} holds, then $x_jx_kx_lx_k$ is a linear combination of greater words, by Remark \ref{corchete}, and
$$ \left[ x_jx_kx_lx_k \right]_c = \left[ \left[ x_jx_kx_l \right]_c , x_k \right]_c = \left[ (\ad x_j) (\ad x_k) x_l, x_k \right]_c . $$

($\Rightarrow$) If $x_jx_kx_lx_k$ is a linear combination of greater words, then the hyperword $ \left[ x_jx_kx_lx_k \right]_c$ is a linear combination of hyperwords corresponding to words greater than $x_jx_kx_lx_k$ (of the same degree, because $\bB$ is homogeneous); this follows by Remark \ref{corchete}. As $(\ad x_k)^2x_j=(\ad x_k)^2x_l=(\ad x_j)x_l=0$, we do not consider hyperwords with $x_jx_k^2$, $x_k^2x_l$ and $x_jx_l$ as factors of the corresponding words. Then, $\left[ x_jx_kx_lx_k \right]_c$ is a linear combination of
\begin{align*}
&\left[ x_kx_lx_k x_j \right]_c = \left[ x_kx_l \right]_c x_kx_j,
 &\left[ x_lx_kx_jx_k \right]_c  = x_l x_k \left[ x_jx_k \right]_c,
\\ &\left[ x_kx_jx_kx_l \right]_c = x_k \left[ x_jx_kx_l \right]_c,
 &\left[ x_lx_k^2x_j \right]_c = x_l x_k^2x_j.
\end{align*}
As $D_j(\left[ x_jx_kx_lx_k \right]_c)= D_j (x_k \left[ x_jx_kx_l \right]_c)= D_j (x_l x_k \left[ x_jx_k \right]_c)=0$, in that linear
combination there are no hyperwords ending in $x_j$; indeed,
$$D_j (\left[ x_kx_l \right]_c x_kx_j) = \left[ x_kx_l \right]_c x_k, \quad D_j(x_l x_k^2x_j) = x_l x_k^2,$$
and $\left[ x_kx_l \right]_c x_k$, $x_l x_k^2$ are linearly independent. Therefore, there
exist $\alpha, \beta \in \kk$ such that
$$ \left[ x_j x_k x_l x_k \right]_c = \alpha x_l x_k \left[ x_j x_k \right]_c + \beta x_k \left[ x_jx_k x_l \right]_c. $$
Applying $D_l$, we have
\begin{align*}
    0 =& \alpha q_{kj}q_{kk}x_l \left[ x_jx_k \right]_c + \alpha (1-q_{kj}q_{jk})x_lx_kx_j + \beta q_{kk}q_{kj}q_{kl} \left[ x_jx_kx_l \right]_c .
\end{align*}
Now, $x_l \left[ x_jx_k \right]_c$, $x_lx_kx_j$ and $\left[
x_jx_kx_l \right]_c$ are linearly independent by Lemma
\ref{lemaAHS}, so $\alpha=\beta=0$.
\smallskip

\medskip

(2) We assume that some quantum Serre relations hold in $\bB$; using them:
\begin{eqnarray*}
    x_jx_kx_lx_k &=& q_{kl}^{-1}(1+q_{kk})^{-1}x_jx_k^2x_l+ q_{kk}q_{kj} (1+q_{kk})^{-1}x_jx_lx_k^2
    \\ &=& q_{kk}^{-1}q_{kj}^{-1}  q_{kl}^{-1}x_kx_jx_kx_l +q_{kk}^{-1}q_{kj}^{-1}q_{kl}^{-1}(1+q_{kk})^{-1}x_k^2x_jx_l
    \\ && +q_{kk}q_{kl}q_{jk} (1+q_{kk})^{-1}x_lx_jx_k^2.
\end{eqnarray*}
It follows that $x_kx_jx_kx_l \notin G_I$, for an order such that $x_j< x_k < x_l$. Also, $x_jx_lx_k^2 \notin G_I$, since $(ad_c x_j)x_l=0$, and \eqref{relA} is valid by
the previous item.

\medskip

(3) If $V$ is a standard braided vector space satisfying the above
conditions, then we consider $V_k$ as the braided vector space
obtained transforming by $s_k$, then $\widetilde{m}_{jl}=0$.
Therefore, $\e_j+\e_l \notin \de^{+} (\bB(V_k))$, so
$s_k(\e_j+\e_l)=2\e_k+\e_j+\e_l \notin \de^{+} (\bB(V))$. It
follows that $x_jx_kx_lx_k$ is a linear combination of greater
words, since it is a Lyndon word when we consider an order such
that $x_j<x_k<x_l$. \edem

\smallskip

We prove now two relations related to the double edge in a Dynkin diagram of standard braiding of type $B_{\theta}$.

\begin{lema}\label{lemaBn}
Let $\bB$ be a graded algebra provided with an inclusion of
braided vector spaces $V \hookrightarrow \bB^1$. Assume that:
\begin{itemize}
    \item there exist $j \neq k \in \{ 1, \ldots, \theta \}$ such that $m_{kj}=2, m_{jk}=1$;
    \item there exist skew derivations as in Proposition \ref{derivations};
    \item the following relations hold in $\bB$:
    \begin{align} \label{Lema1B}
    &(\ad x_k)^3x_j=(\ad x_j)^2x_k=0; &
    \\ \nonumber & x_k^3=x_j^2=0, & \mbox{if }q_{kk}^3=q_{jj}^2=1.
\end{align}
\end{itemize}

(1) If we order the letters $x_1, \ldots, x_{\theta}$ such that $x_k< x_j$, then $x_k^2x_jx_kx_j$ is a linear combination of greater words if and only if
\begin{equation}\label{relB}
    \left[ (\ad x_k)^2x_j, (\ad x_k)x_j \right]_c=0.
\end{equation}

(2) If $V$ is standard, $\ q_{jj} \neq -1$ and $q_{kk}^{2}q_{kj}q_{jk}=1$, then
\eqref{relB} holds in $\bB$.

(3) If $V$ is standard and $\dim \bB(V)< \infty$, then \eqref{relB} holds in $\bB=\bB(V)$.
\end{lema}

\bdem (1) ($\Leftarrow$) If \eqref{relB} holds in $\bB$, then $x_k^2x_jx_kx_j$ is a linear combination of greater words. It follows from \eqref{corchete}, and
$$ \left[x_k^2x_jx_kx_j\right]_c = \left[ \left[x_k^2x_j \right]_c , \left[ x_kx_j \right]_c \right]_c = \left[ (\ad x_k)^2x_j, (\ad x_k)x_j \right]_c.  $$

($\Rightarrow$) If $x_k^2x_jx_kx_j$ is a linear combination of greater words, then $ [ x_k^2x_jx_kx_j ]_c$ is a linear combination of hyperwords corresponding to words greater
than $x_k^2x_jx_kx_j$ (of the same degree, because $\bB$ is homogeneous).

First, there are not hyperwords whose corresponding
words have factors $x_k^3x_j$, $x_kx_j^2$, by \ref{Lema1B}. As $ \left[
x_k^2x_jx_kx_j \right]_c \in \ker D_k$, and
\begin{align*}
D_k( x_j [x_k^{2}x_j]_c x_k) &= x_j [x_k^{2}x_j]_c ,
\\ D_k ([x_kx_j]_c ^2 x_k) &= [x_kx_j]_c ^2,
\\ D_k (x_j [x_kx_j]_c x_k^2) &= (1+q_{kk}) x_j [x_kx_j]_c x_k ,
\end{align*}
in that linear combination there are no hyperwords ending in
$x_k$, except $x_j^2x_k^3$ if $q_{kk} \in \G_3$. We consider
$q_{jj} \neq -1$ if $q_{kk} \in \G_3$, since otherwise
$x_j^2x_k^3=0$ by hypothesis. Then, there exists $\alpha, \alpha'
\in \kk$ such that
$$ \left[ x_k^2x_jx_kx_j \right]_c = \alpha \left[ x_kx_jx_k^2x_j \right]_c + \alpha'x_j^2x_k^3 = \alpha \left[ x_kx_j\right]_c \left[ x_k^2x_j \right]_c + \alpha'x_j^2x_k^3 . $$
We prove by direct calculation that $D_j(\left[ x_k^2x_jx_kx_j
\right]_c)=0$. Then, applying $D_j$ to the previous equality,
\begin{align*}
    0=& \alpha' (1+q_{jj}) x_jx_k^3 + \chi(\e_k+\e_j , 2\e_k+\e_j) \alpha (\ad x_k)^2(x_j)x_k
    \\ &+(1-q_{kj}q_{jk}) (1-q_{kk}q_{kj}q_{jk})\alpha (\ad x_k)(x_j)x_k^2,
\end{align*}
where we use that $(\ad x_k)^3(x_j)=0$ and
$$x_k(\ad x_k)^m(x_j)= (\ad x_k)^{m+1}(x_j)+ q_{kk}^mq_{kj}(\ad x_k)^m(x_j)x_k.$$
As $(1-q_{kj}q_{jk}) (1-q_{kk}q_{kj}q_{jk}) \neq 0$ and $(\ad x_k)^2(x_j)x_k$, $(\ad x_k)(x_j)x_k^2$, $x_jx_k^3$ are linearly
independent, it follows that $\alpha= \alpha'=0$.

\medskip

(2) Using $(\ad x_j)^2 x_k=0$ in the first equality and $(\ad
x_k)^3 x_j=0$ in the last expression,
\begin{align*}
 x_k^2x_jx_kx_j & = (1+q_{jj})^{-1}q_{jk}^{-1} x_k^2x_j^2x_k +
 (1+q_{jj})^{-1}q_{jk}q_{jj}x_k^3x_j^2
 \\ & \in
 (3)_{q_{kk}}(1+q_{jj})^{-1}q_{kj}q_{jk}q_{jj} x_k^2x_jx_kx_j +
 \kk \xx _{>x_k^2x_jx_kx_j}.
\end{align*}
Suppose that $(3)_{q_{kk}}(1+q_{jj})^{-1}q_{kj}q_{jk}q_{jj}=1$;
that is, $(3)_{q_{kk}} = (1+q_{jj})$. Then,
$q_{jj}=q_{kk}+q_{kk}^{2}$, so
$$ 1= q_{jj}q_{kj}q_{jk} =q_{kk}q_{kj}q_{jk} +q_{kk}^{2}q_{kj}q_{jk}=
q_{kk}q_{kj}q_{jk}+1,$$  which is a contradiction since
$q_{kk}q_{kj}q_{jk} \in \kk^{\times}$. It follows that
$x_k^2x_jx_kx_j$ is a linear combination of greater words, so
\eqref{relB} follows by previous item.

\medskip

(3) If $V$ is a standard braided vector space, and we consider
$V_j$ as the braided vector space obtained transforming by $s_j$,
then $\widetilde{m}_{kj}=2$. Therefore, $3\e_k+\e_j \notin \de^{+}
( \bB(V_k))$, so $s_j(3\e_k+\e_j)=3\e_k+2\e_j \notin \de^{+} (
\bB(V))$. As $x_k^2x_jx_kx_j$ is a Lyndon word of degree
$3\e_k+2\e_j$ if $x_k<x_j$, then it is a linear combination of
greater words. \edem

\medskip

\begin{lema}\label{lemaBn2}
Let $\bB$ be a graded algebra provided with an inclusion of
braided vector spaces $V \hookrightarrow \bB^{1}$. Assume that
\begin{itemize}
    \item there exist different $j,k,l \in \{ 1, \ldots, \theta
\}$ such that $m_{kj}=2$, $m_{jk}=m_{jl}=m_{lj}=1$, $m_{kl}=0$;
    \item there exist skew derivations $D_i$ in $\bB$ as in Proposition \ref{derivations};
    \item the following relations hold in $\bB$: \eqref{relA}, \eqref{relB},
    \begin{align} \label{Lema2B}
    &(\ad x_k)^3x_j=(\ad x_j)^2x_k(\ad x_j)^2x_l=(\ad x_k) x_l=0; & \nonumber
    \\  & x_k^3=x_j^2=0, & \mbox{if }q_{kk}^3=q_{jj}^2=1.
\end{align}
\end{itemize}

(1) If we order the letters $x_1, \ldots, x_{\theta}$ such that $x_k< x_j<x_l$, then $x_k^2x_jx_lx_kx_j$ is a linear combination of greater words if and only if
\begin{equation} \label{relB2}
    \left[ (\ad x_k)^2(\ad x_j)x_l, (\ad x_k)x_j \right]_c=0.
\end{equation}

(2) If $V$ is a standard braided vector space and $q_{kk} \notin
\G_3$, $q_{jj} \neq -1$, then \eqref{relB2} holds in $\bB$.

(3) If $V$ is standard and $\dim \bB(V) < \infty$, then \eqref{relB2} holds in $\bB(V)$.
\end{lema}

\bdem (1) ($\Leftarrow$) As in last two Lemmata, if \eqref{relB2} is valid,
then $x_k^2x_jx_lx_kx_j$ is a linear combination of greater words,
by \eqref{corchete}, and
$$ [ x_k^2x_jx_lx_kx_j ]_c  = \left[ (\ad x_k)^2(\ad x_j)x_l, (\ad x_k)x_j \right]_c. $$

($\Rightarrow$) Suppose that $x_k^2x_jx_lx_kx_j$ is a linear combination of greater words. Then,
$[ x_k^2x_jx_lx_kx_j ]_c$ is a linear combination of hyperwords corresponding to words
greater than $x_k^2x_jx_lx_kx_j$ (of the same degree, because
$\bB$ is homogeneous). We discard those words which have $x_kx_l$,
$x_k^3x_j$, $x_kx_j^2$, $x_j^2x_l$, $x_kx_jx_lx_j$ and
$x_k^2x_jx_kx_j$ by the hypothesis about $\bB$.

Also, as $D_k(\left[ x_k^2x_jx_lx_kx_j \right]_c)=0$, the coefficients of those hyperwords corresponding to words ending in $x_k$ are 0 as in Lemma \ref{lemaBn}, except $\left[x_jx_l\right]_c x_jx_k^3$, $x_lx_j^2x_k^3$, if $q_{kk} \in \G_3$. Then,
\begin{align*}
    \left[ x_k^2x_jx_lx_kx_j \right]_c =& \alpha \left[ x_kx_j \right]_c \left[ x_k^2x_jx_l\right]_c + \beta \left[ x_kx_jx_l \right]_c \left[ x_k^2x_j \right]_c
    \\  \quad & + \gamma x_l \left[ x_kx_j \right]_c \left[ x_k^2x_j\right]_c + \mu \left[x_jx_l\right]_c x_jx_k^3+ \nu x_lx_j^2x_k^3.
\end{align*}
By direct calculation, $$D_j(\left[ x_k^2x_jx_lx_kx_j \right]_c) =
D_j(\left[ x_k^2x_jx_l\right]_c)= D_j (\left[ x_kx_jx_l
\right]_c)=0,$$
so applying $D_j$ to the previous equality,
\begin{align*}
    0= & \alpha q_{jk}^2q_{jj}q_{jl} x_j \left[ x_k^2x_jx_l\right]_c + \beta (1-q_{kk}q_{kj}q_{jk}) (1-q_{kk}^2 q_{kj}q_{jk}) \left[ x_kx_jx_l \right]_c x_k^2
    \\  & + \gamma (1-q_{kk}q_{kj}q_{jk}) (1-q_{kk}^2 q_{kj}q_{jk}) x_l \left[ x_kx_j \right]_c x_k^2 + \gamma q_{jk}^2q_{jj} x_l x_j  \left[ x_k^2x_j\right]_c
    \\  & + \mu \left[x_jx_l\right]_c x_k^3+ \nu (1+q_{jj})x_lx_jx_k^3,
\end{align*}
Note that $\nu=0$ if $q_{jj} \neq -1$; otherwise, $x_j^2=0$ by
hypothesis, so we can discard this last summand. The other
hyperwords appearing in this expression are linearly independent,
since the corresponding words are linearly independent by Lemma
\ref{lemaAHS}. Then, $\alpha=\beta=\gamma=\mu=0$.

\medskip

(2) If $q_{kk} \notin \G_3$ and $q_{jj} \neq -1$, then
$x_k^2x_jx_lx_kx_j$ is a linear combination of greater words,
using the quantum Serre relations in a  similar way that in Lemma
\ref{lemaBn2}, so we apply the previous item.

\medskip

(3) If $V$ is a standard braided vector space, and we consider
$V_k$ as the braided vector space obtained transforming by $s_k$,
then $\widetilde{m}_{kj}=2$. Therefore, $\e_k+2\e_j+\e_l \notin
\de^{+} ( \bB(V_k))$ by Lemma \ref{lemaBn}, so
$s_k(\e_k+2\e_j+\e_l)=3\e_k+2\e_j+\e_l \notin \de^{+} ( \bB(V))$.
As $x_k^2x_jx_lx_kx_j$ is a Lyndon word, it follows that it is a
linear combination of greater words, and we apply (1). \edem

\medskip

We give now explicit formulas for the comultiplication of previous hyperwords.

\begin{lema}\label{lemacopro}
Consider the structure of graded braided Hopf algebra of $T(V)$,
given in subsection \ref{subsection:braidedHA}. Then, for all $k
\neq j$,
\begin{eqnarray}\label{coproductQSR}
    \Delta ( (\ad x_k)^{m_{kj}+1}x_j ) &=& (\ad x_k)^{m_{kj}+1}x_j
    \otimes 1 + 1 \otimes (\ad x_k)^{m_{kj}+1}x_j \nonumber
    \\ && + \prod_{1 \leq t \leq m_{kj}} (1-q_{kk}^{t}q_{kj}q_{jk}) x_k^{m_{ij}+1} \otimes
    x_j.
\end{eqnarray}
\end{lema}
\bdem By the definition of $m_{kj}$ and \eqref{paragenerar}, $F_k((\ad x_k)^{m_{kj}+1}x_j)=0$. Also, $F_l((\ad x_k)^{m_{kj}+1}x_j )$ for $l \neq k$ by \eqref{paragenerar2} and
the properties of $F_l$, so we have
$$\de_{1,m_{kj}}((\ad x_k)^{m_{kj}+1}x_j )= \sum^{\theta}_{l=1} x_l \otimes F_l ((\ad x_k)^{m_{kj}+1}x_j ) = 0. $$
Now, $D_k( \left[ x_k^ix_j \right]_c x_k^{s-i})=0$ from \eqref{14}, and from \eqref{15}
$$D_j(\left[ x_k^ix_j \right]_c x_k^{s-i})= \prod_{1 \leq t \leq m_{kj}}
(1-q_{kk}^{t}q_{kj}q_{jk}) x_k^{m_{ij}+1},$$
so we deduce that
$$\de_{m_{kj},1} ((\ad x_k)^{m_{kj}+1}x_j )= \prod_{1 \leq t \leq m_{kj}}
(1-q_{kk}^{t}q_{kj}q_{jk}) x_k^{m_{ij}+1} \ot x_j.$$
As hyperwords form a basis of $T(V)$, we can express for each $1<s<m_{kj}$,
\begin{align*}
\de_{m_{kj}+1-s,s}((\ad x_k)^{m_{kj}+1}x_j ) =& \sum^{m_{kj}+1-s}_{t=0} \epsilon_{st} \left[ x_k^t x_j \right]_c x_k^{m_{kj}+1-s-t} \otimes x_k^s
\\ &+ \sum^{s}_{p=0} \rho_{sp} x_k^{m_{kj}+1-s} \otimes \left[ x_k^{s-p} x_j \right]_c x_k^p,
\end{align*}
for some $\epsilon_{st},\rho_{sp} \in \kk$. Then, for each $0 \leq t \leq m_{kj}+1-s$,
\begin{align*}
    0=& ( (\ad x_k)^{m_{kj}+1}x_j | \left[ x_k^t x_j \right]_c x_k^{m_{kj}+1-t} x_k^s )
    \\ =& \left(( (\ad x_k)^{m_{kj}+1}x_j)_{(1)} | \left[ x_k^t x_j \right]_c x_k^{m_{kj}+1-t-s}\right) \left( ( (\ad x_k)^{m_{kj}+1}x_j)_{(2)} | x_k^s \right)
    \\ =& \epsilon_{st} \left( \left[ x_k^t x_j \right]_c x_k^{m_{kj}+1-t-s} | \left[ x_k^t x_j \right]_c x_k^{m_{kj}+1-t-s}\right) \left( x_k^s | x_k^s \right)
    \\ =& \epsilon_{st} (m_{kj}+1-s-t)_{q_{kk}}! (s)_{q_{kk}}! ( \left[ x_k^t x_j \right]_c | \left[ x_k^t x_j \right]_c ),
\end{align*}
where we use that $(\ad x_k)^{m_{kj}+1}x_j \in I(V)$ for the first equality, \eqref{bilinearprop2} for the second, \eqref{bilinearprop4} and the orthogonality between increasing products of hyperwords for the third, and \eqref{bilformword} for the last. As
$$ (m_{kj}+1-s-t)_{q_{kk}}! (s)_{q_{kk}}! ( \left[ x_k^t x_j \right]_c | \left[ x_k^t x_j \right]_c ) \neq 0,$$
we conclude that $\epsilon_{st}$ for all $0 \leq t \leq m_{kj}+1-s$. In a similar way, $\rho_{sp}=0$ for all $0 \leq p \leq s$, so we obtain \eqref{coproductQSR}.
\edem
\medskip

\begin{lema}\label{lemacoproA}
Let $\bB$ be a braided graded Hopf algebra provided with an
inclusion of braided vector spaces $V \hookrightarrow \cP(\bB)$.
Assume that
\begin{itemize}
    \item there exist $1 \leq j \neq k \neq l \leq \theta$ such that
$m_{kj}=m_{kl}=1$, $m_{jl}=0$;
    \item the following relations hold in $\bB$:
\begin{align*}
&(\ad x_k)^2x_j=(\ad x_k)^2x_l=(\ad x_j)x_l=0;
\\ &x_k^2=0 \mbox{ if } q_{kk}q_{kj}q_{jk} \neq 1 \mbox{ or } q_{kk}q_{kl}q_{lk} \neq 1.
\end{align*}
\end{itemize}
Then, $\quad u:=[(\ad x_j) (\ad x_k) x_l, x_k]_c \in \cP(\bB)$.
\end{lema}

\bdem From \eqref{14}, $D_j(u)=0$. Also, $D_k\left( (\ad x_j) (\ad
x_k) x_l \right)=0$, so
$$ D_k(u)=  \left( 1-q_{kk}^2q_{jk}q_{kj}q_{kl}q_{lk} \right)  (\ad x_j) (\ad x_k)
x_l =0. $$ From \eqref{15} and the properties of $D_l$ we have
\begin{align*}
    D_l(u) =& q_{lk} (1-q_{kl}q_{lk}) [x_jx_k]_c x_k - q_{jk}q_{kk}q_{lk} (1-q_{kl}q_{lk}) x_k [x_jx_k]_c
    \\ =& q_{lk} (1-q_{lk}q_{kl}) [ [x_jx_k]_c, x_k ]_c =0.
\end{align*}
Then, $\de_{31}(u)=0$. Now, from \eqref{paragenerar2} and the properties of $F_k, F_l$ we have $F_k(u)=F_l(u)=0$. Using \eqref{paragenerar}, we have
\begin{align*}
    F_j(u) =& (1-q_{jk}q_{kj}) [x_kx_l]_c x_k - q_{jk}q_{kk}q_{lk}q_{kj} (1-q_{jk}q_{kj}) x_k [x_kx_l]_c
    \\ =& (1-q_{lk}q_{kl}) (1-q_{kj}q_{jk}q_{kk}^2q_{lk}q_{jk}) [x_kx_l]_cx_k  =0.
\end{align*}
Then, we also have $\de_{13}(u)=0$.

Also, we have
$$\de(u)= \de ((\ad x_j) (\ad x_k) x_l) \de(x_k) - q_{\e_k+\e_j+\e_j,\e_j} \de (x_k) \de
((\ad x_j) (\ad x_k) x_l),$$ and looking at the terms in $\bB^2 \otimes \bB^2$,
\begin{align*}
    \de _{2,2}(u) =& (1-q_{lk}q_{kl}) [x_jx_k]_c \otimes \left( x_lx_k-q_{kj}q_{jk}q_{kk}^2q_{lk} x_kx_l \right)
    \\ & + (1-q_{kj}q_{jk}) q_{lk}q_{kk} \left( x_jx_k-q_{jk}x_kx_j
    \right) \otimes [x_kx_l]_c
    \\ =& \left( 1-q_{kj}q_{jk} - \left(1-q_{lk}q_{kl}\right)q_{kk}q_{jk}q_{kj}
    \right) q_{lk}q_{kk} [x_jx_k]_c \otimes [x_kx_l]_c.
\end{align*}
Now, we calculate
\begin{align*}
    1-q_{kj}q_{jk} -&(1-q_{lk}q_{kl})q_{kk}q_{jk}q_{kj} = 1-q_{kj}q_{jk} - q_{kk}q_{jk}q_{kj}+q_{kk}^{-1}
    \\ &= q_{kk}^{-1}(1+q_{kk})(1-q_{kk}q_{kj}q_{jk})=0,
\end{align*}
so $u \in \cP(\bB)$. \edem

\medskip

\begin{lema}\label{lemacoproB}
Let $\bB$ be a braided graded Hopf algebra provided with an
inclusion of braided vector spaces $V \hookrightarrow \cP(\bB)$.
Assume that

\begin{itemize}
    \item there exist $1 \leq k \neq j \leq \theta$ such that $m_{kj}=2,m_{jk}=1$;
    \item the following relations hold in $\bB$:
\begin{itemize}
    \item[$\ast$] $(\ad x_s)^{m_{st}+1}x_t=0$, for all $1 \leq s \neq t \leq \theta$;
    \item[$\ast$] $x_s^{m_{st}+1}=0$ for each $s$ such that $q_{ss}^{m_{st}}q_{st}q_{ts} \neq 1$, for some $t \neq s$.
\end{itemize}
\end{itemize}

\emph{(a)} If $v:= \left[ (\ad x_k)^2x_j, (\ad x_k)x_j \right]_c$, then there exists $b \in \kk$ such that
\begin{equation} \label{coproductGSR2}
\Delta (v) = v \otimes 1 + 1 \otimes v + b
(1-q_{kk}^2q_{kj}^2q_{jk}^2q_{jj}) x_k^{3} \otimes x_j^{2}.
\end{equation}

\emph{(b)} Assume that there exists $l \neq j,k$ such that $m_{jl}=m_{lj}=1$,  $m_{kl}=m_{lk}=0$,  and that \eqref{relA} is valid in $\bB$. Call $$ w := \left[ (\ad x_k)^2(\ad x_j)x_l, (\ad x_k)x_j \right]_c,$$
then there exist constants $b_1, b_2 \in \kk$, such that
\begin{eqnarray} \label{coproductGSR3}
\Delta (w) &=& w \otimes 1 + 1 \otimes w + b_1 v \otimes x_l
\\ && + b_2(1-q_{kk}^2q_{kj}q_{jk}) x_k^3 \otimes \left( (\ad x_j)x_l
\right)x_j. \nonumber
\end{eqnarray}
\end{lema}

\bdem (a) Note that $F_j(v)=0$, since $v$ is a braided
commutator of two elements in $\ker F_j$. Also, using
\eqref{idjac},
$$ [(\ad x_k)^2 x_j, x_j]_c = q_{kj} ( q_{jj} - q_{kk}) [x_kx_j]_c^2, $$
so we calculate
\begin{align*}
   F_k (v) =& (1+q_{kk}) (1-q_{kk}q_{kj}q_{jk})[x_kx_j]_c^2+
    q_{kk}^2q_{jk} (1-q_{kj}q_{jk}) [x_k^2x_j]_c x_j
    \\ & - q_{kk}^2q_{kj}^2q_{jk}q_{jj} (1-q_{kj}q_{jk})x_j
    [x_k^2x_j]_c
    \\ & - q_{kk}^3q_{kj}^2q_{jk}^2q_{jj}(1+q_{kk}) (1-q_{kk}q_{kj}q_{jk})[x_kx_j]_c^2
    \\ =& (q_{kk}^2q_{jk}q_{kj} (1-q_{kj}q_{jk})( q_{jj} - q_{kk}) + (1+q_{kk})
    \\ & (1-q_{kk}q_{kj}q_{jk}) (1- q_{kk}^3q_{kj}^2q_{jk}^2q_{jj}) [x_kx_j]_c^2 =0
\end{align*}
since the coefficient of $[x_kx_j]_c^2$ is zero for each possible braiding. Then, $$\de_{1,4}(v) = x_k \otimes
F_k(v)=0.$$

Also, $D_k(v)=0$, and we calculate
\begin{align*}
    D_j(v) =& (1-q_{kj}q_{jk}) \left( [x_k^2x_j]x_k + (1-q_{kk}q_{kj}q_{jk})q_{jk}q_{jj} x_k^2
    [x_kx_j]_c \right.
    \\  & \left. - q_{kk}^2q_{kj}^2q_{jk}q_{jj} (1-q_{kk}q_{kj}q_{jk}) [x_kx_j]_c x_k^2  - q_{kk}^2q_{kj}^2q_{jk}^3q_{jj}^2 x_k[x_k^2x_j]\right)
    \\  =&  \left( 1+ (1+q_{kk}) (1-q_{kk}q_{kj}q_{jk})q_{kk}q_{kj}q_{jk}q_{jj} - q_{kk}^4q_{kj}^3q_{jk}^3q_{jj}^2 \right)
    \\ & \ \ (1-q_{kj}q_{jk})[x_k^2x_j]x_k,
\end{align*}
where we reorder the hyperwords and use that $(\ad x_k)^3 x_j=0$;
also,
\begin{equation}\label{coefficient}
1+ (1+q_{kk}) (1-q_{kk}q_{kj}q_{jk})q_{kk}q_{kj}q_{jk}q_{jj} - q_{kk}^4q_{kj}^3q_{jk}^3q_{jj}^2 =0,
\end{equation}
by calculation for each possible braiding. Then, $$\de_{4,1}(v) = D_j(v) \otimes x_j=0.$$

To finish, we use that
\begin{align*}
\de(v) =& \de ((\ad_cx_k)^2x_j) \de ((\ad_cx_k)x_j)
 \\  & - \chi (2e_k+e_j,e_k+e_j ) \de ((\ad_cx_k)x_j) \de ((\ad_cx_k)^2x_j).
\end{align*}
Looking at the terms in $\bB^3 \otimes \bB^2$ and $\bB^2 \otimes
\bB^3$, and using the definition of braided commutator, we obtain
\begin{align*}
    \de _{32}(v) =& (1- q_{kk}^4q_{kj}^3q_{jk}^3q_{jj}^2)
    [x_k^2x_j]_c \otimes \left[ x_kx_j \right]_c + (1+q_{kk})(1-q_{kk}q_{kj}q_{jk})
    \\ & q_{kk}q_{kj}q_{jk}q_{jj} \left( x_k \left[ x_kx_j \right]_c - q_{kk}q_{kj} \left[ x_kx_j \right]_c x_k
    \right)\otimes \left[ x_kx_j \right]_c
    \\ & + 1-q_{kj}q_{jk})^2 (1-q_{kk}^2q_{kj}q_{jk}) (1-q_{kk}^2q_{kj}^2q_{jk}^2q_{jj}) x_k^{3} \otimes x_j^{2}
    \\ =& \left(1 + (1+q_{kk})(1-q_{kk}q_{kj}q_{jk})
    q_{kk}q_{kj}q_{jk}q_{jj}-
    q_{kk}^4q_{kj}^3q_{jk}^3q_{jj}^2 \right) \left[ x_k^2 x_j
    \right]_c
    \\ & \otimes \left[ x_kx_j \right]_c + b_1 (1-q_{kk}^2q_{kj}^2q_{jk}^2q_{jj}) x_k^{3} \otimes x_j^{2}.
\end{align*}
Also,
\begin{align*}
\de _{23}(v) =& (1-q_{kk}q_{kj}q_{jk})(1-q_{kj}q_{jk}) x_k^2
\otimes \left( (1+q_{kk}) q_{kk}q_{jk} \left[ x_kx_j \right]_c
    x_j \right.
    \\ & - (1+q_{kk}) q_{kk}^2q_{kj}^2q_{jk}^2q_{jj} x_j \left[ x_kx_j \right]_c  + x_j\left[ x_kx_j \right]_c
    \\ & \left. - q_{kk}^4q_{kj}^2q_{jk}^3q_{jj} \left[ x_kx_j \right]_c x_j \right)
    \\ = &  \left(1- q_{kk}^4q_{kj}^3q_{jk}^3q_{jj}^2 +
(1+q_{kk})(1-q_{kk}q_{kj}q_{jk}) q_{kk}q_{kj}q_{jk}q_{jj} \right)
    \\ & (1-q_{kk}q_{kj}q_{jk})(1-q_{kj}q_{jk})x_k^2 \otimes x_j \left[
x_kx_j \right]_c.
\end{align*}
Using \eqref{coefficient}, we obtain \eqref{coproductGSR2}.
\medskip

(b) We call $y= (\ad x_k)^2 (\ad x_j)x_l, \ z=(\ad
x_k)x_j$. Observe that $\de(w)= \de
(y)\de(z)-\chi(2\e_k+\e_j+\e_l, \e_k+\e_j)\de(z)\de(y)$, and
\begin{align*}
    \de(y) & = y \otimes 1+(1-q_{jl}q_{lj}) (\ad x_k)^2x_j \otimes x_l
    \\ & + (1-q_{kj}q_{jk}) (1-q_{kk}q_{kj}q_{jk}) x_k^2 \otimes (\ad x_j)x_l
    \\ & +  (1+q_{kk}) (1-q_{kk}q_{kj}q_{jk}) x_k \otimes (\ad x_k)(\ad x_j)x_l +1 \otimes y,
    \\ \de(z)&= z \otimes 1 + (1-q_{kj}q_{jk}) x_k \otimes x_j + 1 \otimes z.
\end{align*}

From \eqref{14}, $D_k(w)=0$, and from \eqref{15},
\begin{align*}
D_l(w) &= (1-q_{lj}q_{jl})q_{lk}q_{lj} \left[ (\ad x_k)^2x_j, (\ad
x_k)x_j \right]_c,
\\ D_j(w) &= -(1-q_{kj}q_{jk})q_{kk}^{-2}q_{kj}^{-1}q_{kl}^{-1} (\ad x_k)^3
(\ad x_j)x_l
\\ &= -(1-q_{kj}q_{jk})q_{kk}^{-2}q_{kj}^{-1}q_{kl}^{-1} [ (\ad x_k)^3 x_j,
x_l ]_c =0.
\end{align*}
where we use that $\left[x_k, x_l \right]_c=0$ and \eqref{idjac}
for the last equality. It follows that
$$\Delta_{51}(w)=(1-q_{lj}q_{jl})q_{lk}q_{lj} \left[ (\ad
x_k)^2x_j, (\ad x_k)x_j \right]_c \otimes x_l.$$

Also, $F_j(z)=F_j(y)=F_l(z)=F_l(y)=0$ by \eqref{paragenerar2} and
the properties of these skew derivations, so $F_j(w)=F_l(w)=0$. We
calculate
\begin{align*}
F_k(w) =& (1+q_{kk}) (1-q_{kk}q_{kj}q_{jk}) [x_kx_jx_l]_c
[x_kx_j]_c + q_{kk}^2q_{jk}q_{lk} (1-q_{kj}q_{jk})
\\ & [x_k^2x_jx_l]_c
x_j - \chi(2\e_k+\e_j+\e_l, \e_k+\e_j) \left( (1-q_{kj}q_{jk}) x_j
[x_k^2x_jx_l]_c \right.
\\ & \left. +(1+q_{kk}) (1-q_{kk}q_{kj}q_{jk})q_{kk}q_{jk}
[x_kx_jx_l]_c [x_kx_j]_c \right)
\\ = & q_{kk}^2q_{jk}q_{lk} (1-q_{kj}q_{jk}) \left[ [
x_k^2x_jx_l]_c , x_j \right]_c - (1+q_{kk}) (1-q_{kk}q_{kj}q_{jk})
\\ & q_{kk}^3q_{kj}^2q_{jk}^2q_{jj}q_{lj}q_{lk} \left[ [x_kx_j]_c,
[x_kx_jx_l]_c \right]_c
\\ = & q_{kk}^2q_{kj}q_{jk}q_{jj}q_{lj}q_{lk} \left(1-q_{kj}q_{jk}- (1+q_{kk}) (1-q_{kk}q_{kj}q_{jk}) q_{kk}q_{kj}q_{jk} \right)
\\ & \left[ [x_kx_j]_c, [x_kx_jx_l]_c \right]_c =0
\end{align*}
where we use \eqref{idjac} and \eqref{relA} in the third equality, and calculate that
\begin{equation}\label{coefficient2}
1-q_{kj}q_{jk}- (1+q_{kk}) (1-q_{kk}q_{kj}q_{jk}) q_{kk}q_{kj}q_{jk}=0,
\end{equation}
for each possible standard braiding. It follows that
$\Delta_{15}(w)=0$.

We calculate each of other terms of $\de(w)$ by direct
calculation. First,
\begin{align*}
    \Delta_{42}(w) =& \left( 1- \chi(2\e_k+\e_j+\e_l,\e_k+\e_j) \chi(\e_k+\e_j, 2\e_k+\e_j+\e_l)
    \right)y \otimes z
    \\ & + (1-q_{kj}q_{jk}) (1-q_{lj}q_{jl}) \left( q_{lk} [x_k^2x_j]_c x_k \otimes
    x_lx_j \right.
    \\ & \ \ \left. - \chi(2\e_k+\e_j+\e_l,\e_k+\e_j) q_{jk}^2q_{jj} x_k [x_k^2x_j]_c \otimes x_jx_l \right)
    \\ & + (1-q_{kj}q_{jk}) (1-q_{kk}q_{kj}q_{jk}) \left( \chi(\e_j+\e_l,\e_k+\e_j)
    x_k^2z \right.
    \\ & \ \ \left. - \chi(2\e_k+\e_j+\e_l,\e_k+\e_j) zx_k^2 \right)
    \otimes [x_jx_l]_c
    \\ = & (1-q_{kj}q_{jk}) q_{lk} \left( 1-q_{jk}q_{kj} + (1+q_{kk}) (1-q_{kk}q_{kj}q_{jk})q_{kk}q_{kj}q_{jk} \right)
    \\ & \ [x_k^2x_j]_c x_k \otimes [x_jx_l]_c =0.
\end{align*}

In a similar way we calculate
\begin{align*}
    \Delta_{33}(w) =& (1-q_{lj}q_{jl}) [x_k^2x_j] \otimes \left( x_l z -
    \chi(2\e_k+\e_j+\e_l,\e_k+\e_j)\right.
    \\ & \ \left. \chi(\e_k+\e_j, \e_k+\e_j+\e_l) zx_l  \right)+ (1+q_{kk}) (1-q_{kk}q_{kj}q_{jk})
    \\ & \ \chi(\e_k+\e_j+\e_l,\e_k+\e_j) \left( x_kz-q_{kk}q_{kj}zx_k
    \right) \otimes [x_kx_jx_l]_c
    \\ & + (1-q_{kk}q_{kj}q_{jk}) (1-q_{kj}q_{jk})^2 x_k^3 \otimes
    \left( \chi(\e_j+\e_l,\e_k)[x_jx_l]_c x_j \right.
    \\ & \ \left. - \chi(2\e_k+\e_j+\e_l,\e_k+\e_j)\chi(\e_j,2\e_k)x_j [x_jx_l]_c \right)
    \\ = & \left( (1+q_{kk})(1-q_{kk}q_{kj}q_{jk})-q_{kk}q_{kj}q_{jk}q_{jj}(1-q_{lj}q_{jl}) \right)
    \\ & \  \chi(\e_k+\e_j+\e_l,\e_k+\e_j)[x_k^2x_j]_c \otimes [ x_kx_jx_l]_c
    \\ & + (1-q_{kk}q_{kj}q_{jk}) (1-q_{kj}q_{jk})^2
    (1-q_{kk}^2q_{kj}q_{jk}) x_k^3 \otimes [x_jx_l]_c x_j,
\end{align*}
and the coefficient of $[x_k^2x_j]_c \otimes [ x_kx_jx_l]_c$ is
zero (we calculate it for each possible standard braiding). Also,
\begin{align*}
    \Delta_{24}(w) =& (1-q_{kk}q_{kj}q_{jk}) (1-q_{kj}q_{jk})
    x_k^2 \otimes\left( (1+q_{kk}) \chi(\e_k+\e_j+\e_l,\e_k)
    \right.
    \\ & \ \ [x_kx_jx_l]_c x_j  -  (1+q_{kk}) \chi(2\e_k+\e_j+\e_l,\e_k+\e_j)q_{jk}x_j
    [x_kx_jx_l]_c
    \\ & \ \left. - \chi(2\e_k+\e_j+\e_l,\e_k+\e_j) \chi(\e_k+\e_j, 2\e_k) \left[ [x_kx_j]_c, [x_jx_l]_c \right]_c \right)
    \\ = & (1-q_{kk}q_{kj}q_{jk}) (1-q_{kj}q_{jk}) \chi(\e_j+\e_l,
    \e_k+\e_j)q_{kj}
    \\ & \ \left( q_{kk}(1-q_{kk}q_{kj}q_{jk})-q_{jj}(1-q_{jl}q_{lj})  \right)
    x_k^2 \otimes x_j [x_kx_jx_l]_c =0
\end{align*}

From the above calculations, we obtain \eqref{coproductGSR3}. \edem

\subsection{Presentation when the type is $A_{\theta}$}\label{subsection:An}
\

In this subsection we shall consider $V$ a standard braided vector
space of type $A_{\theta}$, and $\bB$ a $\zt$-graded algebra, provided
with an inclusion of vector spaces $V \hookrightarrow \bB^1= \oplus_{1 \leq j \leq \theta} \bB^{\e_j}$. We can extend the braiding to $\bB$ by
    \[ c(u \otimes v)= \chi(\alpha, \beta) v \otimes u, \quad u \in \bB^{\alpha}, v \in \bB^{\beta}, \ \alpha, \beta \in \N^{\theta}.
\]
We assume that
\begin{align*}
    x_i^{2}&=0  &\mbox{if }q_{ii}=-1,
    \\ \ad_c x_i (x_j) &= 0 &\mbox{if } \left|j-i\right|>1,
    \\ (\ad_c x_i)^2 (x_j)_c &= 0 &\mbox{if } \left|j-i\right|=1,
    \\ \left[ (\ad_c x_i) (ad_c x_{i+1})x_{i+2}, x_{i+1} \right]_c &= 0 & 2 \leq i \leq \theta-1,
\end{align*}
are valid on $\bB$. Using the same notation as in Subsection \ref{subsection:generators},
    \[   x_{\e_i}=x_i, \qquad x_{\ub_{i,j}}:=\left[ x_i, x_{\ub_{i+1,j}} \right]_c \quad (i<j). \]

\begin{lema}\label{lemaA2}
Let $1 \leq i \leq j < p \leq r \leq \theta$. The following relations hold in $\bB$:
\begin{eqnarray}
    \left[ x_{\ub_{ij}}, x_{\ub_{pr}} \right]_c &=& 0, \quad  p-j \geq 2; \label{A1}
    \\ \left[ x_{\ub_{ij}}, x_{\ub_{j+1,r}} \right]_c &=& x_{\ub_{ir}}. \label{A2}
\end{eqnarray}
\end{lema}
\bdem Note that $x_{\ub_{pr}}$ belongs to the subalgebra generated by
$x_p, \ldots, x_r$, and $\left[ x_{\ub_{ij}}, x_s \right]_c=0$,
for each $p \leq s \leq r$. Then, \eqref{A1} is deduced from this
fact.
\smallskip

We prove \eqref{A2} by induction on $j-i$: if $i=j$, it is exactly
the definition of $x_{\ub_{ir}}$. To prove the inductive step, we
use the inductive hypothesis, \eqref{A1} and \eqref{idjac} (the braided Jacobi identity) to obtain
\begin{align*}
    \left[ x_{\ub_{i,j+1}}, x_{\ub_{j+2,r}} \right]_c =& \left[ \left[ x_{\ub_{ij}}, x_{i+1} \right]_c , x_{\ub_{j+2,r}} \right]_c = \left[ x_{\ub_{ij}}, \left[ x_{i+1} , x_{\ub_{j+2,r}} \right]_c  \right]_c
    \\ =& \left[ x_{\ub_{ij}}, x_{\ub_{j+1,r}} \right]_c = x_{\ub_{ir}},
\end{align*}
and \eqref{A2} is also proved.
\edem

\begin{lema} \label{lemaA3}
If $i<p \leq r <j$, the following relation holds in $\bB$:
\begin{equation}
    \left[ x_{\ub_{ij}}, x_{\ub_{pr}} \right]_c = 0. \label{A3}
\end{equation}
\end{lema}
\bdem When $p=r=j-1$ and $i=j-2$, note that this is exactly $$\left[
(\ad_c x_i) (ad_c x_{i+1})x_{i+2}, x_{i+1} \right]_c = 0.$$
Then, we have by \eqref{idjac}
$$[ x_{\ub_{i-1,j}}, x_{j-1} ]_c = [ [x_{i-1},x_{\ub_{i,j}}]_c, x_{j-1}]_c = [ x_{i-1},
[x_{\ub_{i,j}}, x_{j-1}]_c ]_c.$$ We assume that $j-i>2$, so $ [x_{i-1},
x_{j-1}]_c =0$ by hypothesis on $\bB$. Then, we prove the case $p=r=j-1$ by induction on
$p-i$.

Using \eqref{idjac} and \eqref{A2}, we also have
\begin{eqnarray*}
[x_{\ub_{i,j+1}}, x_{p} ]_c &=& [ [x_{\ub_{i,j}}, x_{j+1}]_c,
x_{p}]_c = [ x_{\ub_{i,j}}, [x_{j+1}, x_{p}]_c ]_c
\\ && + q_{j+1,p}[x_{\ub_{i,j}}, x_{j-1}]_c x_{j+1} - \chi (\ub_{i,j}, \e_{j+1}) x_{j+1} [x_{\ub_{i,j}},
x_{j-1}]_c,
\end{eqnarray*}
so using that $ [x_{j+1}, x_{p}]_c =0$ if $j>p$, by induction on
$j-p$ we prove \eqref{A3} for the case $p=r$.

For the general case, we use \eqref{idjac} one more time as follows
\begin{align*}
[x_{\ub_{i,j}}, x_{\ub_{p,r+1}} ]_c &=  [ x_{\ub_{i,j}},
[x_{\ub_{pr}}, x_{r+1}]_c ]_c = [ [x_{\ub_{i,j}}, x_{\ub_{pr}}]_c,
x_{r+1}]_c
\\ & - \chi (\ub_{pr}, \e_{r+1}) [x_{\ub_{ij}}, x_{r+1}]_c x_{\ub_{pr}} + \chi (\ub_{ij}, \ub_{pr}) x_{\ub_{pr}} [x_{\ub_{ij}},
x_{r+1}]_c,
\end{align*}
and we prove \eqref{A3} by induction on $r-p$.
\edem

\begin{lema} \label{lemaA4}
The following relations hold in $\bB$:
\begin{eqnarray}
    \left[ x_{\ub_{ij}}, x_{\ub_{ip}} \right]_c &=& 0, \quad i \leq j <p; \label{A4}
    \\ \left[ x_{\ub_{ij}}, x_{\ub_{pj}} \right]_c &=& 0, \quad i <p \leq j. \label{A5}
\end{eqnarray}
\end{lema}
\bdem To prove \eqref{A4}, note that if $i=j=p-1$, we have
$$ \left[ x_{\ub_{ii}}, x_{\ub_{i,i+1}} \right]_c = \left[ x_i, [ x_i, x_{i+1} ]_c \right]_c = (\ad x_i)^2x_{i+1}=0. $$
As $[x_i, x_{\ub_{i+2,p}}]_c=0$ for each $p > i+1$ by \eqref{A1},
we use \eqref{idjac}, the previous case and \eqref{A2} to obtain
$$ \left[ x_{\ub_{ii}}, x_{\ub_{ip}} \right]_c = \left[ x_{\ub_{ii}}, [x_{\ub_{i,i+1}}, x_{\ub_{i+2,p}} ]_c \right]_c = 0.$$
Now, if $i<j<p$, from \eqref{A1} and the relations between the
$q_{st}$ we obtain
$$ \left[ x_{\ub_{i+1,j}}, x_{\ub_{ip}} \right]_c = - \chi ( \ub_{ip}, \ub_{i+1,j})  \left[ x_{\ub_{ip}}, x_{\ub_{i+1,j}} \right]_c=0. $$
Using \eqref{idjac} and the previous case we conclude
$$ \left[ x_{\ub_{ij}}, x_{\ub_{ip}} \right]_c = \left[ \left[ x_{\ub_{ii}}, x_{\ub_{i+1,j}}  \right]_c, x_{\ub_{ip}} \right]_c =0 .$$

The proof of \eqref{A5} is analogous.\edem

\begin{lema} \label{lemaA5}
If $i<p \leq r <j$, the following relation holds in $\bB$:
\begin{equation}
    \left[ x_{\ub_{ir}}, x_{\ub_{pj}} \right]_c = \chi(\ub_{ir}, \ub_{pr}) \left( 1- q_{r,r+1} q_{r+1,r} \right) x_{\ub_{pr}}x_{\ub_{ij}}. \label{2B}
\end{equation}
\end{lema}
\bdem We calculate
\begin{align*}
\left[ x_{\ub_{ir}}, x_{\ub_{pj}} \right]_c &= \left[ x_{\ub_{ir}}
, \left[ x_{\ub_{pr}}, x_{\ub_{r+1,j}} \right]_c  \right]_c
\\ &= \chi(\ub_{ir},\ub_{pr})x_{\ub_{pr}} x_{\ub_{ij}} - \chi(\ub_{pr}, \ub_{r+1,j})x_{\ub_{ij}} x_{\ub_{pr}}
\\ &= \left( \chi(\ub_{ir},\ub_{pr})-\chi(\ub_{ij},\ub_{pr})\chi(\ub_{pr}, \ub_{r+1,j}) \right) x_{\ub_{pr}}x_{\ub_{ij}}
\\ &= \chi(\ub_{ir}, \ub_{pr}) \left( 1- \chi(\ub_{pr}, \ub_{r+1,j})\chi(\ub_{r+1,j} , \ub_{pr} ) \right) x_{\ub_{pr}}x_{\ub_{ij}},
\end{align*}
where we use \eqref{A2} in the first equality, \eqref{idjac} in
the second, \eqref{A5} in the third and the relation between the
$q_{ij}$ in the last. \edem

We prove the main Theorem of this subsection, namely, the presentation by generators and
relations of the Nichols algebra associated to $V$.

\begin{theorem}\label{presentationA}
Let $V$ be a standard braided vector space of type $A_{\theta}$,
$\theta= \dim V$, and $C=(a_{ij})_{i,j \in \unon}$ the
corresponding Cartan matrix of type $A_{\theta}$.

The Nichols algebra $\bB(V)$ is presented by generators $x_{i}$,
$1\le i \le \theta$, and relations
\begin{eqnarray*}
x_{\alpha}^{N_{\alpha}} &=& 0, \quad \alpha \in \de^{+};
\\ ad_{c}(x_{i})^{1-a_{ij}}(x_{j})  &=& 0, \quad i\neq j;
\\  \left[ (\ad x_{j-1})(\ad x_j)x_{j+1}, x_j \right]_c &=& 0, \quad 1 <
j < \theta, \, q_{jj}=-1.
\end{eqnarray*}

Moreover, the following elements constitute a basis of $\bB(V)$:
\begin{equation}\label{baseAn}
x_{\beta_{1}}^{h_{1}} x_{\beta_{2}}^{h_{2}} \dots
x_{\beta_{P}}^{h_{P}}, \qquad 0 \le h_{j} < N_{\beta_j},
\text{ if }\, \beta_j \in S_I, \quad  1\le j \le P.
\end{equation}
\end{theorem}
\bdem From Corollary \ref{corollary:genPBW} and the definitions of
$x_{\alpha}$'s, we know that the last statement about the PBW
basis is true.
\medskip

Now, let $\bB$ be the algebra presented by generators $x_1,\ldots,
x_{\theta}$ and the relations of the Theorem. From Lemmata
\ref{buscarmij}, \ref{lemaAn} and Proposition
\ref{heigthgenerators} we have a canonical epimorphism $\phi: \bB
\rightarrow \bB(V)$. Note that last relation also holds in $\bB$
for $q_{jj} \neq 1$, by Lemma \ref{lemaAn}, (2).

The proof is similar to the ones of \cite[Lemma 3.7]{AD} and
\cite[Lemma 6.12]{AS5}. Consider the subspace $\cI$ of $\bB$ generated by the
elements in \eqref{baseAn}. Using Lemmata \ref{lemaA2},
\ref{lemaA3}, \ref{lemaA4} and \ref{lemaA5} we prove that $\cI$ is an ideal.
But $1 \in \cI$, so $\cI=\bB$.

The image of the elements in \eqref{baseAn} by $\phi$ are a basis
of $\bB(V)$, so $\phi$ is an isomorphism.\edem

Note that the presentation and the dimension of $\bB(V)$ agree
with the results presented in \cite{AD} and \cite{AS5}.

\subsection{Presentation when the type is $B_{\theta}$} \label{subsection:Bn}
\

Now, we shall consider $V$ a standard braided vector space of type $B_{\theta}$, and $\bB$ a $\zt$-graded algebra, provided with an inclusion of vector spaces $V \hookrightarrow \bB^1= \oplus_{1 \leq j \leq \theta} \bB^{\e_j}$. Then, we can extend the braiding to $\bB$. We assume that the following relations hold in $\bB$
\begin{align*}
x_i^{2}&=0 & \mbox{if }q_{ii}=-1,
\\ x_1^{3}&=0  &\mbox{if }q_{11}\in \G_3,
\\ (\ad_c x_i) x_j &=0,  & \left|j-i\right|>1;
\\ (\ad_c x_i)^2x_j &= 0,  & \left|j-i\right|=1, \ i\neq 1;
\\ \left[ (\ad_c x_i)(ad_c x_{i+1})x_{i+2}, x_{i+1} \right]_c &= 0, & 2 \leq i \leq \theta;
\\ (\ad_c x_1)^3x_2 &= 0;
\\ \left[ (\ad_c x_1)^2x_2, (ad_c x_1)x_2 \right]_c &= 0,
\\ \left[ (\ad x_1)^2(\ad x_2)x_3, (\ad x_1)x_2 \right]_c &=0.
\end{align*}
Using the same notation as in Subsection \ref{subsection:generators},
    \[   x_{\vb_{ij}}= \left[ x_{\ub_{1i}} , x_{\ub_{1j}}\right]_c, \qquad  1 \leq i<j \leq \theta. \]

From the proof of relations corresponding the $A_{\theta}$ case,
relations \eqref{A1}, \eqref{A2}, \eqref{A3}, \eqref{A5} and
\eqref{2B} are valid for $i \geq 1$, but for relation \eqref{A4}
we must consider $i>1$.

\begin{lema}\label{reorder1}
Let $1 \leq s <t$, $1<k \leq j$. The following relations hold in $\bB$:
\begin{equation*}
[x_{\vb_{st}}, x_{\ub_{kj}} ]_c \begin{cases} = 0 & t+1<k;
\\  = x_{\vb_{sj}} & t+1=k < j;
\\  = 0 & s+1<k \leq j \leq t;
\\  =\chi(\vb_{st}, \ub_{kt}) (1-q_{t,t+1}q_{t+1,t}) x_{\ub_{kt}}
x_{\vb_{sj}} & s+1 <k \leq t <j;
\\  = \chi(\ub_{1t},\ub_{s+1,j}) x_{\vb_{jt}} & s+1 =k \leq j <t ;
\\  = ( \chi(\ub_{1t},\ub_{s+1,t}) - \chi(\ub_{1s},\ub_{1t}) ) x_{\ub_{1t}}^{2} & s+1 =k, j=t ;
\\  \in \kk x_{\vb_{tj}} + \kk x_{\ub_{1j}} x_{\ub_{1t}} +\kk x_{\ub_{s+1,j}} x_{\vb_{sj}} & s+1 =k \leq t <j ;
\\  =  \gamma^{kj}_{st} x_{\ub_{ks}} x_{\vb_{jt}} & k \leq s < j \leq t ;
\\  \in \kk x_{\ub_{ks}}x_{\vb_{tj}} + \kk x_{\ub_{ks}}x_{\ub_{1j}} x_{\ub_{1t}} +\kk x_{\ub_{kt}} x_{\vb_{sj}} & k\leq s < t <j;
\\  0 & k \leq j \leq s,
\end{cases}
\end{equation*}
where $\gamma^{kj}_{st} = \chi(\ub_{1t},\ub_{kj})
\chi(\ub_{1s},\ub_{ks}) (1-q_{s,s+1}q_{s+1,s})$.
\end{lema}
\bdem The first, the third and the last cases follow from the fact
that
$$ [x_{\ub_{1s}}, x_{\ub_{kj}}]_c=[x_{\ub_{1t}}, x_{\ub_{kj}}]_c= 0
$$ using \eqref{A1}, \eqref{A3}, \eqref{A4} or \eqref{A5}(depending on each case), and
\eqref{idjac}.

For the second, use that $\left[x_{\ub_{1s}},
x_{\ub_{t+1,j}}\right]_c=0$, \eqref{A2} and \eqref{idjac}, to
obtain
\begin{align*}
x_{\vb_{sj}} &= \left[ x_{\ub_{1s}}, x_{\ub_{1j}} \right]_c =
\left[ x_{\ub_{1s}}, \left[ x_{\ub_{1t}}, x_{\ub_{t+1,j}}
\right]_c \right]_c
\\ &= \left[ \left[ x_{\ub_{1s}},  x_{\ub_{1t}} \right]_c , x_{\ub_{t+1,j}} \right]_c  = \left[ x_{\vb_{st}}, x_{\ub_{t+1,j}} \right]_c.
\end{align*}

For the fourth, use \eqref{idjac} and the third case to calculate
\begin{align*}
[x_{\vb_{st}},  x_{\ub_{kj}} ]_c =& [x_{\vb_{st}}, [x_{\ub_{kt}}, x_{\ub_{t+1,j}} ]_c ]_c
\\ = & \chi(\vb_{st},\ub_{kt}) x_{\ub_{kt}} x_{\vb_{sj}} - \chi(\ub_{kt},\ub_{t+1,j}) x_{\vb_{sj}} x_{\ub_{kt}}
\\ = & \chi(\vb_{st},\ub_{kt}) (1- \chi(\ub_{kt},\ub_{t+1,j})\chi(\ub_{t+1,j}, \ub_{kt}) ) x_{\ub_{kt}} x_{\vb_{sj}}.
\end{align*}

For the fifth, note that $\chi(\ub_{1t},\ub_{s+1,j})^{-1} =
\chi(\ub_{s+1,j},\ub_{1t})$. Then, use \eqref{A2}, \eqref{A3} and
\eqref{idjac} to prove that
\begin{align*}
[x_{\vb_{st}},  x_{\ub_{s+1,j}} ]_c =& [ [x_{\ub_{1s}},
x_{\ub_{1t}} ]_c, x_{\ub_{s+1,j}} ]_c
\\ = & \chi(\ub_{1t},\ub_{s+1,j}) x_{\ub_{1j}} x_{\ub_{1t}}
- \chi(\ub_{1s},\ub_{1t}) x_{\ub_{1t}} x_{\ub_{1s}}
\\ = & \chi(\ub_{1t},\ub_{s+1,j}) (x_{\ub_{1j}} x_{\ub_{1t}}
- \chi(\ub_{1j},\ub_{1t}) x_{\ub_{1t}} x_{\ub_{1s}} ).
\end{align*}
The sixth case follows in a similar way.

For the seventh case, we use \eqref{idjac}, \eqref{der} and the
previous case to calculate
\begin{align*}
[x_{\vb_{st}},  & x_{\ub_{s+1,j}} ]_c = [ x_{\vb_{st}}, [
x_{\ub_{s+1,t}} , x_{\ub_{t+1,j}} ]_c ]_c
\\ = & ( \chi(\ub_{1t},\ub_{s+1,t}) - \chi(\ub_{1s},\ub_{1t}) ) [ x_{\ub_{1t}}^{2}, x_{\ub_{t+1,j}}]
\\ & +\chi(\vb_{st},\ub_{s+1,t}) x_{\ub_{s+1,t}} x_{\vb_{sj}} - \chi(\ub_{s+1,t},\ub_{t+1,j}) x_{\vb_{sj}} x_{\ub_{s+1,t}}
\\ = & ( \chi(\ub_{1t},\ub_{s+1,t}) - \chi(\ub_{1s},\ub_{1t}) ) ( (x_{\vb_{tj}} + \chi(\ub_{1t},\ub_{1j}) x_{\ub_{1j}}x_{\ub_{1t}})
\\ & + \chi(\ub_{1t},\ub_{t+1,j}) x_{\ub_{1j}}x_{\ub_{1t}} ) - \chi(\ub_{s+1,t},\ub_{t+1,j}) x_{\vb_{tj}}
\\ & +( \chi(\vb_{st},\ub_{s+1,t})- \chi(\ub_{s+1,t},\ub_{t+1,j})\chi(\vb_{sj},\ub_{s+1,t}) ) x_{\ub_{s+1,t}}
x_{\vb_{sj}}.
\end{align*}
We use the previous cases, \eqref{A3} and \eqref{2B} to calculate
for the eighth case
\begin{align*}
[x_{\vb_{st}},  x_{\ub_{kj}} ]_c =& [ [x_{\ub_{1s}}, x_{\ub_{1t}}
]_c, x_{\ub_{kj}} ]_c
\\ = & \chi(\ub_{1t},\ub_{kj}) ( \chi(\ub_{1s},\ub_{ks}) (1-q_{s,s+1}q_{s+1,s})  x_{\ub_{ks}} x_{\ub_{1j}})  x_{\ub_{1t}}
\\ & - \chi(\ub_{1s},\ub_{1t}) x_{\ub_{1t}}  ( \chi(\ub_{1s},\ub_{ks}) (1-q_{s,s+1}q_{s+1,s})  x_{\ub_{ks}} x_{\ub_{1j}})
\\ = & \gamma^{kj}_{st} x_{\ub_{ks}} (x_{\ub_{1j}} x_{\ub_{1t}} - \chi(\ub_{1j},\ub_{1t}) x_{\ub_{1t}} x_{\ub_{1j}} ).
\end{align*}
To finish, we prove the ninth case in a similar way as follow
\begin{align*}
[x_{\vb_{st}},   x_{\ub_{kj}} ]_c = &[ x_{\vb_{st}}, [
x_{\ub_{kt}} , x_{\ub_{t+1,j}} ]_c ]_c
\\ = & \gamma^{kt}_{st}(1-q_{\vb_{1t}}) [ x_{\ub_{ks}} x_{\ub_{1t}}^2 ,
x_{\ub_{t+1,j}}]+\chi(\vb_{st},\ub_{k,t}) x_{\ub_{kt}}
x_{\vb_{sj}}
\\ & - \chi(\ub_{kt},\ub_{t+1,j}) x_{\vb_{sj}} x_{\ub_{kt}}.
\end{align*}
\edem
We consider the remaining commutator $\left[ x_{\vb_{st}}, x_{\ub_{jk}} \right]_c$: when $j=1$.

\begin{lema}
Let $s<t$ in $\unon$. The following relations hold in $\bB$
\begin{eqnarray}\label{3B}
\left[ x_{\vb_{st}} , x_{\ub_{1k}} \right]_c &=& 0, \quad s < k
\leq t;
\\ \left[ x_{\ub_{1s}} , x_{\vb_{st}}  \right]_c &=& 0. \label{3B'}
\end{eqnarray}
\end{lema}
\bdem By hypothesis we have
\begin{align*}
    \left[ x_{\vb_{12}} , x_{\ub_{12}} \right]_c&= \left[ (\ad_c x_1)^2 x_2, (\ad_c x_1) x_2 \right]_c = 0,
    \\ \left[ x_{\vb_{13}} , x_{\ub_{12}} \right]_c&= \left[ (\ad_c x_1)^2 (\ad_c x_2)x_3, (\ad_c x_1) x_2 \right]_c =0.
\end{align*}
For $t \geq 4$, $\left[x_{\ub_{4t}}, x_{\ub_{12}} \right]_c=0$ by
\eqref{A1}, and using  \eqref{idjac},
$$ \left[ x_{\vb_{1t}}, x_{\ub_{12}} \right]_c= \left[\left[ x_{\vb_{13}}, x_{\ub_{4t}} \right]_c , x_{\ub_{12}} \right]_c=0. $$
For each $k \leq t$ we have $\left[ x_{\ub_{1t}}, x_{\ub_{3k}} \right]_c= \left[ x_1, x_{\ub_{3k}} \right]_c=0$, so $\left[ x_{\vb_{1t}}, x_{\ub_{3k}} \right]_c=0$. Using \eqref{idjac} and \eqref{A2} we have
$$ \left[ x_{\vb_{1t}}, x_{\ub_{1k}} \right]_c = \left[ x_{\vb_{1t}}, \left[ x_{\ub_{12}},  x_{\ub_{3k}} \right]_c \right]_c  =0. $$

Now, consider $2 \leq s \leq k$. As $\left[ x_{\vb_{1t}},
x_{\ub_{1k}} \right]_c = \left[ x_{\ub_{2s}}, x_{\ub_{1k}}
\right]_c=0$ by previous results and \eqref{A3}, we conclude from
\eqref{der} and Lemma \ref{reorder1} that \eqref{3B} is valid in
the general case.

To prove \eqref{3B'}, for $s=1, t=2$ we have
$$ \left[ x_{\ub_{11}} , x_{\vb_{12}} \right]_c= \left[ x_1, x_{\vb_{12}} \right]_c = (\ad_c x_1)^3x_2=0. $$
Using that $\left[x_1, x_{\ub_{3t}} \right]_c=0$ if $t \geq 3$ and
\eqref{idjac}, we deduce that
$$ \left[ x_{\ub_{11}} , x_{\vb_{1t}} \right]_c= \left[x_1, \left[ x_{\vb_{12}}, x_{\ub_{3t}} \right]_c \right]_c=0. $$
If $1<s<t$, by the previous case we have
$$ \left[ x_{\ub_{1s}} , x_{\vb_{1t}} \right]_c= - \chi (x_{\ub_{1s}}, x_{\vb_{1t}})\left[ x_{\vb_{1t}}, x_{\ub_{1s}} \right]_c=0. $$
By \eqref{A5}, $\left[ x_{\ub_{1s}} , x_{\ub_{2s}} \right]_c=0$.
Also, $ [x_{\vb_{1t}}, x_{\ub_{2s}} ]_c = \chi(\ub_{1t},\ub_{2s})
x_{\vb_{st}}$, by Lemma \ref{reorder1}. Then, \eqref{3B'} follows
by \eqref{idjac} and the last three equalities.\edem

\begin{lema}\label{reorder2}
Let $s<k<t$. The following relations hold in $\bB$:
\begin{eqnarray}\label{4B}
    \left[ x_{\vb_{sk}}, x_{\ub_{1t}} \right]_c &=& \chi (\vb_{sk}, \ub_{1k} ) ( 1-q_{k,k+1}q_{k+1,k} )
    x_{\ub_{1k}}x_{\vb_{st}},
    \\ \left[ x_{\ub_{1s}}, x_{\vb_{kt}} \right]_c &=& \chi(\ub_{1s}, \ub_{1k}) (1+q_{\ub_{1k}})( 1-q_{k,k+1}q_{k+1,k} )
    x_{\ub_{1k}}x_{\vb_{st}}. \label{6B}
\end{eqnarray}
\end{lema}
\bdem The proof follows by \eqref{idjac}, the second
case of Lemma \ref{reorder1} and \eqref{3B'},
\begin{align*}
[ x_{\vb_{sk}}, & x_{\ub_{1t}} ]_c = \left[ x_{\vb_{sk}}, \left[
x_{\ub_{1k}} , x_{\ub_{k+1,t}} \right]_c \right]_c
\\  = & \chi (\vb_{sk}, \ub_{1k})x_{\ub_{1k}} x_{\vb_{st}} - \chi(\ub_{1k}, \ub_{k+1,t}) x_{\vb_{st}} x_{\ub_{1k}}
\\  = & \chi (\vb_{sk}, \ub_{1k}) \left( 1- \chi(\ub_{1k}, \ub_{k+1,t})\chi(\ub_{k+1,t}, \ub_{1k}) \right) x_{\ub_{1k}}
x_{\vb_{st}};
\\ [ x_{\ub_{1s}}, & x_{\vb_{kt}} ]_c = [ x_{\ub_{1s}}, [ x_{\ub_{1k}} , x_{\ub_{1t}} ]_c
]_c= \left[ x_{\vb_{sk}}, x_{\ub_{1t}} \right]_c \\ & + \chi
(\ub_{1s}, \ub_{1k} ) x_{\ub_{1k}}x_{\vb_{st}}- \chi(\ub_{1k},
\ub_{1t}) x_{\vb_{st}} x_{\ub_{1k}}
\\ = & \chi(\ub_{1s}, \ub_{1k}) (q_{\ub_{1k}}( 1-q_{k,k+1}q_{k+1,k} )+
1-q_{k,k+1}q_{k+1,k})x_{\ub_{1k}}x_{\vb_{st}}.
\end{align*}
\edem

We deal with the expression of the commutator of two words of
type $x_{\vb_{st}}$.

\begin{lema}\label{reorder3}
Let $s<t, k<j$, $s \leq k$, with $k \neq s$ or $j \neq t$. The following relations hold in $\bB$:
\begin{equation*}
    \left[ x_{\vb_{st}}, x_{\vb_{kj}} \right]_c \begin{cases}
    =0, & k<j \leq t,
\\ & k=s, t < j
\\ =  \chi (\vb_{st}, \vb_{kt})
(1-q_{t,t+1}q_{t+1,t})x_{\vb_{kt}} x_{\vb_{sj}}, & k <t <j;
\\ =  \chi (\vb_{st}, \ub_{1t})^2
(1-q_{t,t+1}q_{t+1,t}) &
\\ \quad (1-q_{\ub_{1t}} q_{t,t+1}q_{t+1,t})
x_{\ub_{1t}}^{2} x_{\vb_{sj}}, & k =t <j;
\\ \in \kk x_{\vb_{tj}} x_{\vb_{sk}}+\kk x_{\vb_{tk}} x_{\vb_{sj}}+ \kk x_{\ub_{1k}}x_{\ub_{1t}}x_{\vb_{sj}}, & t <k <j.
\end{cases}
\end{equation*}
\end{lema}
\bdem The first and the second cases follow from \eqref{idjac} and
\eqref{3B}, \eqref{3B'}, respectively. For the third case, we use the previous one and \eqref{idjac},
\begin{align*}
[ x_{\vb_{st}},  x_{\vb_{kj}} ]_c = & \left[ x_{\vb_{st}}, \left[
x_{\ub_{1k}} , x_{\ub_{1j}} \right]_c \right]_c
\\  = & \chi (\vb_{st}, \ub_{1k})x_{\ub_{1k}} ( \chi (\vb_{st}, \ub_{1t}) (1-q_{t,t+1}q_{t+1,t}) x_{\ub_{1t}}  x_{\vb_{sj}} )
\\  & - \chi(\ub_{1k}, \ub_{1j}) ( \chi (\vb_{st}, \ub_{1t}) (1-q_{t,t+1}q_{t+1,t}) x_{\ub_{1t}}  x_{\vb_{sj}} ) x_{\ub_{1k}}
\\  = & (1-q_{t,t+1}q_{t+1,t}) \left( \chi (\vb_{st}, \ub_{1k})\chi (\vb_{st}, \ub_{1t}) x_{\ub_{1k}} x_{\ub_{1t}}  x_{\vb_{sj}}
 \right.
\\  & \left. - \chi(\ub_{1k}, \ub_{1j}) \chi (\vb_{st}, \ub_{1t}) \chi (\vb_{sj}, \ub_{1k}) x_{\ub_{1t}}
x_{\ub_{1k}}x_{\vb_{sj}} \right)
\\  = & \chi (\vb_{st}, \ub_{1k}) \chi (\vb_{st}, \ub_{1t})
(1-q_{t,t+1}q_{t+1,t})
\\  & ( x_{\ub_{1k}}x_{\ub_{1k}} - \chi(\ub_{1k},
\ub_{1t}) x_{\ub_{1k}}x_{\ub_{1k}} ) x_{\vb_{sj}}.
\end{align*}
The fourth case is similar to the previous one.

To prove the last case, use \eqref{idjac} and Remark
\ref{reorder2} to express
\begin{align*}
[ x_{\vb_{st}},  x_{\vb_{kj}} ]_c = & \left[ x_{\vb_{st}}, \left[
x_{\ub_{1k}} , x_{\ub_{1j}} \right]_c \right]_c
\\  = & [ \chi (\vb_{st}, \ub_{1t}) (1-q_{t,t+1}q_{t+1,t}) x_{\ub_{1t}}  x_{\vb_{sk}}, x_{\ub_{1j}}  ]_c
\\  & + \chi (\vb_{st}, \ub_{1k})x_{\ub_{1k}} ( \chi (\vb_{st}, \ub_{1t}) (1-q_{t,t+1}q_{t+1,t}) x_{\ub_{1t}}  x_{\vb_{sj}} )
\\  & - \chi(\ub_{1k}, \ub_{1j}) ( \chi (\vb_{st}, \ub_{1t}) (1-q_{t,t+1}q_{t+1,t}) x_{\ub_{1t}}  x_{\vb_{sj}} ) x_{\ub_{1k}}.
\end{align*}
and the proof finish using \eqref{der} and the previous
identities. \edem

\medskip

\begin{theorem}\label{presentationB}
Let $V$ be a standard braided vector space of type $B_{\theta}$,
$\theta= \dim V$, and $C=(a_{ij})_{i,j \in \unon}$ the
corresponding Cartan matrix of type $B_{\theta}$.

The Nichols algebra $\bB(V)$ is presented by generators $x_{i}$,
$1\le i \le \theta$, and relations
\begin{eqnarray*}
x_{\alpha}^{N_{\alpha}} &=& 0, \quad \alpha \in \de^{+};
\\ ad_{c}(x_{i})^{1-a_{ij}}(x_{j})  &=& 0, \quad i\neq j;
\\  \left[ (\ad x_{j-1})(\ad x_j)x_{j+1}, x_j \right]_c &=& 0, \quad 1 <
j < \theta, \, q_{jj}=-1;
\\  \left[ (\ad x_1)^2x_2, (\ad x_1)x_2 \right]_c &=& 0, \quad q_{11} \in \G_3 \mbox{ or } q_{22}=-1 ; \label{Brel}
\\  \left[ (\ad x_1)^2 (\ad x_2) x_3, (\ad x_1)x_2 \right]_c &=& 0, \quad q_{11} \in \G_3 \mbox{ or } q_{22}=-1.\label{Brel2}
\end{eqnarray*}

Moreover, the following elements constitute a basis of $\bB(V)$:
\begin{equation}\label{baseBn}
x_{\beta_{1}}^{h_{1}} x_{\beta_{2}}^{h_{2}} \dots
x_{\beta_{P}}^{h_{P}}, \qquad 0 \le h_{j} \le N_{\beta_j} - 1,
\text{ if }\, \beta_j \in S_I, \quad  1\le j \le P.
\end{equation}
\end{theorem}

\bdem The proof is analogous to the corresponding of Theorem
\ref{presentationA}, since by previous Lemmata we express the commutator of two generators $x_{\alpha} < x_{\beta}$ as a linear combination of monotone hyperwords, whose greater hyperletter is great or equal than $x_{\beta}$.
\edem

\subsection{Presentation when the type is  $G_2$}\label{subsection:presentationG2}
\

We consider now standard braidings of $G_2$ type, with $m_{12}=3, m_{21}=1$.

\begin{lema}\label{lemaG2}
Let $\bB:=T(V)/I$, for some $I \in \bS$, such that in $\bB$ hold
\begin{align}
    (\ad x_1)^4x_2&=(\ad x_2)^2x_1=0; \label{condG}
    \\ x_i^{N_i}&=0, \qquad i=1,2, \ N_i:= \ord q_{ii}. \label{condG2}
\end{align}

\textbf{(a)} $\left[ x_1^3x_2x_1x_2 \right]_c=0
\Longleftrightarrow 4e_1+2e_2 \notin \de^{+}(\bB);$

\smallskip

Assume now that \textbf{(a)} holds in $\bB$. Then

\textbf{(b)} $\left[ (\ad x_1)^3x_2, (\ad x_1)^2x_2 \right]_c=0
\Longleftrightarrow 5e_1+2e_2 \notin \de^{+}(\bB);$

\smallskip

\textbf{(c)} $\left[ \left[ x_1^2x_2x_1x_2 \right]_c  ,
\left[x_1x_2\right]_c  \right]_c=0  \Longleftrightarrow 4e_1+3e_2
\notin \de^{+}(\bB);$

\smallskip

Assume also that \textbf{(b), (c)} holds in $\bB$, then

\textbf{(d)} $\left[  \left[x_1^{2}x_2\right]_c  , \left[
x_1^2x_2x_1x_2 \right]_c \right]_c=0  \Longleftrightarrow
5e_1+3e_2 \notin \de^{+}(\bB).$

\smallskip

In particular, all these relations hold when $V$ is a standard
braiding and $\bB=\bB(V)$ is finite dimensional.
\end{lema}

\bdem Order the letters $x_1< x_2$, and consider a PBW basis
as in Theorem \ref{basePBW}. We denote $\gamma_k:= \prod_{0 \leq j \leq k-1} (1-q_{11}^jq_{12}q_{21})$.
\smallskip

\emph{\textbf{(a)}} If the first assertion is true, then
$4e_1+2e_2 \notin \de^{+}(\bB)$ since there are no possible Lyndon words
in $S_I$: $x_1^3x_2x_1x_2$ is the unique Lyndon word such that
$x_1^3x_2, x_1x_2^2$ are not factors, and it is not in $S_I$
because of the hypothesis.

Reciprocally, if $4e_1+2e_2 \notin \de^{+}(\bB)$, then $\left[
x_1^3x_2x_1x_2 \right]_c$ is a linear combination of greater
hyperwords, and $\left[x_1x_2x_1^3x_2\right]_c, \left[
x_1^2x_1^2x_2 \right]$ are the unique greater hyperwords that are
not in $S_I$ and do not end in $x_1$ (we discard words ending in
$x_1$ since $\left[ x_1^3x_2x_1x_2 \right]_c$ is in $\ker D_1$).
So, taking their Shirshov decomposition, there exist $\alpha,
\beta \in \kk$ such that
\begin{equation}\label{palabra}
\left[ x_1^3x_2x_1x_2 \right]_c- \alpha \left[ x_1x_2\right]_c \left[x_1^3 x_2 \right]_c -\beta \left[ x_1^2x_2 \right]_c^2=0 .
\end{equation}
Note that $ \left[ x_1^3x_2x_1x_2 \right]_c = \ad x_1 \left( \left[ x_1^2x_2x_1x_2 \right]_c \right)$, so by direct calculation,
$$ D_2 \left( \left[ x_1^2x_2x_1x_2 \right]_c \right)= 0. $$
Then, we apply $D_2$ to both sides of equality \eqref{palabra} and express the result as a
linear combination of $\left[x_1^3x_2 \right]_c x_1$, $\left[x_1^2x_2
\right]_c x_1^2$ and $\left[x_1x_2 \right]_c x_1^3$, then the
coefficient of $\left[x_1x_2 \right]_c x_1^3$ is
$$\alpha (1-q_{12}q_{21})(1-q_{11}q_{12}q_{21}),$$
so $\alpha=0$.  Then, note also that
$$D_1^2D_2 \left(\left[ x_1^3x_2x_1x_2 \right]_c \right) =0,$$
but
$$D_1^2D_2 \left( \left[ x_1^2x_2 \right]_c^2 \right) = (1-q_{12}q_{21})(1-q_{11}q_{12}q_{21}) (1+q_{11})(q_{2e_1+e_2}+1) \left[ x_1^2x_2 \right]_c. $$
Looking at the proof of Proposition \ref{dimensionG2}, $q_{2e_1+e_2} \neq -1$, so
$\beta=0$.

\medskip

\emph{\textbf{(b)}} Under the conditions \emph{\textbf{(a)}}, \eqref{condG} and \eqref{condG2}, the unique possible Lyndon word of degree $5e_1+2e_j$ is $x_1^3x_2x_1^2x_2$,
and
    $$ \left[ x_1^2x_2x_1x_2x_1x_2 \right]_c = \left[ (\ad x_1)^3x_2, (\ad x_1)^2x_2 \right]_c.$$
Then we proceed as before. One implication is clear. For the other, if
$5e_1+2e_j \notin \de^{+}(\bB)$, then there exists $\alpha \in
\kk$ such that
    $$ \left[ (\ad x_1)^3x_2, (\ad x_1)^2x_2 \right]_c = \alpha (\ad x_1)^2x_2 (\ad x_1)^3x_2.$$
Then, we apply $D_2$ and express the equality as a linear
combination of $(\ad x_1)^3x_2 x_1^2$ and $(\ad x_1)^2x_2 x_1^3$ (we use that $(\ad x_1)^4x_2=0$ by hypothesis); the coefficient of $(\ad
x_1)^2x_2 x_1^3$ is $\alpha \gamma_3$, so $\alpha=0$.
\medskip

\emph{\textbf{(c)}} The proof is similar. Since we consider Lyndon
words without $x_1^3x_2$, $x_1x_2^2$ as factor, the unique
possible Lyndon word of degree $4e_1+3e_j$ is
$x_1^2x_2x_1x_2x_1x_2$, and
    $$ \left[ x_1^2x_2x_1x_2x_1x_2 \right]_c = \left[ \left[ x_1^2x_2x_1x_2 \right]_c  , \left[x_1x_2\right]_c  \right]_c.$$
    If $4e_1+3e_j \notin \de^{+}(\bB)$, then there exist $\alpha_i \in \kk$ such that
\begin{eqnarray*}
    \left[ x_1^2x_2(x_1x_2)^2 \right]_c &=& \alpha_1 \left[ x_1x_2 \right]_c \left[ x_1^2x_2x_1x_2 \right]_c + \alpha_2 \left[ x_1x_2 \right]_c^2 \left[ x_1^2x_2 \right]_c
    \\ && + \alpha_3 x_2\left[ x_1^2x_2 \right]_c^2 + \alpha_4 x_2 \left[ x_1x_2 \right]_c \left[ x_1^3x_2
    \right]_c,
\end{eqnarray*}
since we discard words greater than $x_1^2x_2x_1x_2x_1x_2$
ending in $x_1$ as above; we also discard words with factors $x_1^4x_2$, $x_1x_2^2$, $x_1^3x_2x_1^2x_2$, by the hypothesis about $\bB$. We apply $D_2$ to this equality. Using
the definition of braided commutator, express it as a linear combination of elements of PBW basis, which have degree $4e_1+2e_2$.

The coefficient of $x_2 \left[ x_1x_2 \right]_cx_1^3$ is $\alpha_4 \gamma_3$ since this PBW generator appears only in the expression of $D_2(x_2 \left[ x_1x_2 \right]_c \left[ x_1^3x_2 \right]_c)$. Then, $\alpha_4=0$.

Using this fact, the coefficient of $x_2 \left[ x_1^3 x_2 \right]_cx_1$ is
    $$\alpha_3 \gamma_2 (1+q_{11})q_{11}^2q_{12}q_{21}^2q_{22},$$
since it appears only in the expression of $D_j(x_2\left[ x_1^2x_2 \right]_c^2)$. Then, $\alpha_3=0$.

Now, look at the coefficient of $[x_1x_2]_c^2x_1^2$. It is $\alpha_2 \gamma_2$, so $\alpha_2=0$. Then, we calculate the coefficient of $[x_1^2x_2]_c^2$:
$$ \alpha_1 \gamma_2 \left(\chi(\e_1, 5\e_1+\e_2) - \chi(2\e_1+\e_2, \e_1+\e_2) \right)= \alpha_1 \gamma_2 q_{11}q_{12} \left( q_{11}^3-q_{22}q_{12}q_{21} \right). $$
As $q_{11}^3 \neq q_{22}q_{12}q_{21}$ for each standard braiding, we conclude $\alpha_1=0$.
\medskip

\emph{\textbf{(d)}} If \textbf{\emph{(b), (c)}} holds, then the unique possible Lyndon word which does not have factors $x_1^4x_2$ and $x_1x_2^2$ of degree $5e_1+3e_2$ is $x_1^2x_2x_1^2x_2x_1x_2$, and
  $$ \left[ x_1^2x_2x_1^2x_2x_1x_2 \right]_c = \left[  \left[x_1^{2}x_2\right]_c  , \left[ x_1^2x_2x_1x_2 \right]_c \right]_c. $$
Then this hyperword is not in $S_I$ iff there exist $\nu_i \in \kk$ such that
\begin{align}
\left[ x_1^2x_2x_1^2x_2x_1x_2 \right]_c =& \nu_1 \left[x_1^2x_2x_1x_2\right]_c \left[x_1^2x_2\right]_c + \nu_2 \left[x_1x_2\right]_c\left[x_1^2x_2\right]_c^2 \nonumber
\\ &+ \nu_3 \left[x_1x_2\right]_c^2 \left[x_1^3x_2\right]_c + \nu_4 x_2\left[x_1^2x_2\right]_c \left[x_1^3x_2\right]_c.
\label{combination}
\end{align}

Apply $D_2$ and note that $D_2(\left[ x_1^2x_2x_1^2x_2x_1x_2 \right]_c)=0$ under the hypothesis of $\bB$. Then, express the resulting sum as a linear combination of elements of the PBW basis, which have degree $5e_1+2e_2$.

The hyperword $x_2[x_1^2x_2]x_1^3$ appears only for $D_2  (x_2\left[x_1^2x_2\right]_c \left[x_1^3x_2\right]_c)$, and its coefficient is $\nu_4 \gamma_3$, and as $\gamma_3 \neq 0$ we conclude that $\nu_4=0$.

Analogously, $\left[x_1x_2\right]_c^2 x_1^3$ appears only for $\left[x_1x_2\right]_c^2 \left[x_1^3x_2\right]_c$ (due to $\nu_4=0$). Its coefficient is $\nu_3 \gamma_3$, so $\nu_3=0$.

Note that $D_1^2D_2 (\left[x_1^2x_2x_1x_2\right]_c)=0$. We apply $D_1^2D_2$ to the expression \eqref{combination}, and obtain
$$ 0= \nu_1 \gamma_2 (1+q_{11}) \left[x_1^2x_2x_1x_2\right]_c + \nu_2 \gamma_2 (1+q_{11}) (1+q_ {2\e_1+\e_2}) [x_1x_2]_c [x_1^2x_2]_c. $$
$\left[x_1^2x_2x_1x_2\right]_c$ and $[x_1x_2]_c [x_1^2x_2]_c$
are linearly independent, since they are linearly independent in $\bB(V)$, and we have a surjection $\bB \rightarrow \bB(V)$. Then,
$$ \nu_1 \gamma_2 (1+q_{11})= \nu_2 \gamma_2 (1+q_{11}) (1+q_{2\e_1+\e_2}) =0. $$
But for standard braidings of type $G_2$ we note that $q_{11},
q_{2\e_1+\e_2} \neq -1$ and $\gamma_2 \neq 0$, so $\nu_1= \nu_2=
0$.
\bigskip

The last statement is true since
$$\Delta^{+}(\bB(V))=\left\{ e_1,e_1+e_2,2e_1+e_2,3e_1+e_2, 3e_1+2e_2,e_2 \right\},$$
if the braiding is standard of type $G_2$.
\edem

\begin{obs}\label{obsG2}
Let $V$ be a standard braided vector space of $G_2$ type, and
$\bB$ a braided graded Hopf algebra satisfying the hypothesis of
the Lemma above. In a similar way to Lemma \ref{lemaBn}, if $q_{11}
\notin \G_4, \ q_{22} \neq -1$, then $5e_1+2e_2, 4e_1+2e_2
4e_1+3e_2, 5e_1+3e_2 \notin \de^{+}(\bB)$.

It follows because $x_1^3x_2x_1^2x_2, x_1^2x_2x_1x_2x_1x_2,
x_1^2x_2x_1^{2}x_2x_1x_2 \notin S_I$, using the quantum Serre
relations as in cited Lemma.
\end{obs}

\begin{theorem}\label{presentationG2}
Let $V$ be a standard braided vector space of type $G_2$.

The Nichols algebra ${\mathfrak B}(V)$ is presented by generators $x_1,x_2$, and relations
\begin{eqnarray}
x_{\alpha}^{N_{\alpha}} &=& 0, \quad \alpha \in \de^{+},
\label{alturaPBWG2}
\\  ad_{c}(x_1)^4(x_2) = ad_{c}(x_2)^2(x_1)  &=& 0,
\label{serreG2}
\end{eqnarray}
and if $q_{11} \in \G_4$ or $q_{22}=-1$,
\begin{eqnarray}
\left[ (\ad x_1)^3x_2, (\ad x_1)^2x_2 \right]_c&=&0, \label{G21}
\\  \left[x_1,  \left[x_1^2x_2x_1x_2 \right]_c \right]_c &=&0, \label{G22}
\\  \left[ \left[ x_1^2x_2x_1x_2 \right]_c  , \left[x_1x_2\right]_c  \right]_c
&=&0, \label{G23}
\\ \left[  \left[x_1^{2}x_2\right]_c  , \left[ x_1^2x_2x_1x_2 \right]_c \right]_c &=&
0. \label{G24}
\end{eqnarray}

Moreover, the following elements constitute a basis of ${\mathfrak
B}(V)$:
\begin{eqnarray}
    x_2^{h_{e_2}} \left[x_1 x_2 \right]_c ^{h_{e_1+e_2}} \left[ x_1^2x_2x_1x_2 \right]_c ^{h_{3e_1+2e_2}} \left[x_1^2 x_2 \right]_c ^{h_{2e_1+e_2}} \left[x_1^3 x_2 \right]_c ^{h_{3e_1+e_2}} x_{1}^{h_{e_1}}, \nonumber
    \\ 0 \le h_{\alpha} \le N_{\alpha} - 1. \label{baseG2}
\end{eqnarray}
\end{theorem}

\bdem The statement about the PBW basis follows from Corollary \ref{corollary:genPBW} and the definitions of $x_{\alpha}$'s.
\medskip

Now, let $\bB$ be the algebra presented by generators $x_1,x_2$
and relations \eqref{serreG2}, \eqref{alturaPBWG2}, \eqref{G21}, \eqref{G22}, \eqref{G23} and \eqref{G24}. From Lemma \ref{lemaG2} and Proposition \ref{heigthgenerators}, we have a canonical
epimorphism of algebras $\phi: \bB \rightarrow \bB(V)$.

Consider the subspace $\cI$ of $\bB$ generated by the elements in \eqref{baseG2}. We prove by induction in the sum $S$ of the $h_{\alpha}$'s of a such product $M$ that $x_1M \in \cI$; moreover, we prove that it is a linear combination of products which first hyperletter is least or equal than the first
hyperletter of $M$. If $S=0$, we have $M=1$. Then

$\bullet$ If $M=x_1^{N_1}$, then $x_1M=x_1^{N_1+1}$, that is zero if $N_1=\ord x_1-1$.

\smallskip
$\bullet$ If $M=\left[x_1^3 x_2 \right]_c M'$, then we use that $x_1 \left[x_1^3 x_2 \right]_c= q_{11}^3q_{12}\left[x_1^3 x_2 \right]_cx_1$ to prove that $x_1M \in \cI$, and it is zero or begins with $\left[x_1^3 x_2 \right]_c$.

\smallskip
$\bullet$ If $M=\left[x_1^2 x_2 \right]_c M'$, then we use that
$$x_1 \left[x_1^2 x_2 \right]_c= \left[x_1^3 x_2 \right]_c+q_{11}^2q_{12}\left[x_1^2 x_2 \right]_cx_1.$$
So, we use the inductive step and relation \eqref{G21} to prove that $x_1M \in \cI$, and it is zero or a linear combination of hyperwords that begin with an hyperletter least or equal than $\left[x_1^2 x_2 \right]_c$.

\smallskip
$\bullet$ If $M=\left[x_1^2 x_2x_1x_2 \right]_c M'$, then we deduce from \eqref{G22}
$$x_1 \left[x_1^2 x_2x_1x_2 \right]_c = \chi (\e_1,3\e_1+2\e_2) \left[x_1^2 x_2x_1x_2 \right]_c x_1   ,$$
and using also relations \eqref{G23}, \eqref{G24}, we prove that
$x_1M \in \cI$, and it is zero or a linear combination of
hyperwords that begin with an hyperletter least or equal than
$\left[x_1^2 x_2x_1x_2 \right]_c$.

\smallskip
$\bullet$ If $M=\left[x_1 x_2 \right]_c M'$, then observe that
    $$x_1 \left[x_1 x_2 \right]_c= \left[x_1^2 x_2 \right]_c+q_{11}q_{12}\left[x_1 x_2 \right]_cx_1.$$
In this case, using inductive step, relations \eqref{G22}, \eqref{G23}, and that
\begin{align*}
    \left[x_1^2 x_2 \right]_c\left[x_1 x_2 \right]_c =& \left[ \left[x_1^2 x_2 \right]_c , \left[x_1 x_2 \right]_c \right]_c
    \\ & +\chi(2\e_1+\e_2,\e_1+\e_2)\left[x_1 x_2 \right]_c\left[x_1^2 x_2 \right]_c,
\end{align*}
by definition of braided commutator, we prove that $x_1M \in \cI$, and it is zero or a linear combination of hyperwords that begin with an hyperletter least or equal than $\left[x_1 x_2 \right]_c$.

\smallskip
$\bullet$ If $M=x_2M'$, then we use that $x_1x_2= \left[x_1x_2\right]_c+ q_{12}x_2x_1$, and also $\left[\left[x_1x_2\right]_c, x_2 \right]_c=0$ to prove that $x_1M \in \cI$, and it is zero or a linear combination of hyperwords.
\smallskip

Now, $x_2M$ is a product of non increasing hyperwords or is zero, for each element in \eqref{baseG2}, so $\cI$ is an ideal of $\bB$ containing 1, and then $\cI=\bB$. As the elements in \eqref{baseG2} are a basis of $\bB(V)$, we have that $\phi$ is an isomorphism.
\edem

\subsection{Presentation when the braiding is of Cartan type} \label{subsection:presentation}
\

In this subsection, we present the Nichols algebra of a diagonal
braiding vector space of Cartan type with matrix $(q_{ij})$, by
generators and relations. This was established in \cite[Th.
4.5]{AS3} assuming that $q_{ii}$ has odd order and that order is
not divisible by 3 if $i$ belongs to a component of type $G_2$.
The proof in loc. cit. combines a reduction to symmetric
$(q_{ij})$ by twisting, with results from \cite{AJS} and
\cite{dCP}. We also note that some particular instances were
already proved earlier in this section.

\smallskip

Fix a standard braided vector space $V$ with connected Dynkin diagram and a natural $i \in \unon$. Suppose that $\bB$ is a quotient by an ideal $I \in \bS$ of $T(V)$. \emph{We assume moreover that}
\begin{equation}\label{conditionsalgebra}
    \begin{cases}
        \eqref{qserre}\mbox{ holds in }\bB, \mbox{ if } 1 \leq i \neq j \leq \theta;
        \\ \eqref{relA}\mbox{ holds in }\bB, \mbox{ if } m_{kj}=m_{kl}=1, \ m_{jl}=0;
        \\ \eqref{relB}\mbox{ holds in }\bB, \mbox{ if } m_{kj}=2, \ m_{jk}=1;
        \\ \eqref{relB2}\mbox{ holds in }\bB, \mbox{ if } m_{kj}=2, \ m_{jk}==m_{jl}=1, \ m_{kl}=0;
        \\ V \mbox{ is not of type }G_2.
    \end{cases}
\end{equation}
Note that if \eqref{qserre} holds in an algebra with derivations $D_i$, then \eqref{conditiontransformation} holds also, by Lemma \ref{lemaAHS}. By \ref{transfnichols}, we have an algebra $s_i(\bB)$ provided with skew derivations $D_i$. We call $\tilde{x}_{k}= (ad_c x_i)^{m_{ik}}(x_k) \# 1 \in s_i(\bB)$, for $k \neq i$, and $\tilde{x}_{i}= 1 \# y$. They generate $s_i(\bB)^{1}$ as vector space.
\smallskip

Under these conditions, we prove:

\begin{lema}\label{relationssiB}
The graded algebra $s_i(\bB)$ satisfies \eqref{conditionsalgebra}.
\end{lema}

\bdem \emph{\textsc{Step I:}} we prove that $s_i(\bB)$ satisfies \eqref{qserre}.

\bdem
Each $m \e_k+e_j$, $0 \leq m \leq m_{kj}$ is an element of $ \de( \bB(V_i))$, so
$s_i(m\e_k+\e_j) \in \de(\bB(V))$. As we have a surjective morphism of braided graded Hopf algebras $\bB \rightarrow \bB(V)$, we have $\de(\bB(V)) \subseteq \de(\bB)$.

From Lemma \ref{buscarmij}, $(\ad_c \tilde{x}_{k})^{m} \tilde{x}_{j}=0$ if and only if $x_k^mx_j$ is a linear combination of greater words, for an order in which $x_k<x_j$. Then, by the relation between the Hilbert series of $\bB$ and $s_i(\bB)$ established in Theorem \ref{transfnichols}, \eqref{qserre} for $s_i(\bB)$ is equivalent to
$$s_i \left( (m_{kj}+1)\e_k+\e_j \right) \notin \de^{+} (\bB).$$

When $k=i \neq j$, this says $-\e_i+\e_j \notin \de^{+} (\bB)$ so
\eqref{qserre} holds.

\medskip

To prove \eqref{qserre} for $s_i(\bB)$ when $j=i$, we prove that
$$ (m_{ki}+1)\e_k+((m_{ki}+1)m_{ik}-1)\e_i \notin \de^{+}(\bB).$$
We analyze several cases.
\begin{itemize}
    \item If $m_{ki}=m_{ik}=0$, we have $\e_k-\e_i \notin \de^{+}(\bB)$.

    \item If $m_{ki}=m_{ik}=1$, then $2\e_k+\e_i \notin \de^{+}(\bB)$, because $(\ad x_k)^2x_i=0$.

    \item If $m_{ki}=1, m_{ik}=2$, then $2\e_k+3\e_i \notin
    \de^{+}(\bB)$, since we apply Lemma \ref{lemaBn} to $\bB$ and it satisfies relation
    \eqref{relB} by hypothesis.

    \item If $m_{ki}=2, m_{ik}=1$, then $3\e_k+2\e_i \notin \de^{+}(\bB)$, as before.

\end{itemize}
Then \eqref{qserre} holds, for each $k \neq i$.

\medskip

Now, consider $\theta \geq 3$, and $k,j \neq i$.

If $m_{ik}=m_{ij}=0 $, then $s_i(m\e_k+\e_j)=m\e_k+\e_j$, and
$(m_{kj}+1)\e_k+\e_j \notin \Delta^{+}(\bB)$, since the quantum Serre relation holds in $\bB$.
\smallskip

If $m_{ik}=1 , m_{ij}=0 $, then $s_i(m\e_k+\e_j)=
m\e_i+m\e_k+\e_j$. If we consider $x_j < x_i<x_k$ and look at the
possible Lyndon words in $S_I$, from \eqref{qserre}, it has no
factors $x_i^{2}x_k, x_jx_i$, so the unique possibility is
$x_j(x_kx_i)^{m}$.

$\bullet$ If $m_{kj}=0$, then $x_jx_kx_i=q_{jk}x_kx_jx_i$, so
$x_jx_kx_i \notin S_I$.

$\bullet$ If $m_{kj}=1$, then $x_jx_kx_lx_k \notin S_I$ when
$m_{ki}=1$, since \eqref{relA} is valid in $\bB$, or when
$m_{ki}=2$ we have $q_{kk} \neq -1$ and
\begin{align*}
    x_j(x_kx_i)^2 =& (1+q_{kk})^{-1}q_{ki}^{-1} \ x_jx_k^2x_i^2+ (1+q_{kk})^{-1}q_{ki}q_{kk}^2 \ x_jx_ix_k^2x_i
    \\ =& q_{ki}^{-1}q_{kj}^{-1}q_{kk}^{-2} \ x_kx_jx_kx_i^2 +
    (1+q_{kk})^{-1}q_{ki}^{-1}q_{kj}^{-2}q_{kk}^{-2} \ x_k^2x_jx_i^2
    \\ & +(1+q_{kk})^{-1}q_{ki}q_{kk}^2q_{ji} \ x_ix_jx_k^2x_i,
\end{align*}
so in both cases, $x_j(x_kx_i)^2 \notin S_I$.

$\bullet$ If $m_{kj}=2$, then $m_{ki}=m_{jk}=1$ and $q_{kk} \neq
-1$. The proof is similar to the previous case.
\smallskip

If $m_{ik}=2 , m_{ij}=0 $, then $s_i(m\e_k+\e_j)
=2m\e_i+m\e_k+\e_j$ and $m_{kj}=0,1$. When $m_{kj}=0$, the proof
is clear as above. When $m_{kj}=1$, for $j<k<i$ and considering
only the quantum Serre relations, the unique possible Lyndon word
is $x_j(x_kx_i^2)^2$. But from $[[x_i^2x_k]_c, [x_ix_k]_c]_c=0$,
we deduce that such word is not in $S_I$.
\smallskip

If $m_{ik}=0 , m_{ij}=1 $, then $s_i(m\e_k+\e_j)
=\e_i+m\e_k+\e_j$. If $k<i<j$, note that from $x_kx_i,
x_k^{m_{kj}+1}x_j \notin S_I$, there are not Lyndon words of
degree $\e_i+(m_{kj}+1)\e_k+\e_j$ in $S_I$.
\smallskip

If $m_{ik}=0 , m_{ij}=2 $, then $s_i(m\e_k+\e_j)
=2\e_i+m\e_k+\e_j$, and the proof is analogous to the previous
case.
\smallskip

If $m_{ik}=m_{ij}=1 $, then $m_{kj}=0$, and $s_i(\e_k+\e_j)
=2\e_i+\e_k+\e_j$, which is not in $\de^{+}(\bB)$ from Lemma
\ref{lemaAn}.
\smallskip

If $m_{ik}=2 , m_{ij}=1 $ (it is analogous to $m_{ik}=1, m_{ij}=2
$), then $m_{kj}=0$ and $s_i(\e_k+\e_j) =3\e_i+\e_k+\e_j$. In this
way, $q_{ii} \neq -1$ and if $x_k<x_i<x_j$, then the unique Lyndon
word without $x_i^2x_j$, $x_kx_i^3$ as factors is
\begin{align*}
    x_kx_i^2x_jx_i=& (1+q_{ii})^{-1}q_{ij}^{-1}x_kx_i^3x_j+
    (1+q_{ii})^{-1} q_{ii}^2q_{ij}x_kx_ix_jx_i^2
    \\ \in & \kk(x_ix_kx_i^2x_j) + \kk (x_i^2x_kx_ix_j) + \kk (x_i^3x_kx_j) + \kk
    (x_kx_ix_jx_i^2),
\end{align*}
using the quantum Serre relations, and then there are not Lyndon
words of degree $3\e_i+\e_k+\e_j$ in $S_I$.
\smallskip

So, \eqref{qserre} holds, for each $k,j \neq i, k \neq j$. \edem
\bigskip

\emph{\textsc{Step II:}} $s_i(\bB)$ satisfies \eqref{relA}.
\bdem Consider $m_{kj}=m_{kl}=1$. We prove case-by-case that
$$s_i(2\e_k+\e_j+\e_l) \notin \de^{+} (\bB).$$

$\bullet$ $m_{ij}=m_{ik}=m_{il}=0$, then $s_i(2\e_k+\e_j+\e_l)
=2\e_k+\e_j+\e_l$, so it follows from Lemma \ref{lemaAn}, because
$2\e_k+\e_j+\e_l \notin \Delta^{+}(\bB)$.

\smallskip
$\bullet$  $m_{ij}\neq 0$ (analogously, $m_{il} \neq 0$), so
$m_{ik}=m_{il}=0$, because there are no cycles in the Dynkin
diagram. Then $s_i(2\e_k+\e_j+\e_l) =2\e_k+\e_j+\e_l+m_{ij}\e_i$,
and if we consider $x_k<x_l<x_j<x_i$, using that
$x_kx_i=q_{ki}x_ix_k, x_jx_l=q_{jl}x_lx_j$ and
$x_lx_i=q_{li}x_ix_l$, and also that $x_k^{2}x_l, x_k^{2}x_j
\notin S_I$, we conclude that all possible Lyndon words of degree
$2e_k+e_j+e_l+m_{ij}e_i$ are not elements of $S_I$, except
$x_kx_lx_kx_jx_i^{m_{ij}}$, but also it is not an element of
$S_I$, because $x_kx_lx_kx_j \notin S_I$. Then,
$2\e_k+\e_j+\e_l+m_{ij}\e_i \notin \Delta^{+}(\bB)$.

\smallskip
$\bullet$ $m_{ik}=1$, and therefore $m_{ij}=m_{il}=0$, then,
$s_i(2\e_k+\e_j+\e_l)=2\e_k+\e_j+\e_l+2m_{ik}\e_i$, and if we
consider $x_l<x_i<x_k<x_j$, using that $x_jx_i=q_{ji}x_ix_j,
x_jx_l=q_{jl}x_lx_j$ and $x_lx_i=q_{li}x_ix_l$, and also that
$x_k^{2}x_l, x_k^{2}x_j \notin S_I$, we discard as before all
possible Lyndon words of degree $2\e_k+\e_j+\e_l+2m_{ik}\e_i$,
except $x_lx_kx_jx_kx_i^{2m_{ij}}$, but it is not an element of
$S_I$, because $x_kx_lx_kx_j \notin S_I$. Then
$2\e_k+\e_j+\e_l+2m_{ij}\e_i \notin \Delta^{+}(\bB)$.

\smallskip
$\bullet$ $i=j$ (or analogously, $i=l$):
$s_j(2\e_k+\e_j+\e_l)=2\e_k+\e_j+\e_l \notin \Delta^{+}(\bB)$ if
$m_{jk}=1$ by Lemma \ref{lemaAn}, or $s_j(2\e_k+\e_j+\e_l)
=2\e_k+3\e_j+\e_l \notin \Delta^{+}(\bB)$ if $m_{jk}=2$, by Lemma
\ref{lemaBn}.

\smallskip
$\bullet$ $i=k$: $s_k(2\e_k+\e_j+\e_l)=\e_j+\e_l \notin
\Delta^{+}(\bB)$, since $m_{jl}=0$.

\smallskip

Also, if $\ub \in \{ \e_k+\e_j, \e_k+\e_l,\e_k,\e_j,\e_l\}$, then
$\ub \in \de(\bB(V_i))$, so $s_i(\ub) \in \de(\bB(V))$. The
canonical surjective algebra morphisms from $T(V)$ to $\bB$ and
$\bB(V)$ induce a surjective algebra morphism $\bB \rightarrow
\bB(V)$, so $\de(\bB(V)) \subseteq  \de(\bB)$; in particular, each
$s_i (\ub) \in \de(\bB)$.

Consider a basis as in Proposition \ref{firstPBWbasis} for an
order such that $x_j<x_k<x_l$. From Lemma \ref{lemaAHS}, $x_jx_k$,
$x_kx_l$, $x_jx_kx_l$ are elements of this basis, since they are
not linear combination of greater words modulo $I_i$, the ideal of
$T(V_i)$ such that $s_i(\bB)=T(V_i)/I_i$. In the same way,
$(x_kx_l)(x_jx_k)$, $x_lx_k(x_jx_k)$, $(x_kx_l)x_kx_j$,
$x_k(x_jx_kx_l)$, $x_lx_k^2x_j$ (if $x_k^2 \neq 0$) are elements
of such basis, where the parenthesis indicates the Lyndon
decomposition as non increasing products of Lyndon words. Also,
$x_jx_l$, $x_jx_k^2$, $x_k^2x_l$ are not in such basis by
\eqref{qserre}. By the relations between the Hilbert series in
Theorem \ref{transfnichols} and the fact that $2\e_k+\e_j+\e_l
\notin s_i\left(\de^{+} (\bB)\right)$, we note that $x_jx_kx_lx_k$
is not an element of such basis. Then, this word is a linear
combination of greater words. By Lemma \ref{lemaAn}, this implies
that \eqref{relA} holds in $s_i(\bB)$.
\edem
\bigskip

\emph{\textsc{Step III:}} $s_i(\bB)$ satisfies \eqref{relB}.
\bdem As before, we prove first that $s_i(3\e_k+2\e_j) \notin \de^{+}
(\bB)$ case by case:

$\bullet$ $m_{ik}=m_{ij}=0$, then $s_i(3\e_k+2\e_j)=3\e_k+2\e_j
\notin \Delta^{+}(\bB)$ by hypothesis.

\smallskip
$\bullet$ $m_{ik}=0, m_{ij}=1$, then
$s_i(3\e_k+2\e_j)=2\e_i+3\e_k+2\e_j$. If consider the order in the
letters $x_k < x_i < x_j$, a Lyndon word of degree
$2e_i+3e_k+2e_j$ in $S_I$ begins with $x_k$, and $x_kx_i$ is not a
factor, because $x_kx_i=q_{ki}x_ix_j$. Then the possible Lyndon
words with these conditions are $x_k^{2}x_jx_ix_kx_jx_i$ and
$x_k^{2}x_jx_kx_jx_i^{2}$; the first is not in $S_I$ because from
\eqref{relA} for $j,k,i$ we can express $x_jx_ix_kx_j$ as a linear
combination of greater words, and the second is not in $S_I$
because $x_k^{2}x_jx_kx_j \notin S_I$.

\smallskip
$\bullet$ $m_{ik}=1, m_{ij}=0$, then
$s_i(3\e_k+2\e_j)=3\e_i+3\e_k+2\e_j$. If consider the order in the
letters $x_j < x_i < x_k$, a Lyndon word of degree
$3e_i+3e_k+2e_j$ in $S_I$ begins with   $x_j$, and $x_jx_i$ is not
a factor. Using that also $x_i^{2}x_k, x_j^{2}x_k \notin S_I$, the
possible Lyndon word with these conditions is
$x_jx_kx_ix_jx_kx_ix_kx_i$. But from the condition in the
$m_{rs}$'s, we are in cases $C_{\theta}$ or $F_4$, and we use that
$(\ad x_i)^{2}x_k =0$, $q_{ii} \neq -1$ to replace $x_ix_kx_i$ by
a linear combination of $x_i^{2}x_k$ and $x_kx_i^{2}$, and also
use $x_jx_i=q_{ji}x_ix_j$, so we conclude that
$x_jx_kx_ix_jx_kx_ix_kx_i \notin S_I$.

\smallskip
$\bullet$ $i=j$: $s_j(3\e_k+2\e_j)=3\e_k+\e_j \notin
\Delta^{+}(\bB)$, since $m_{kj}=2$.

\smallskip
$\bullet$ $i=k$: $s_k(3\e_k+2\e_j)=\e_k+2\e_j \notin
\Delta^{+}(\bB)$, since $m_{jk}=1$.
\smallskip

\smallskip

If $\vb \in \{ \e_k+\e_j, 2\e_k+\e_j,\e_k,\e_j \}$, then $\vb \in
\de(\bB(V_i))$, so $s_i(\vb) \in \de(\bB(V))$. As $\de(\bB(V))
\subseteq  \de(\bB)$, in particular we have that each $\vb \in s_i
\left( \de(\bB) \right)$.

As in \textit{a)}, consider a basis as in Proposition
\ref{firstPBWbasis} for an order such that $x_k<x_j$. In a similar
way, $x_kx_j$, $x_k^2x_j$ are elements of this basis, but
$x_k^3x_j$ and $x_kx_j^2$ are not in such basis by \eqref{qserre}.
From Lemma \ref{lemaAHS}, $(x_kx_j)(x_k^2x_j)$,
$x_j(x_k^2x_j)x_k$, $(x_kx_j)^2x_k$, $x_j(x_kx_j)x_k^2$,
$x_j^2x_k^3$ (the last if $x_j^2, x_k^3 \neq 0$) are not linear
combination of greater words modulo $I_i$, so they are elements of
previous basis. And by the relations between the Hilbert series
and the fact that $3\e_k+2\e_j \notin s_i\left(\de^{+}
(\bB)\right)$, we note that the Lyndon word $x_k^2x_jx_kx_j$ is
not an element of such basis. Then, this word is a linear
combination of greater words, and by Lemma \ref{lemaBn}, this
implies that \eqref{relB} holds in $s_i(\bB)$.
\edem
\bigskip

\emph{\textsc{Step IV:}} $s_i(\bB)$ satisfies \eqref{relB2}.
\bdem We prove case-by-case that
$$ s_i(3\e_k+2\e_j+\e_l) \notin \de^{+}(\bB).$$

$\bullet$ $m_{ik}=m_{ij}=m_{il}=0$:
$s_i(3\e_k+2\e_j+\e_l)=3\e_k+2\e_j+\e_l$, and it is not in
$\de^{+}(\bB)$ by Lemma \ref{lemaBn2}.
\smallskip

$\bullet$ $i \neq j,k,l$ and $m_{ik} \neq 0$: the unique
possibility is $m_{ik}= m_{ki}=1$, so $V$ is of type $F_4$. Then,
$s_i(3\e_k+2\e_j+\e_l)=3\e_i+3\e_k+2\e_j+\e_l$. For the order
$x_l<x_j<x_k<x_i$, the unique possible Lyndon word without factors
$x_l x_j^2$, $x_lx_k$, $x_lx_i$, $x_j^2x_k$, $x_jx_i$, $x_kx_i^2$,
$x_k^2x_i$ is $x_lx_jx_kx_ix_jx_k x_ix_kx_i$. Using the quantum
Serre relations, and the fact that $q_{ii}=q_{kk} \neq -1$, we
obtain that this Lyndon word is not in $S_I$. Then,
$3\e_i+3\e_k+2\e_j+\e_l \notin \de^{+}(\bB)$.
\smallskip

$\bullet$ $i \neq j,k,l$ and $m_{ij} \neq 0$: there are not
standard braided vector spaces with these $m_{st}$'s.
\smallskip

$\bullet$ $i \neq j,k,l$ and $m_{il} \neq 0$: the unique
possibility is $m_{il}=m_{li}=1$. Then,
$s_i(3\e_k+2\e_j+\e_l)=3\e_k+2\e_j+\e_l+\e_i$ If we consider
$x_k<x_j<x_l<x_i$, then the unique possible Lyndon word of such
degree without factors $x_kx_l$, $x_kx_i$, $x_jx_i$, $x_k^3x_j$,
$x_kx_j^2$ is $x_k^2x_jx_lx_ix_kx_i$. But by hypothesis,
$$ \left[ \left[x_k^2x_jx_l\right]_c , \left[x_kx_j\right]_c \right]_c = \left[ x_i , \left[ x_kx_j \right]_c \right]_c =0,$$
so $\left[ x_k^2x_jx_lx_ix_kx_i \right]_c= \left[ \left[
x_k^2x_jx_lx_i\right]_c , \left[x_kx_j \right]_c \right]_c = 0$,
and $x_k^2x_jx_lx_ix_kx_i \notin S_I$.
\smallskip

$\bullet$ $i=k$: $s_i(3\e_i+2\e_j+\e_l)=\e_i+2\e_j+\e_l \notin
\de^{+}(\bB)$, by Lemma \ref{lemaAn}.
\smallskip

$\bullet$ $i=j$: $s_i(3\e_k+2\e_i+\e_l)=3\e_k+2\e_i+\e_l \notin
\de^{+}(\bB)$, by Lemma \ref{lemaBn2}.
\smallskip

$\bullet$ $i=k$: $s_i(3\e_k+2\e_j+\e_i)=\e_k+2\e_j+\e_i \notin
\de^{+}(\bB)$, as before.
\smallskip

\smallskip

Now, if $\wb \in \{ \e_k,\e_j,\e_l,\e_k+\e_j,
\e_k+\e_j+\e_l,2\e_k+\e_j, 2\e_k+\e_j+\e_l, 2\e_k+2\e_j+\e_l\}$,
then $\wb \in \de(\bB(V_i))$, so $s_i(\wb) \in \de(\bB(V))$, and
then $s_i (\wb) \in \de(\bB)$.

Consider a basis as in Proposition \ref{firstPBWbasis} for an
order such that $x_k<x_j<x_l$. Then, $x_jx_k$, $x_kx_l$ are
elements of this basis. We know $x_kx_l$, $x_k^3x_j$, $x_kx_j^2$,
$x_kx_jx_lx_k$, $x_k^2x_jx_kx_j$ are not elements of such basis,
since in $\bB$ hold \eqref{qserre}, \eqref{relA} and \eqref{relB}.
By Lemma \ref{lemaAHS}, the relations between the Hilbert series
in Theorem \ref{transfnichols} and the fact that $3\e_k+2\e_j+\e_l
\notin s_i\left(\de^{+} (\bB)\right)$, we note that the Lyndon
word $x_k^2x_jx_lx_kx_j$ is not an element of such basis. Thus
this word is a linear combination of greater words. By Lemma
\ref{lemaBn2}, this implies that \eqref{relB2} holds in
$s_i(\bB)$. \edem
As $s_i(\bB)$ is of the same type as $\bB$, we conclude the proof.
\edem
\bigskip

Let $V$ of type different of $G_2$. We define the algebra $\hat{\bB}(V):= T(V)/\bI(V)$, where $\bI(V)$ is the 2-sided ideal of $T(V)$ generated by
\begin{itemize}
    \item $(\ad_c x_k)^{m_{kj}+1}x_j$, $k \neq j$;
    \item $\left[ (\ad_c x_j)(\ad_c x_k)x_l, x_k \right]_c$, $l \neq k \neq j$, $q_{kk}=-1$, $m_{kj}=m_{kl}=1$;
    \item $\left[ (\ad_c x_k)^2x_j, (\ad_c x_k)x_j \right]_c$, $k \neq j$, $q_{kk} \in \G_3$ or $q_{jj}=-1$, $m_{kj}=2, m_{jk}=1$;
    \item $\left[ (\ad_c x_k)^2(\ad_c x_j)x_l, (\ad_c x_k)x_j \right]_c$, $k \neq j \neq l$, $q_{kk} \in \G_3$ or $q_{jj}=-1$, $m_{kj}=2, m_{jk}=m_{jl}=1$.
\end{itemize}
Compare with the definitions in \cite[Section 4]{AS3}. Since $V$ is of Cartan type, $\bI(V)$ is a Hopf ideal, by Lemmata \ref{lemacopro}, \ref{lemacoproA} and \ref{lemacoproB}. As also $\bI(V)$ is $\zt$-homogeneous, we have $\bI(V) \in \bS$.

By Lemmata \ref{lemaAn}, \ref{lemaBn} and \ref{lemaBn2}, the
canonical epimorphism $T(V) \rightarrow \bB(V)$ induces a
epimorphism of braided graded Hopf algebras
\begin{equation}\label{projection}
    \pi_{V}: \hat{\bB}(V) \rightarrow \bB(V).
\end{equation}

Also, $\hat{\bB}(V)$ satisfies for each $i \in \unon$
conditions on Theorem \ref{transfnichols}, so we can transform it.

\begin{lema}\label{algtransf}
With the above notation, $s_i(\hat{\bB}(V)) \cong \hat{\bB}(V_i)$.
\end{lema}
\bdem By Lemma \ref{relationssiB}, the relations defining $\bI(V_i)$ are satisfied in
$s_i(\hat{\bB}(V))$. Then, the canonical projections from $T(V_i)$ onto $\hat{\bB}(V_i)$, $s_i
\left(\hat{\bB}(V) \right)$ induce a surjective algebra map $\hat{\bB}(V_i) \rightarrow s_i \left(\hat{\bB}(V) \right)$. Reciprocally, each relation defining $\bI(V)$ is
satisfied in $s_i(\hat{\bB}(V_i))$, so we have the following
situation:
$$
\xymatrix{ \hat{\bB}(V) \ar@{>>}[rr] \ar@{~>}[drr]
&& s_i \left(\hat{\bB}(V_i) \right) \\
\hat{\bB}(V_i) \ar@{>>}[rr] \ar@{~>}[urr] && s_i
\left(\hat{\bB}(V) \right) . }
$$
From the relation between the
Hilbert series in Theorem \ref{transfnichols}, for each $\ub \in
\N^{\theta}$ we have
$$ \dim s_i (\hat{\bB}(V))^{\ub}=  \sum_{k \in \N: \ \ub-k\e_i \in \N^{\theta}, \ s_i (\ub-k \e_i) \in \N^{\theta}} \dim \hat{\bB}(V)^{s_i(\ub - k\e_i)}, $$
and a analogous relation for $\dim s_i
(\hat{\bB}(V_i))^{\ub}$. But from the previous surjections
we have
$$ \dim s_i (\hat{\bB}(V))^{\ub} \leq \dim \hat{\bB}(V_i)^{\ub}, \quad \dim s_i (\hat{\bB}(V_i))^{\ub} \leq \dim \hat{\bB}(V)^{\ub}, $$
for each $\ub \in \N^{\theta}$. Using that $s_i^2=\id$, each of
above inequalities is in fact an equality, and
$s_i(\hat{\bB}(V))= \hat{\bB}(V_i)$. \edem

\medskip

We are now able to prove one of the main results of this paper.

\begin{theorem}\label{presentation}
Let $V$ be a braided vector space of Cartan type, of dimension $\theta$, and $C=(a_{ij})_{i,j \in \unon}$ the corresponding finite Cartan matrix, where $a_{ij}:=-m_{ij}$.

The Nichols algebra ${\mathfrak B}(V)$ is presented by generators
$x_{i}$, $1\le i \le \theta$, and relations
\begin{align}
\label{alturaPBW} x_{\alpha}^{N_{\alpha}} &= 0, \quad \alpha \in
\de^{+};
\\ ad_{c}(x_{k})^{1-a_{kj}}(x_{j})  &= 0, \quad k\neq j;
\label{serrebis}
\end{align}
if there exist $j \neq k \neq l$ such that $\ m_{kj}=m_{kl}=1$,
$q_{kk}=-1$, then
\begin{equation} \left[ (\ad x_k)x_j, (\ad
x_k)x_l \right]_c = 0;\label{Arelation}
\end{equation}
if there exist $k \neq j$ such that $\ m_{kj}=2, m_{jk}=1$, $q_{kk} \in \G_3$ or $q_{jj}=-1$, then
\begin{equation} \left[(\ad x_k)^2x_j, (\ad x_k)x_j \right]_c = 0; \label{Brelation}
\end{equation}
if there exist $k \neq j \neq l$ such that $\ m_{kj}=2, m_{jk}=
m_{jl}=1$, $q_{kk} \in \G_3$ or $q_{jj}=-1$, then
\begin{equation}
\label{Brelation2} \left[ (\ad x_k)^2 (\ad x_j) x_l, (\ad x_k)x_j
\right]_c = 0;
\end{equation}
if $\theta=2$ and $V$ if of $G_2$ type, $q_{11} \in \G_4$ or
$q_{22}=-1$, then
\begin{align}
\left[ (\ad x_1)^3x_2, (\ad x_1)^2x_2 \right]_c=&0, \label{G21'}
\\  \left[x_1,  \left[x_1^2x_2x_1x_2 \right]_c \right]_c =&0, \label{G22'}
\\  \left[ \left[ x_1^2x_2x_1x_2 \right]_c  , \left[x_1x_2\right]_c  \right]_c
=&0, \label{G23'}
\\ \left[  \left[x_1^{2}x_2\right]_c  , \left[ x_1^2x_2x_1x_2 \right]_c \right]_c =&
0. \label{G24'}
\end{align}

Moreover, the following elements constitute a basis of ${\mathfrak
B}(V)$:
$$x_{\beta_{1}}^{h_{1}} x_{\beta_{2}}^{h_{2}} \dots x_{\beta_{P}}^{h_{P}}, \qquad 0 \le h_{j} \le N_{\beta_j} - 1, \text{ if }\, \beta_j \in S_I, \quad  1\le
j \le P.$$
\end{theorem}

\bdem We may assume that $C$ is connected. If $V$ is of type $G_2$, then the Theorem was proved in Theorem \ref{presentationG2}. So we can assume $m_{kj} \neq 3$, $k \neq j$.

The statement about the PBW basis was proved in Corollary \ref{corollary:genPBW} -- see the
definition of the $x_{\alpha}$'s in Subsection \ref{subsection:generators}.

\medskip

Consider the image of $x_{\alpha}$ in $\hat{\bB}(V)$; they
correspond in $\bB(V)$ with $x_{\alpha}$, and are PBW generators
for a basis constructed as in Theorem \ref{basePBW}, considering
the same order in the letters. As we observed in \eqref{projection}, there exists a
surjective morphism of braided Hopf algebras $\hat{\bB}(V)
\rightarrow \bB(V)$, so
$$\Delta(\bB(V)) \subseteq \Delta(\hat{\bB}(V)).$$

Also, $\hat{\bB}(V)$ verifies the conditions in Theorem
\ref{transfnichols} for each $i \in \{1, \ldots, \theta\}$, so we
can transform it. By Lemma \ref{algtransf}, the new algebra is
$\hat{\bB}(V_i)$. So, we can continue. Then, consider the
sets
\begin{eqnarray*}
\hat{\Delta} &:=& \cup \{ \Delta(s_{i_1} \cdots s_{i_k}
\hat{\bB}): k \in \N, 1 \leq i_1, \ldots, i_k \leq \theta
\},
\\ \hat{\Delta}^{+} &:=& \Delta \cap \N^{\theta};
\end{eqnarray*}
and then $\hat{\Delta}$ is invariant by the $s_i$'s. Also,
$\Delta(\bB(V)) \subseteq \Delta$, and $$\Delta(s_{i_1} \cdots
s_{i_k} \hat{\bB}(V))= s_{i_1} \cdots s_{i_k}
\Delta(\hat{\bB}(V)).$$

Consider $\alpha \in \hat{\Delta}^{+} \setminus \Delta^{+}(\bB(V))$.
Suppose that $\alpha \neq m \alpha_i$, for all $m \in\N$ and $i
\in\unon$, of minimal height among these roots. For each $s_i$, as
$\alpha$ is not a multiple of $\alpha_i$, we have $s_i(\alpha) \in
\Delta^{+} \setminus \de^{+} \left( \bB(V) \right)$, and then
$\gr(s_i(\alpha))-\gr(\alpha) \geq 0$. But
$\alpha=\sum_{i=1}^{\theta} b_i\e_i$, so $\sum_{i=1}^{\theta}
b_ia_{ij} \leq 0$, and as $b_i \geq 0$, we have
$\sum_{i,j=1}^{\theta} b_ia_{ij}b_j \leq 0$. This contradicts the
fact that $(a_{ij})$ is definite positive, and $(b_i) \geq 0,
(b_i) \neq 0$.

Also, $m \e_i \in \de^{+} (\hat{\bB}) \Leftrightarrow m=
N_{\e_i}, 1$, since $x_i^{N_{\e_i}} \neq 0$. Then,
$$\label{igualdaddelta} \Delta (\hat{\bB}(V)) =
\Delta(\bB(V)) \cup \{ N_{\alpha} \alpha: \alpha \in
\Delta(\bB(V)) \}. $$ It follows since by Corollary \ref{corollary:genPBW} each $\alpha \in
\de^{+}(\bB(V))$ is of the form $$ \alpha= s_{i_1} \cdots
s_{i_m}(\e_j), \quad i_1, \ldots, i_m , j \in \unon. $$ Now,
$N_{\e_j} \e_j \in \Delta (\hat{\bB}(V))$, so
$$ N_{\alpha} \alpha = N_{\e_j} \alpha = s_{i_1} \cdots s_{i_m}(N_{\e_j} \e_j ) \in \Delta (\hat{\bB}(V)). $$
Also, each degree $N_{\alpha} \alpha$ has multiplicity one too in $\Delta (\hat{\bB}(V))$.

Now, for degrees $N_{\alpha} \alpha$, suppose that there are some
Lyndon words of degree $N_{\alpha} \alpha$, and consider one of
them of minimal height. This word $u$ has a Shirshov decomposition
$$ u=vw, \qquad \beta:=\deg(v), \ \gamma:=\deg(w) \in \Delta^{+} (\hat{\bB}(V)). $$
From the previous assumption, we have $\beta, \gamma \in \Delta^{+}(\bB(V))$. Write
$$ \alpha = \sum^{\theta}_{k=1}a_k \e_k, \quad \beta=\sum^{\theta}_{k=1}b_k \e_k, \quad \gamma=
\sum^{\theta}_{k=1}c_k\e_k, $$ so $N_{\alpha}a_k= b_k+c_k$, for
each $k \in \unon$. We can consider $a_1, a_{\theta} \neq 0$ (if
not, we look at a smaller subdiagram).

Now, if we consider $V$ of type $F_4$ and $\beta= 2\e_1+3\e_2
+4\e_3+3\e_4$, then $c_1=0$, $a_1=1$, $N_{\alpha}=2$, or
$a_1=c_1=1$, $N_{\alpha}=3$, since $\alpha, \gamma \neq \beta$.
\begin{itemize}
    \item If $N_{\alpha}=3$, then $3a_2=3+c_2$. Then, $c_2=0$, so $c_3=c_4=0$, or
$c_2=3$, and  $c_3=4$, $c_4=2$. But in both cases we have a
contradiction to $\alpha \in \N^4$.
    \item If $N_{\alpha} =2$, $c_1=0$, then $c_2, c_4$ are odd,
    and $c_3$ is even, not zero. The unique possibility is $\gamma = \e_2+2 \e_3+
    \e_4$, so $\alpha= \e_1+2 \e_2+ 3 \e_3+2 \e_4$. But $q_{\alpha}=q \neq
    -1$, so $N_{\alpha} \neq 2$, which is a contradiction.
\end{itemize}

Thus, we can consider $b_1, c_1 \leq 1$ or $b_{\theta}, c_{\theta}
\leq 1$, so $a_1=b_1=c_1=1$ or $a_{\theta}=b_{\theta} =c_{\theta}
=1$, and in both cases, $N_{\alpha}=2$. For each possible $\beta$
with $b_1 \neq 0$ (by assumption of $a_1 \neq 0$, we have $b_1
\neq 0$ or $c_1 \neq 0$, we look for $\gamma$ such that
$\beta+\gamma$ has even coordinates. In types $A$, $D$ and $E$
there are not such pairs of roots. For the other types,
\begin{enumerate}
    \item $B_{\theta}$: $\beta=\vb_{i\theta}$, $\gamma=\ub_{i+1,\theta}$. Then, $\alpha=\ub_{1\theta}$, but $q_{\alpha}=q_{11} \neq
    -1$, which is a contradiction.
    \item $C_{\theta}$: $\beta=\wb_{11}$, $\gamma=\e_{\theta}$. Then, $\alpha=\ub_{1\theta}$, but $q_{\alpha}=q_{\theta \theta} \neq
    -1$, which is a contradiction.
    \item $F_4$: $\beta=\e_1+\e_2+2 \e_3+2 \e_4$, $\gamma=\e_1+\e_2$, or $\beta=\e_1+2\e_2+2 \e_3+2 \e_4$, $\gamma=\e_1$. In both cases, $\alpha=\e_1+\e_2+\e_3+\e_4$, but $q_{\alpha}=q \neq
    -1$, which is a contradiction.
\end{enumerate}

Then, each root $N_{\alpha} \alpha$ corresponds to
$x_{\alpha}^{N_{\alpha}}$, and each $x_{\alpha}$ as before has
infinite height. The elements
$$x_{\beta_{1}}^{h_{1}} x_{\beta_{2}}^{h_{2}} \dots
x_{\beta_{P}}^{h_{P}}, \qquad 0 \le h_{j} <
\infty , \text{ if }\, \beta_j \in S_I, \quad
1\le j \le P,$$ constitute a basis of $\hat{\bB}(V)$ as vector space.

Now, let $\bar{I}(V)$ be the ideal of $T(V)$ generated by
relations \eqref{serrebis}, \eqref{Arelation}, \eqref{Brelation},
\eqref{Brelation2}, and also \eqref{alturaPBW}. Then we have
$\bI(V) \subseteq \bar{I}(V) \subseteq I(V)$, so the
corresponding projections induce a surjective morphisms of
algebras $\phi: \bB \rightarrow \bB(V)$, where $\bB:= T(V)/\bar{I}(V)$.

$$
\xymatrix{ T(V) \ar@{>>}[d] \ar@{>>}[r] & \hat{\bB}(V)
\ar@{>>}[dl] \ar@{>>}[d]
\\ \bB(V)  & \ar@{>>}[l]^{\phi} \bB }
$$
Also, the elements
$$x_{\beta_{1}}^{h_{1}} x_{\beta_{2}}^{h_{2}} \dots
x_{\beta_{P}}^{h_{P}}, \qquad 0 \le h_{j} < N_{\beta_j}, \text{ if
}\, \beta_j \in S_I, \quad  1\le j \le P,$$ generate $\bB$ as
vector space, because they correspond to the image of elements
that generate $\hat{\bB}(V)$, and are not zero (each non
increasing product of hyperwords as before such that $h_j \geq N_{
\beta_j}$ is zero in $\bB$). But also $\phi$ is surjective, and
the corresponding images of these elements constitute a basis of
$\bB(V)$, so $\phi$ is an isomorphism.

\edem

\subsection*{Acknowledgments}
I thank N. Andruskiewitsch for his guidance, his important
suggestions and the careful reading of this work.

\end{document}